%% file: JamesThesis.tex
\documentclass{BristolPhDThesis11} %Class file includes the front matter - based on Matt Grime's template.  Change 11 to 12 for font size

\usepackage{natbib}
\bibpunct{(}{)}{;}{a}{,}{,} %punctuation for the natbib package
\usepackage{amsmath}
\usepackage{amsthm}
\usepackage{amsfonts}
\usepackage{amssymb}

\usepackage{labelfig}
\usepackage{epsf}
\usepackage{graphicx}
\usepackage{multirow}
\usepackage{pstricks,pst-node,pst-plot}

\author{\Large{James Springham}}
\title{Ergodic properties of linked-twist maps}
\date{September 2008}
\TwoSided
\FrontHeadingFontStyle{\Large\textbf}

\SetFronttitleoffset{2cm}
\SetInsideMargin{2.75cm}
\SetOutsideMargin{2.75cm}
\SetTopMargin{3cm}
\SetBottomMargin{4cm}
\SetBindingOffset{1.5cm}
\ResetGeometry

\begin{document}

\input input/defs %shortcut commands

\ThesisTitle{\Fronttitleoffset}{2cm}
\ThesisFrontMatter
\Abstract{\Fronttitleoffset}{ \input input/abstract}
\Dedication{\Fronttitleoffset}{ \input input/dedication}
\Acknowledgements{\Fronttitleoffset}{ \input input/acknowledgements }
\Declaration{\Fronttitleoffset}
\tableofcontents
%\listoftables
\listoffigures
\cleardoublepage %\newpage %

\ThesisMainMatter

%define equation and figure numbering
\renewcommand{\theequation}{\thesection.\arabic{equation}}
\renewcommand{\thefigure}{\thechapter.\arabic{figure}}
\renewcommand{\thetable}{\thechapter.\arabic{table}}

\input intro/introduction

\input chapter1/chapter1

\input chapter2/chapter2
\input chapter3/chapter3

\input chapter4/chapter4

%
\input input/concl

\renewcommand{\theequation}{\thechapter.\arabic{equation}}

\ThesisBib{input/JamesBib} %bib file

\end{document}

%% file: input/defs.tex
%%User-defined macros

\def\beql{\smallskip
\begin{equation}}
\def\eeql{\smallskip
\end{equation}}
\def\beqal{\smallskip
\begin{eqnarray}}
\def\eeqal{\smallskip
\end{eqnarray}}

\def\beq*{\smallskip
\begin{equation*}}
\def\eeq*{\smallskip
\end{equation*}}
\def\beqa*{\beq*
\begin{array}{ll}}
\def\eeqa*{\end{array}
\eeq*}

\def\sn{\textup{sn}}
\def\cn{\textup{cn}}
\def\dn{\textup{dn}}
\def\sgn{\textup{sgn}}

\def\d{\textup{d}}
\def\pd{\partial}
\def\eps{\varepsilon}

\def\and{\quad\text{and}\quad}
\def\andand{\qquad\text{and}\qquad}

\def\leqs{\leqslant}
\def\geqs{\geqslant}

\def\ltm{linked-twist map }
\def\ltms{linked-twist maps }
\def\Ltm{Linked-twist map }
\def\Ltms{Linked-twist maps }

\newtheorem{thm}{Theorem}[section]
\newtheorem{lem}[thm]{Lemma}
\newtheorem{prop}[thm]{Proposition}
\newtheorem{cor}[thm]{Corollary}
\newtheorem*{defn}{Definition}
\newtheorem{cond}[thm]{Condition}
\newtheorem*{conj}{Conjecture}

%% file: input/abstract.tex
We study a class of homeomorphisms of surfaces collectively known as linked-twist maps. We introduce an abstract definition which enables us to give a precise characterisation of a property observed by other authors, namely that such maps fall into one of two classes termed co- and counter-twisting. We single out three specific linked-twist maps, one each on the two-torus, in the plane and on the two-sphere and for each prove a theorem concerning its ergodic properties with respect to the invariant Lebesgue measure.

For the map on the torus we prove that there is an invariant, zero-measure Cantor set on which the dynamics are topologically conjugate to a full shift on the space of symbol sequences. Such features are commonly known as topological horseshoes. For the map in the plane we prove that there is a set of full measure on which the dynamics are measure-theoretically isomorphic to a full shift on the space of symbol sequences. This is commonly known as the Bernoulli property and verifies, under certain conditions, a conjecture of Wojtkowski's. We introduce the map on the sphere and prove that it too has the Bernoulli property.

We conclude with some conjectures, drawn from our experience, concerning how one might extend the results we have for specific linked-twist to the abstract linked-twist maps we have defined.

%% file: input/dedication.tex
\begin{center}
For my parents. Read it well, there will be questions...
\end{center}

%% file: input/acknowledgements.tex
First and foremost it is my great pleasure to thank my supervisor, Prof.\ Stephen Wiggins. His patience and guidance over the course of four years have been appreciated far more than I have ever told him. He has given more of his time than I had any right to expect. And he has helped me to become a better mathematician, for which I will always be grateful.

I am most grateful to Prof.\ Jens Marklof and Dr.\ Mark Holland for their time taken in reviewing this work and for their many helpful suggestions for improvements. I also thank Dr.\ Rob Sturman, Dr.\ Holger Waalkens and Dr.\ Isaac Chenchiah for their time, advice and interest in my work.

I thank EPSRC who have funded me throughout.

My close friends and office mates have immeasurably improved my time in Bristol. It has been a pleasure to share the experience with them and I would like to thank Dan Bailey, Alice Baker, Hung Manh Bui, Laura Dennis, Laura Hutchinson, David Jessop, Jack Kuipers, Socratis Mouratidis, Jaime Norwood, Dave Oziem, Ben Sandground, Henrik Ueberschaer, Ian Williams and Johanna Ziegler.

I thank my family, my mother Sandra, my father Ernie and my brother Matt who have supported me throughout, and last but by no means least, my girlfriend Michelle for all the love, laughter and lasagne.

%% file: intro/introduction.tex
%\renewcommand{\thechapter}{\Alph{chapter}}
%\setcounter{chapter}{0}
%\chapter*{Introduction} 
%\addcontentsline{toc}{chapter}{Introduction} 
%\chaptermark{Introduction}
%\setcounter{equation}{0}
%\setcounter{figure}{0}
%\label{intro}
%\renewcommand{\thechapter}{\arabic{chapter}}

\chapter{Introduction}
\setcounter{equation}{0}
\setcounter{figure}{0}
\label{intro}

The work in this thesis can be categorised as \emph{dynamical systems} or \emph{non-linear dynamics}. This huge field, in broad terms, studies the trajectories of the points which constitute some space, given some rule which governs the evolution of that space as time progresses. It has strong connections to many of the major fields in pure and applied mathematics, to the natural sciences and to engineering. The present work is primarily of a pure-mathematical nature and relies heavily upon the results and techniques of \emph{ergodic theory}.

Ergodic theory studies dynamical systems with an \emph{invariant measure}. We discuss ergodic theory in greater detail in Section~\ref{LitReview.ErgTh}. Ergodic theory is built upon measure theory, itself one of the cornerstones of mathematical analysis. Its influence is felt in two crucial ways: it allows us to describe and to prove certain \emph{limiting behaviour}, which provides us with information about the evolution of our dynamical system; and it allows us to disregard certain points which evolve in a manner that is atypical and inconvenient for us.

Similarly important is \emph{hyperbolicity}, which we discuss in Section~\ref{LitReview.Hyp}. Hyperbolic behaviour in our dynamical systems is of critical importance insofar as all of our techniques for demonstrating ergodic properties rely upon it. In essence (and of course, we give rigorous definitions later) hyperbolicity concerns the behaviour of those points `close to' some reference point whose evolution we are following. Depending upon the direction of the displacement, these nearby points either approach or move away from our reference point as we evolve the system, but crucially they do not stay at a fixed distance. This behaviour can lead to initial conditions being perpetually thrown apart and back together and result in a \emph{mixing} of the ambient space.

In the remainder of this introduction we will introduce the maps that we shall study and state the three main theorems we shall prove. We do not do so by the most direct route however, preferring first to motivate the concept of hyperbolicity in a simple example. This occupies Section~\ref{intro.motiv}. In Section~\ref{intro.abs} we introduce the reader to the \emph{linked-twist maps} with whose properties this work is concerned. We do this first in an abstract setting which enables us highlight what unites them all and classify them in an important way. Finally Section~\ref{Intro.Theorems} is divided into three parts, in each of which we define a linked-twist map and state a theorem we shall prove for that map.

Following on from this, the remainder of our thesis is organised as follows. In Chapter~\ref{chapter1} we provide a literature review which is divided into four sections. In Section~\ref{LitReview.ErgTh} we discuss ergodic theory, providing the definitions we will need throughout this work, in particular of the Bernoulli property. In Section~\ref{LitReview.Hyp} we discuss hyperbolicity and describe some important results we will use. In Section~\ref{LitReview.Ltm} we survey those results already known for the maps we shall study. Lastly in Section~\ref{LitReview.App} we shall discuss a number of applications which can be modelled by linked-twist maps.

Chapters~\ref{chapter2}, \ref{chapter3} and \ref{chapter4} are where we prove the new results. In each case we define the map and state the theorem later in this introduction, then give a detailed breakdown of the method at the start of the chapter. In Chapter~\ref{chapter2} we show that a linked-twist map defined on a subset of $\mathbb{T}^2$ has an invariant, zero-measure Cantor set on which the dynamics are topologically conjugate to a full shift on $N$ symbols. For further details see Section~\ref{Intro.Theorems.Torus}. In Chapter~\ref{chapter3} we show that a linked-twist map defined on a subset of the plane has the Bernoulli property on a set of full Lebesgue measure. This verifies (under certain conditions) a conjecture of Wojtkowski's (\citeyear{woj}), a precise statement of which is postponed until that chapter, where we establish the required notation. We give more details in Section~\ref{Intro.Theorems.Plane}. Finally in Chapter~\ref{chapter4} we prove the Bernoulli property for a linked-twist map defined on a subset of $\mathbb{S}^2$. We introduce this map in Section~\ref{Intro.Theorems.Sphere}.

We conclude in Chapter~\ref{concl} by analysing the results we have established and discussing the strengths and weaknesses of our methods. There are some obvious generalisations which suggest themselves as well as some different directions one could take whilst still building upon the work we have done, so we consider both. Based on what we have learned we feel confident in making some conjectures and we include these here.

\input intro/section1
\input intro/section2

\input intro/section3

%% file: intro/section1.tex
\section{Motivation}\label{intro.motiv}
\setcounter{equation}{0}

We begin by describing a system which illustrates hyperbolicity in perhaps the simplest non-trivial setting. We will use some of the language of ergodic theory and hyperbolic theory to be introduced in Sections~\ref{LitReview.ErgTh} and \ref{LitReview.Hyp}. The reader who is unfamiliar with these terms is encouraged to skip forward to these definitions as necessary, although we have tried to keep the exposition as elementary as is possible.

\subsection{A hyperbolic toral automorphism}\label{intro.motiv.catmap}

\begin{figure}[htp]
\centering
(a)\includegraphics[totalheight=0.12\textheight]{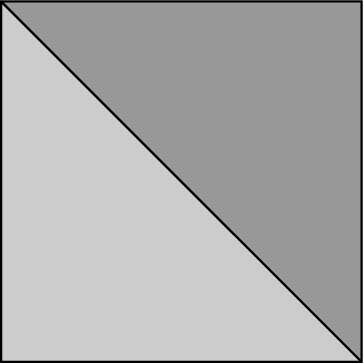}\quad
(b)\includegraphics[totalheight=0.36\textheight]{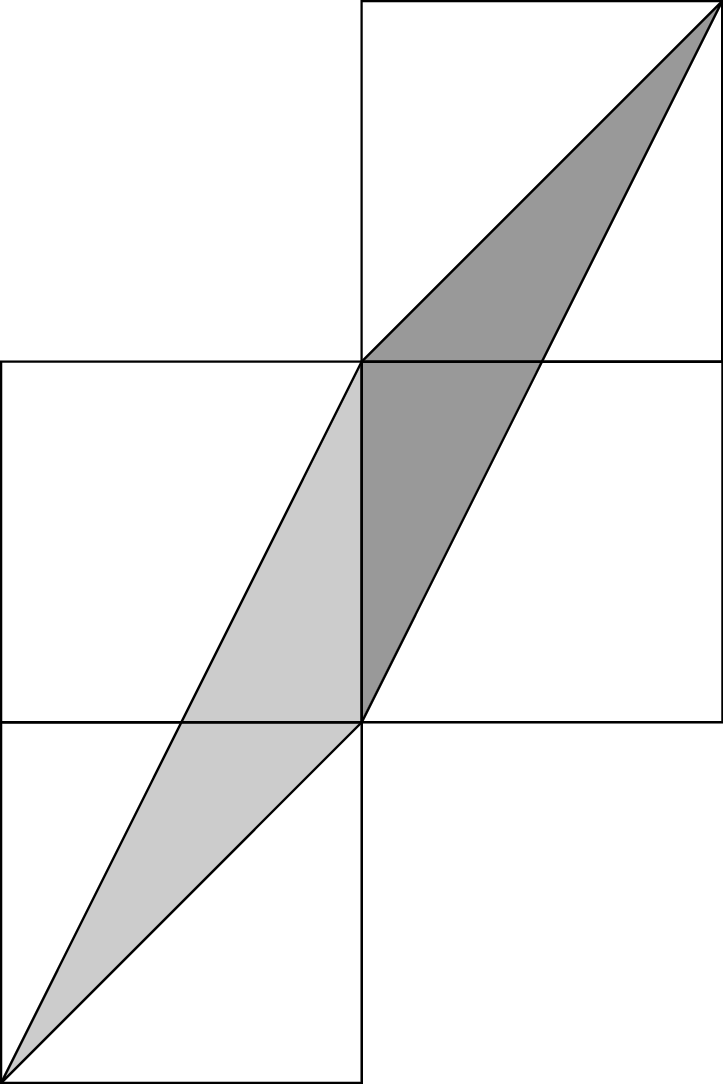}\quad
(c)\includegraphics[totalheight=0.12\textheight]{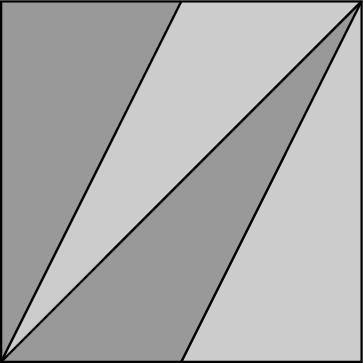}
\caption[Hyperbolic toral automorphism]{The `cat map'. Part~(a) illustrates $\mathbb{T}^2$ which we represent as the unit square in the plane. The shading will help us to illustrate the map. Part~(b) shows the image of the unit square under $H$ if we consider $H$ as a linear map of the plane (i.e.\ without taking the image modulo $\mathbb{Z}^2$). Part~(c) shows how this image looks upon projection to $\mathbb{T}^2$.}
\label{fig:catmap}
\end{figure}

\emph{Hyperbolic toral automorphisms} are canonical examples of dynamical systems displaying hyperbolic behaviour. We describe one here,  commonly known as the \emph{cat map}. More details can be found in most dynamical systems text; we recommend \citet{kh} or \citet{bs}. Given the two-torus $\mathbb{T}^2=\mathbb{R}^2\backslash\mathbb{Z}^2$, the cat map is the linear diffeomorphism $H:\mathbb{T}^2\to\mathbb{T}^2$ given by
\footnote{It is perhaps more common to define the map as $(x,y)\mapsto(2x+y,x+y)$ but this is merely a matter of personal taste and the results we will quote hold for any hyperbolic toral automorphism. When we introduce linked-twist maps on $\mathbb{T}^2$ we will wish to emphasise the cat map as a special case, and for this purpose our definition is more convenient.}
\beq*
H(x,y)=(x+y,x+2y)\mod\mathbb{Z}^2.
\eeq*
We naturally think of $\mathbb{T}^2$ as the unit square in the plane with opposing sides identified. In Figure~\ref{fig:catmap} we illustrate $H$ by first viewing it as a linear map of the plane and then seeing how the `pieces' fit back together on $\mathbb{T}^2$.

Let us describe, without giving the general definition, what we mean when we say that the cat map is hyperbolic. The Jacobian matrix is given by
\beq*
DH_z=\left( \begin{array}{cc}
1 & 1 \\ 1 & 2
\end{array}\right)
\eeq*
and is independent of $z\in\mathbb{T}^2$. It has distinct real eigenvalues $0<\lambda_-<1<\lambda_+=1/\lambda_-$ and corresponding eigenvectors $v_{\pm}=\left(1,\lambda_{\pm}-1\right)$. Using only elementary linear algebra we can draw some simple conclusions about the dynamics of $H$.

Suppose that $z\in\mathbb{T}^2$ and consider the line through $z$ having gradient $v_-$; we call this line the \emph{stable manifold of $z$}. It is easily checked that the gradient is irrational and so the line extends indefinitely and never self-intersects. Let $z'=z+kv_-$, where $k\in\mathbb{R}$, be on this line.
\footnote{For the benefit of a cleaner exposition we are not appending `modulo $\mathbb{Z}^2$' to our points.} 
Then $H(z')=H(z+kv_-)=H(z)+H(kv_-)=H(z)+k\lambda_-v_-$, i.e.\ $H(z')$ is in the unstable manifold of $H(z)$. Moreover the distance between the points (as measured along the unstable manifold) is smaller by a factor of $\lambda_-$ than the corresponding distance between $z$ and $z'$.

We can repeat this construction using $v_+$ in place of $v_-$ to obtain the \emph{unstable manifold of $z$}. In this case the distance between points is increased by a factor of $\lambda_+$. These facts together show that the cat map is \emph{hyperbolic}; in fact we can say more than this. The picture to have in mind is of the two distinct (in fact, orthogonal) directions experiencing stretching and contraction respectively. The constructions we have given hold for any $z\in\mathbb{T}^2$ and the growth rates established hold uniformly at each point, so in fact we say $H$ is \emph{uniformly hyperbolic} or even an \emph{Anosov diffeomorphism}. 

It transpires that from these few facts one can establish a great deal about the dynamics of the cat map. In particular it is \emph{ergodic, mixing} and has the \emph{Bernoulli property}. In Section~\ref{LitReview.Hyp} we will describe a theorem due to \citet{ks} which gives sufficient criteria for a map to have all of these properties. One could certainly use this theorem to establish them for the cat map; however one would, metaphorically speaking, be using a sledgehammer to crack a walnut. For more elegant ways to prove such results we recommend the book of \citet{bs}, in which many more than these three properties are established for hyperbolic toral automorphisms.

%% file: intro/section2.tex
\section{Abstract linked-twist map theory}\label{intro.abs}
\setcounter{equation}{0}

In this section we describe what we will call an \emph{abstract linked-twist map}. The results presented in this thesis are all for \emph{specific} linked-twist maps and the reader who is eager to understand the maps we have studied and the results we have proven can safely overlook this section on first reading. We would encourage her to return to this material later though, for two reasons.

First, this section is our attempt to formalise what precisely it is that the different maps on different surfaces that are all referred to as linked-twist maps have in common. This is perhaps a simple exercise but nevertheless it serves to draw together the results we present.

Second, perhaps more interestingly, we define a property of linked-twist maps which divides them into two classes, namely the \emph{co-twisting} and \emph{counter-twisting} classes. The distinction can have great implications for the dynamics of otherwise similar maps. Other authors have noticed this distinction but have treated it as something which must be \emph{determined} for a given map; conversely we define it for an abstract linked-twist map and later \emph{prove} that a given linked-twist map is either co- or counter-twisting. We are grateful to Prof.\ Robert MacKay for his helpful suggestion, from which this idea was born.

\subsection{Review section: smooth embeddings}\label{intro.abs.rev}

We begin with a number of definitions from the field of differential geometry. The terminology will be necessary in order to define an abstract linked-twist map; the reader who is already comfortable with the definition of a \emph{smooth manifold} and an \emph{orientation-preserving embedding} can safely skip these. Our definitions are taken from the excellent book of \citet{docarmo}. We also recommend the book of \citet{spivak} or the short review section given by \citet{bs}.

\begin{defn}[Smooth manifold of dimension 2]
A smooth manifold is a set $S$ together with a family of one-to-one maps $\phi_{\alpha}:U_{\alpha}\to S$ of open sets $U_{\alpha}\subset\mathbb{R}^2$ into $S$ such that
\begin{enumerate}
\item $\bigcup_{\alpha}\phi_{\alpha}(U_{\alpha})=S$, and
\item for each pair $\alpha,\beta$ with
\beq*
W=\phi_{\alpha}(U_{\alpha})\cap\phi_{\beta}(U_{\beta})\neq\emptyset
\eeq*
we have that
\begin{enumerate}
\item $\phi_{\alpha}^{-1}(W)$ and $\phi_{\beta}^{-1}(W)$ are open sets in $\mathbb{R}^2$, and
\item $\phi_{\beta}^{-1}\circ\phi_{\alpha}$ and $\phi_{\alpha}^{-1}\circ\phi_{\beta}$ are differentiable maps.
\end{enumerate}
\end{enumerate}
\end{defn}

The pair $(U_{\alpha},\phi_{\alpha})$ with $p\in\phi_{\alpha}(U_{\alpha})$ is called a \emph{coordinate system of $S$ around $p$}. The image $\phi_{\alpha}(U_{\alpha})$ is called a \emph{coordinate neighbourhood} and if $q=\phi_{\alpha}(u_{\alpha},v_{\alpha})\in S$, we say that $(u_{\alpha},v_{\alpha})$ are the \emph{coordinates} of $q$ in this coordinate system.

\begin{defn}[Orientable; oriented]
A smooth manifold $S$ is called orientable if it is possible to cover it with a family of coordinate neighbourhoods in such a way that if $p\in S$ belongs to two such neighbourhoods then the change of coordinates has positive Jacobian. The choice of such a family is called an orientation of $S$ and $S$ is called oriented.
\end{defn}

Familiar examples of orientable surfaces include the two-torus $\mathbb{T}^2$ and the two-sphere $\mathbb{S}^2$. Conversely the M\"{o}bius strip is not orientable.

We now extend the notion of a differentiable map in the context of smooth manifolds of dimension 2.

\begin{defn}[Differentiable map]
Let $S_1$ and $S_2$ be smooth manifolds of dimension 2. A map $f:S_1\to S_2$ is differentiable at $p\in S_1$ if given a parametrization $\psi:V\subset\mathbb{R}^2\to S_2$ around $f(p)$ there exists a parametrization $\phi:U\subset\mathbb{R}^2\to S_1$ around $p$ such that $f(\phi(U))\subset\psi(V)$ and the map
\beq*
\psi^{-1}\circ f\circ\phi:U\subset\mathbb{R}^2\to\mathbb{R}^2
\eeq*
is differentiable at $\phi^{-1}(p)$. The map $f$ is differentiable on $S_1$ if it is differentiable at every $p\in S_1$.
\end{defn}

\begin{defn}[Immersion]
A differentiable map $f:S\to\mathbb{R}^3$ of $S$, a smooth manifold of dimension 2, is an immersion if the differential
\beq*
Df_p:T_p(S)\to T_{f(p)}(\mathbb{R}^3)
\eeq*
is injective for each $p\in S$.
\end{defn}

We can now state the definition of an embedding.

\begin{defn}[Embedding]
Let $S$ be a smooth manifold of dimension 2. A differentiable map $f:S\to\mathbb{R}^3$ is an embedding if it is an immersion and a homeomorphism onto its image.
\end{defn}

\begin{figure}[htp]
\centering
(a)\includegraphics[totalheight=0.12\textheight]{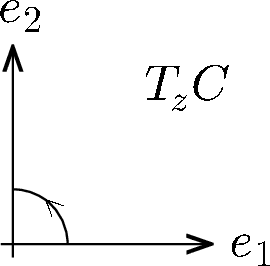}\quad
(b)\includegraphics[totalheight=0.12\textheight]{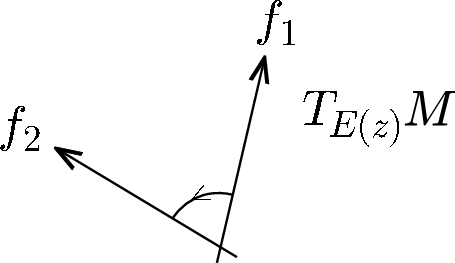}\quad
(c)\includegraphics[totalheight=0.12\textheight]{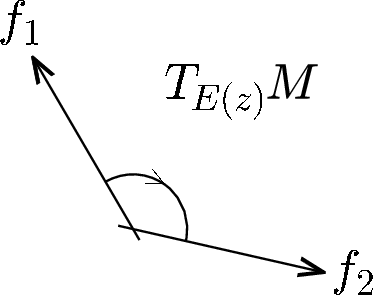}
\caption[Orientation-preserving and orientation-reversing embeddings]{Given a differentiable map $E:C\to M$, the figure shows two possibilities for the image of the standard basis of $\mathbb{R}^2$, denoted $(e_1,e_2)$ and shown in part~(a), under the differential $DE_z$. Parts~(b) and~(c) show bases $(f_1,f_2)$ of $\mathbb{R}^2$, where $f_j=DE_z(e_j)$ for $j=1,2$. In part~(b) $DE_z$ has preserved the orientation, or `handedness', of the standard basis, as shown by the arrow. The corresponding map $E:C\to M$ is called orientation-preserving. Conversely in part~(c) $DE_z$ reverses the orientation of the basis. In this case $E:C\to M$ is called orientation-reversing.}
\label{fig:orientation}
\end{figure}

Finally, an embedding is called \emph{orientation-preserving} if its Jacobian has positive determinant, and \emph{orientation-reversing} otherwise. We illustrate the situation in Figure~\ref{fig:orientation}.

\subsection{Abstract linked-twist maps}\label{Intro.Background.Abstract}

Let $\mathbb{S}^1$ be the circle. Without loss of generality we assume a coordinate $x\in[0,1]$ on $\mathbb{S}^1$, where $0$ and $1$ are identified. In some situations it will be convenient to use some other interval in place of $[0,1]$; in that case obvious amendments should be made to our definitions.

Let $I=[i_0,i_1]\subset\mathbb{R}$ be a closed interval. Moreover we will want $I\subset[0,1]$ (or in the closed interval we use in place of $[0,1]$ as the case may be). The Cartesian product $C=\mathbb{S}^1\times I$ is called a \emph{cylinder} or an \emph{annulus} and consists of pairs $(x,y)$ such that $x\in\mathbb{S}^1$ and $y\in I$. $C$ is an oriented smooth manifold (with boundary) and we identify the tangent space $T_zC$ at a point $z\in C$ with $\mathbb{R}^2$. We give $T_zC$ the standard basis $(e_1,e_2)$, where $e_1=(1,0)$ and $e_2=(0,1)$ in the usual Cartesian coordinates.

We define a class of homeomorphisms of $C=\mathbb{S}^1\times I$:

\begin{defn}[Twist map; twist function]
A twist map $T:C\to C$ is a map of the form
\beq*
T(x,y)=(x+t(y),y),
\eeq*
where $t:I\to\mathbb{S}^1$, called a twist function, satisfies the following conditions:
\begin{enumerate}
\item $t$ is continuous on $[i_0,i_1]$ and differentiable on $(i_0,i_1)$,
\item $t(i_0)=0$ and $t(i_1)=1$ (or an equivalent condition using a different interval for $\mathbb{S}^1$),
\item $\d t/\d y>0$ on $(i_0,i_1)$.
\end{enumerate}
\end{defn}

\begin{figure}[htp]
\centering
(a)\includegraphics[totalheight=0.18\textheight]{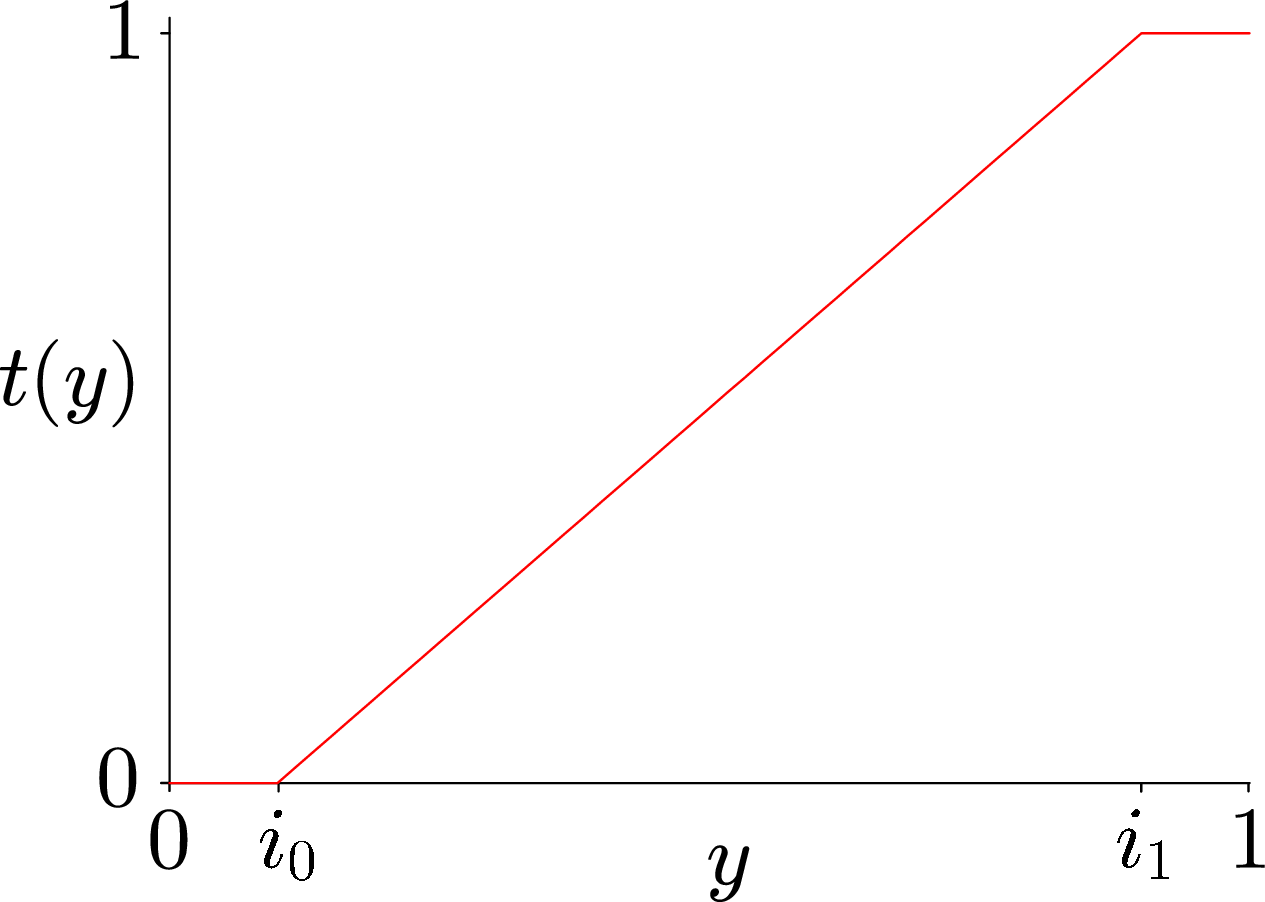}\qquad
(b)\includegraphics[totalheight=0.18\textheight]{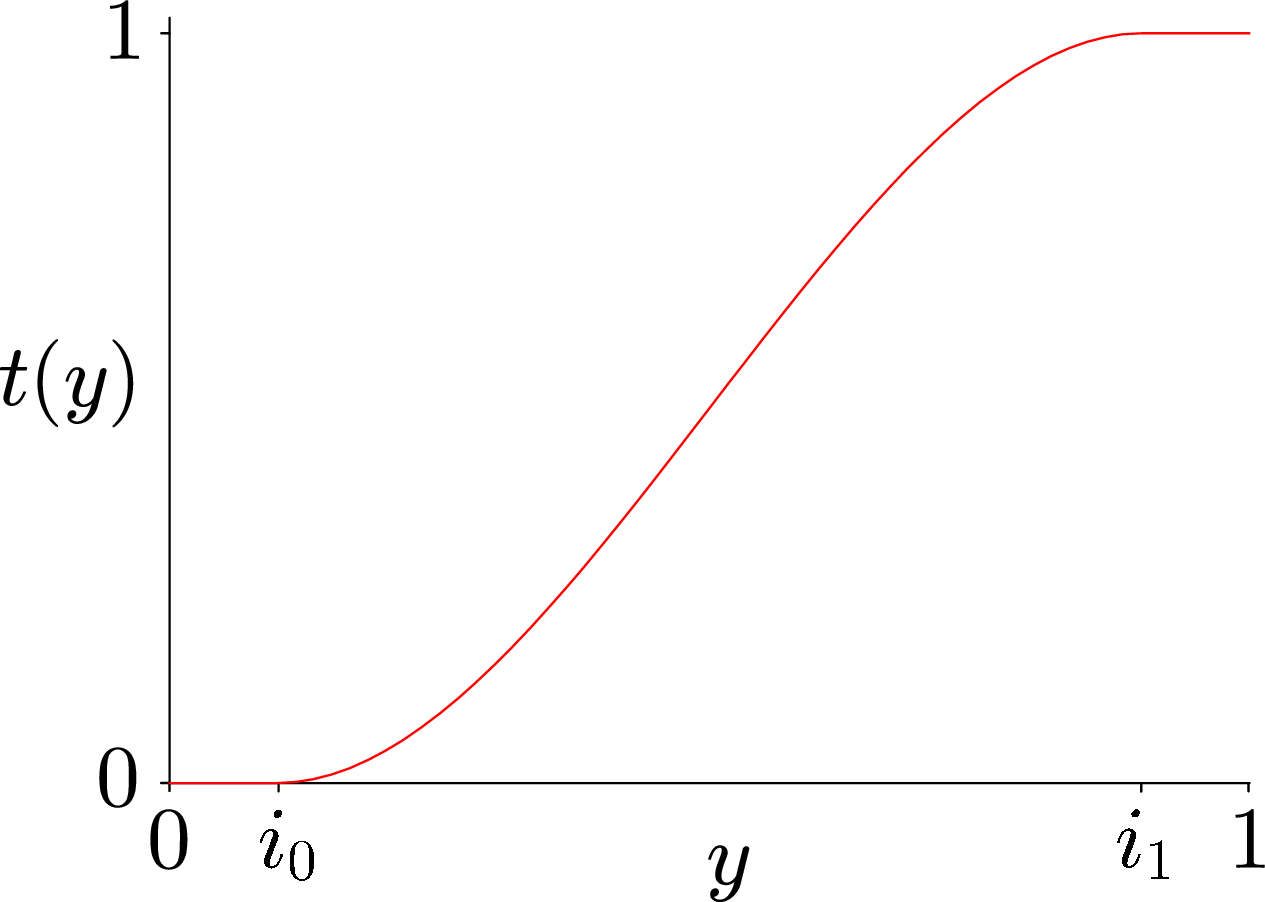}
\caption[Linear and smooth twist functions]{Linear and smooth twist functions respectively (recall that $0$ and $1$ are identified in $\mathbb{S}^1$). Our twist maps (and hence our linked-twist maps) will be constructed from the former. Each introduces a different problem into the analysis; the former because of the non-differentiable points and the latter because there is no lower bound on the derivative.}
\label{fig:twist_function}
\end{figure}

We comment that other authors call $T$ an \emph{integrable} twist map. $T$ preserves area (Lebesgue measure) and orientation; see \citet{kh}. Two possibilities for the twist function $t$ are shown in Figure~\ref{fig:twist_function}. Part~(a) of the figure illustrates a \emph{linear} twist (we should properly call this an affine twist, of course), of the kind our twist maps will be constructed from. It is defined by
\beql
\label{eqn:lin_twist_fn}
t(y)=\left\{
\begin{array}{r@{\quad}l}
(y-i_0)/(i_1-i_0) &\text{if }y\in [i_0,i_1], \\
0 & \text{otherwise.} \end{array} \right.
\eeql
The function is not differentiable at $y\in\{i_0,i_1\}$.

Part~(b) shows a \emph{smooth} (i.e.\ everywhere differentiable) twist of the kind studied by \citet{be}. It is defined by a cubic equation in $y$. Smooth twists require a different kind of analysis to that which we shall conduct and we do not intend to discuss them in this thesis; see the original paper or \citet{sturman} for further details.

We  now define an abstract linked-twist map on a subset $R$ of a two-dimensional smooth manifold $M$. We do so with reference to the cylinders we will embed in $M$ to create $R$. Later on, when we define the linked-twist maps to be studied in this thesis, we more commonly do so directly on $R\subset M$. We introduce an important definition.

\begin{defn}[Transversal embedded cylinders]
Consider two embedded cylinders in some two-dimensional manifold $M$, i.e.\ we have cylinders $C_i$ and diffeomorphisms $E_i:C_i\to M$ for $i=1,2$. Suppose that $E_1(C_1)\cap E_2(C_2)\neq\emptyset$ and let $z_i\in C_i$ be such that $E_1(z_1)=E_2(z_2)\in M$. We will say that such embedded cylinders are transversal if and only if the vectors $(DE_1)_{z_1}(e_1)$ and $(DE_2)_{z_2}(e_1)$, which lie in $T_{E_1(z_1)}M=T_{E_2(z_2)}M$, are themselves transversal in the usual sense (i.e.\ they form a basis for the tangent space).
\end{defn}

We call the connected region(s) $E_1(C_1)\cap E_2(C_2)$ the \emph{intersection region(s)}. For examples of pairs of transversal embedded cylinders the reader is encouraged to look ahead to Figures~\ref{fig:R_torus}, \ref{fig:A_plane} and~\ref{fig:R_sphere}. We can now define a linked-twist map.

\begin{defn}[Linked-twist map]
Let $M$ be a two-dimensional oriented smooth manifold and let $E_i:C_i\to M$, for $i=1,2$,  be a pair of transversal embeddings of cylinders $C_i=\mathbb{S}^1\times I_i$ into $M$. Denote $R=E_1(C)\cup E_2(C)\subset M$. Let $T_i:C_i\to C_i$ for $i=1,2$ be two twist maps given by $T_i(x,y)=(x+t_i(y),y)$ where the twist functions $t_i: I\to\mathbb{S}^1$ satisfy the conditions in the definition above.

For $i=1,2$ and $p\in M$ define $H_i:R\to R$ by
\beql
\label{eqn:abstract_ltm}
H_i(p)=\left\{
\begin{array}{r@{\quad}l}
E_i\circ T_i\circ E_i^{-1}(p) &\text{if }p\in E_i(C), \\
id & \text{otherwise,} \end{array} \right.
\eeql
where $id$ denotes the identity map. A linked-twist map $H:R\to R$ is given by the composition $H=H_2^k\circ H_1^j$ where $j$ and $k$ are positive integers.
\end{defn}

All linked-twist maps of this form can be categorised as either co- or counter-twisting. The definition is as follows:

\begin{defn}[Co-twisting; counter-twisting]
Let $H$ be a linked-twist map as above and let $E_1,E_2:C\to M$ be the transversal embeddings with which it is defined. If both $E_1$ and $E_2$ are orientation-preserving, or both $E_1$ and $E_2$ are orientation-reversing, then we say that $H$ is counter-twisting. Conversely if one of $E_1,E_2$ is orientation-preserving and the other orientation-reversing, then we say that $H$ is co-twisting.
\end{defn}

We have some comments to make regarding the definition.

First and foremost it might seem to the reader counter-intuitive to give the definition as we have, with the co-twisting systems defined as those where the embedded cylinders have different orientations; in fact, in light of our definition, we agree. However there is a considerable literature for linked-twist maps and we would like our definition to agree with it. The terminology seems to have been introduced by \citet{sturman} and their reasoning can be best understood once we have defined linked-twist maps on the torus; we do this in the next section.

Second, the counter-twisting maps (at least, all of the explicit examples of which we are aware) are more difficult to analyse than the corresponding co-twisting maps. In this thesis we deal exclusively with co-twisting linked-twist maps so we do not intend to say too much about why this is so, but when we survey the literature in Section~\ref{LitReview.Ltm} we will see that, where corresponding co- and counter-twisting maps can be shown to have strong ergodic properties, the criteria are more restrictive in the latter case.

Third, we will dispense entirely with the other notion introduced by \citet{sturman} of co- and counter-\emph{rotating} linked-twist maps. This notation was intended to explain the relative sense of rotation of the two twist maps acting on the embedded cylinders, but leads to the somewhat uncomfortable situation whereby planar linked-twist maps (introduced in the next section) are simultaneously co-twisting and counter-rotating or \emph{vice versa}.

%% file: intro/section3.tex
\section{Definitions and statements of theorems}\label{Intro.Theorems}
\setcounter{equation}{0}

We now introduce the linked-twist maps to be studied and state the theorems we shall prove. All of these maps fit the abstract definition we have given above, though we shall not prove so in every case. In the first two cases it is quite obvious. In the case of the third map the embedding uses functions with which the reader may not be familiar so we will provide all the details. We shall prove that each map is co-twisting.

\subsection{Linked-twist maps on the two-torus}\label{Intro.Theorems.Torus}

The simplest linked-twist maps to define and analyse are those on the torus. In this section we will define a toral linked-twist map and state a theorem to be proven in Chapter~\ref{chapter2}. We give an overview of results in the literature for toral linked-twist maps in Section~\ref{LitReview.Ltm.Torus}. In Section~\ref{LitReview.App.Torus} we discuss a situation where toral linked-twist maps can be used to model the behaviour of certain physical phenomena. As mentioned above we will define the map directly on the torus.

Let $\mathbb{S}^1$ denote the closed unit interval $[0,1]$ with opposite ends identified. We identify the two torus, denoted $\mathbb{T}^2$ with the Cartesian product $\mathbb{S}^1\times\mathbb{S}^1$. This gives us two angular coordinates $(x,y)$.

\begin{figure}[htp]
\centering
\includegraphics[totalheight=0.3\textheight]{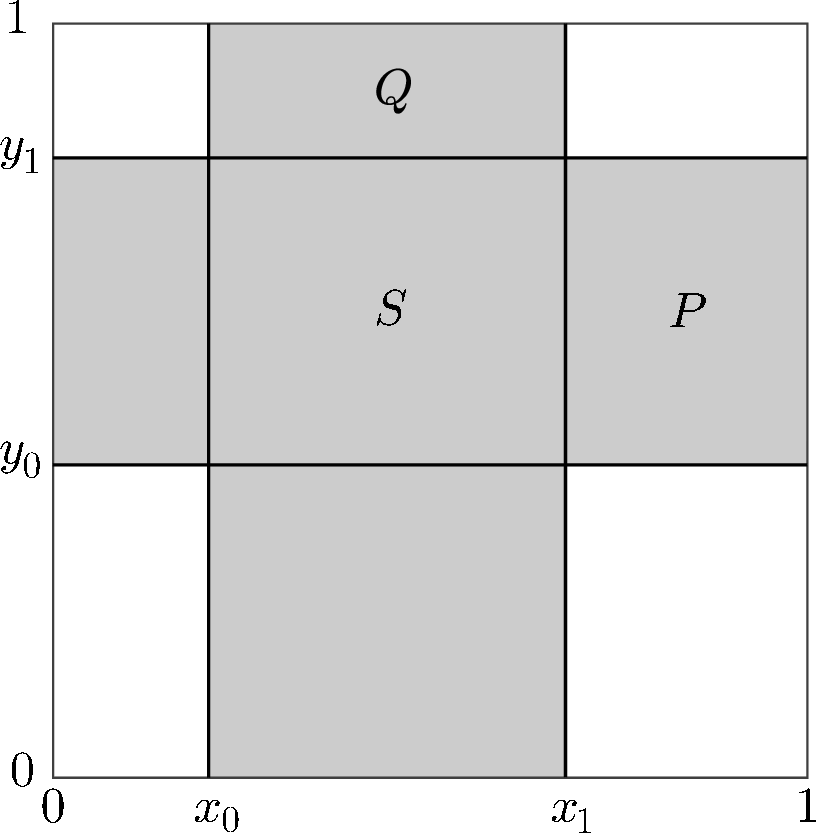}
\caption[The manifold $R\subset\mathbb{T}^2$]{The manifold $R\subset\mathbb{T}^2$ (shaded).}
\label{fig:R_torus}
\end{figure}

Fix four constants $0<x_0<x_1<1$ and $0<y_0<y_1<1$. We define two embedded cylinders $P,Q\subset\mathbb{T}^2$ as follows:
\beq*
P=\{(x,y):x\in\mathbb{S}^1,y_0\leqslant y\leqslant y_1\}\and Q=\{(x,y):x_0\leqslant x\leqslant x_1,y\in\mathbb{S}^1\}.
\eeq*
We shall call $P$ a `horizontal' annulus and $Q$ a `vertical' annulus. We denote by $R=P\cup Q$ the manifold on which our linked-twist map will be defined and by $S=P\cap Q$ the `intersection region'. See Figure~\ref{fig:R_torus}.

The set $\pd P_0=\{(x,y):x\in\mathbb{S}^1,y=y_0\}$ denotes the `lower' boundary of $P$, with the `upper' boundary $\pd P_1$ defined similarly. The `left-hand' boundary of $Q$ is denoted $\pd Q_0$ and the `right-hand' boundary denoted $\pd Q_1$. Again, these are defined similarly. Finally we denote $\pd P=\pd P_0\cup\pd P_1$ and $\pd Q=\pd Q_0\cup\pd Q_1$.

It is convenient to define the twist functions $f$ and $g$ from which our twist maps will be constructed on all of $\mathbb{S}^1$, as opposed to just on $[y_0,y_1]$ and $[x_0,x_1]$ respectively (as would most naturally fit in with the abstract definition above). Let  $f:\mathbb{S}^1\to\mathbb{S}^1$ be given by
\beq*
f(y)=\left\{
\begin{array}{r@{\quad}l}
(y-y_0)/(y_1-y_0) &\text{if }y\in [y_0,y_1], \\
0 & \text{otherwise,} \end{array} \right.
\eeq*
and similarly $g:\mathbb{S}^1\to\mathbb{S}^1$ by
\beq*
g(x)=\left\{
\begin{array}{r@{\quad}l}
(x-x_0)/(x_1-x_0) &\text{if }x\in [x_0,x_1], \\
0 & \text{otherwise.} \end{array} \right.
\eeq*
Both of these functions have the form (\ref{eqn:lin_twist_fn}) illustrated in Figure~\ref{fig:twist_function}(a) (recall that $0$ and $1$ are identified in $\mathbb{S}^1$). They are differentiable for $y\in\mathbb{S}^1\backslash\{y_0,y_1\}$ and $x\in\mathbb{S}^1\backslash\{x_0,x_1\}$ respectively.

A \emph{horizontal twist map} $F:\mathbb{T}^2\to\mathbb{T}^2$ is given by $F(x,y)=(x+f(y),y)$ and it follows that $F$ is continuous on $\mathbb{T}^2$ and differentiable on $\mathbb{T}^2\backslash\partial P$. We remark that $F$ is a homeomorphism of $\mathbb{T}^2$ and is the identity map outside of $P$. We say that $F$ is \emph{linear} because of the piecewise linearity of $f$.

Analogously we define a \emph{vertical twist map} $G:\mathbb{T}^2\to\mathbb{T}^2$ by $G(x,y)=(x,y+g(x))$ and similar comments apply; in particular $G=id$ outside of $Q$.

A \emph{linear linked-twist map} $H_{j,k}:\mathbb{T}^2\to\mathbb{T}^2$ is given by the composition $G^k\circ F^j$ for positive integers $j$ and $k$. We consider the restriction of $H_{j,k}$ to the invariant set $R$. Both twist maps preserve the Lebesgue measure (see \citet{kh} or \citet{sturman}), so the composition $H_{j,k}$ does also. We denote the Lebesgue measure on $R$ by $\mu$.

If we take $x_0=y_0=0$ and $x_1=y_1=1$ and also $j=k=1$ then $H_{j,k}$ is precisely the cat map we have mentioned in Section~\ref{intro.motiv.catmap}.

Finally, let us consider $H_{j,k}$ as an abstract linked-twist map. \citet{sturman} call the map co-twisting because $jk$ is positive. We can take $F=E_1\circ T_1\circ E_1^{-1}$ where $T_1:C_1\to C_1$ is the linear twist map (defined by (\ref{eqn:lin_twist_fn})) on $C_1=\mathbb{S}^1\times[y_0,y_1]$, and where
\beq*
E_1(x,y)=\left(x,y\right).
\eeq*
Similarly we have $G=E_2\circ T_2\circ E_2^{-1}$ where $T_2:C_2\to C_2$ is the linear twist map (again defined by (\ref{eqn:lin_twist_fn})) on $C_2=\mathbb{S}^1\times[x_0,x_1]$ and where
\beq*
E_2(x,y)=\left(y,x\right).
\eeq*
The Jacobians of $E_1$ and $E_2$ are given by
\beq*
\left( \begin{array}{cc}
1 & 0 \\ 0 & 1
\end{array}\right) \quad\text{and}\quad
\left( \begin{array}{cc}
0 & 1 \\ 1 & 0
\end{array}\right)
\eeq*
respectively. The former has determinant $1$ and the latter has determinant $-1$. The fact that the signs are opposite shows that $H_{j,k}$ is co-twisting.

Chapter~\ref{chapter2} is devoted to proving the following:
\begin{thm}
If $j$ and $k$ are each at least 2, and one of them is at least 3, then the manifold $R\subset\mathbb{T}^2$ has an invariant Cantor set  on which the linked-twist map $H_{j,k}$ is topologically conjugate to a full shift on $N=(j-1)(k-1)$ symbols.
\end{thm}

\subsection{Linked-twist maps in the plane}\label{Intro.Theorems.Plane}

Linked-twist maps in the plane have been studied by a number of authors. We provide the definition and state a theorem to be proven in Chapter~\ref{chapter3}. In Section~\ref{LitReview.Ltm.Plane} we discuss the existing literature on planar linked-twist maps. In Sections~\ref{LitReview.App.Plane} and~\ref{LitReview.App.Other} we discuss some physical systems for which planar linked-twist maps provide a natural model.

When dealing with linked-twist maps in the plane it will be convenient to denote $\mathbb{S}^1=[-\pi,\pi]$ where the opposite ends of the interval are identified. Let $L$ be an annulus in the plane, centred at the origin and having inner and outer radii of $r_0$ and $r_1$ respectively (where of course $r_0<r_1$), i.e.
\beq*
L=\{(r,\theta):r_0\leqslant r\leqslant r_1\}
\eeq*
where $(r,\theta)\in\mathbb{R}^+_0\times\mathbb{S}^1$ are the usual polar coordinates. For convenience in what will follow, we assume that $r_1<\pi$. We observe that $L$ is a cylinder as in our previous discussions.

Define functions $M_{\pm}:\mathbb{R}^+_0\times\mathbb{S}^1\to\mathbb{R}^2$ by
\beq*
M_{\pm}(r,\theta)=\pm(r\cos\theta-1,r\sin\theta).
\eeq*
The images $M_{\pm}(L)$ are annuli of the `same size' in the plane, centred at $(-1,0)$ and at $(1,0)$ respectively. The annuli in the plane are expressed in Cartesian coordinates, which we denote by $(u,v)$. Let $A_{\pm}=M_{\pm}(L)$ denote these annuli.

\begin{figure}[htp]
\centering
\includegraphics[totalheight=0.25\textheight]{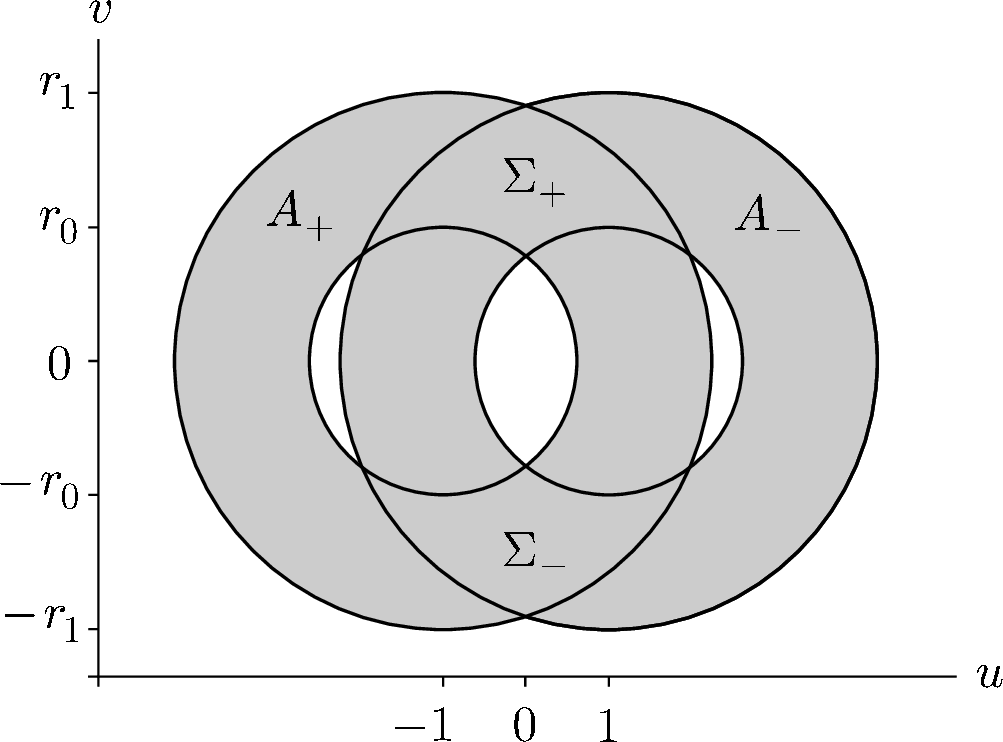}
\caption[The manifold $A\subset\mathbb{R}^2$]{The manifold $A\subset\mathbb{R}^2$ (shaded).}
\label{fig:A_plane}
\end{figure}

Under certain restrictions on $r_0,r_1$ the annuli intersect in two distinct regions; this will be a necessary though not a sufficient condition for what follows and we will say more on the sizes of annuli later. We denote the intersection region in which the $v$ coordinate is positive by $\Sigma_+$ and the other by $\Sigma_-$. See Figure~\ref{fig:A_plane}. Let $A=A_+\cup A_-$ and let $\Sigma=\Sigma_+\cup\Sigma_-$.

Inverse functions $M_{\pm}^{-1}:\mathbb{R}^2\to\mathbb{R}_0^+\times\mathbb{S}^1$ are given by
\beq*
M_{\pm}^{-1}=\left(\sqrt{(1\pm u)^2+v^2},\tan^{-1}\frac{v}{u\pm 1}\right).
\eeq*
A twist map $\Lambda:L\to L$ is defined in polar coordinates:
\beq*
\Lambda(r,\theta)=(r,\theta+2\pi(r-r_0)/(r_1-r_0)).
\eeq*
The twist function $r\mapsto 2\pi(r-r_0)/(r_1-r_0)$ has derivative $c=2\pi/(r_1-r_0)$ and is affine; as before we abuse the notation slightly and call it `linear'. It has the form (\ref{eqn:lin_twist_fn}) illustrated in Figure~\ref{fig:twist_function}(a).

\begin{figure}[htp]
\centering
(a)\includegraphics[totalheight=0.13\textheight]{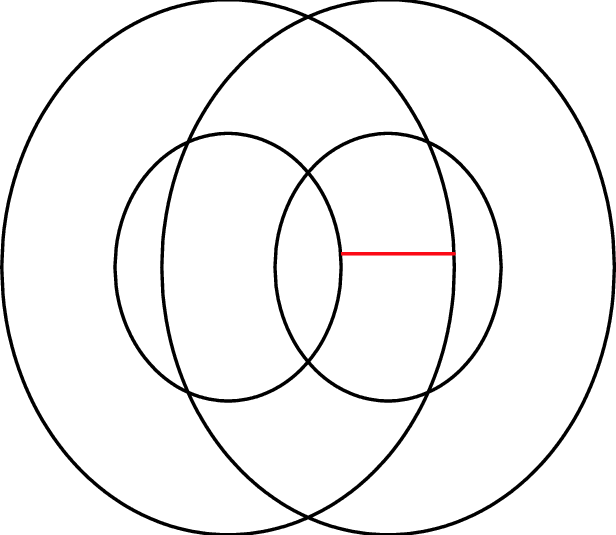}
(b)\includegraphics[totalheight=0.13\textheight]{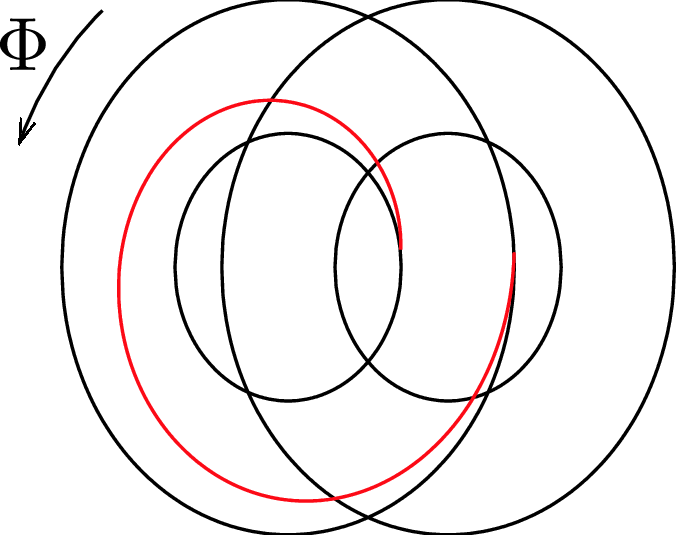}
(c)\includegraphics[totalheight=0.13\textheight]{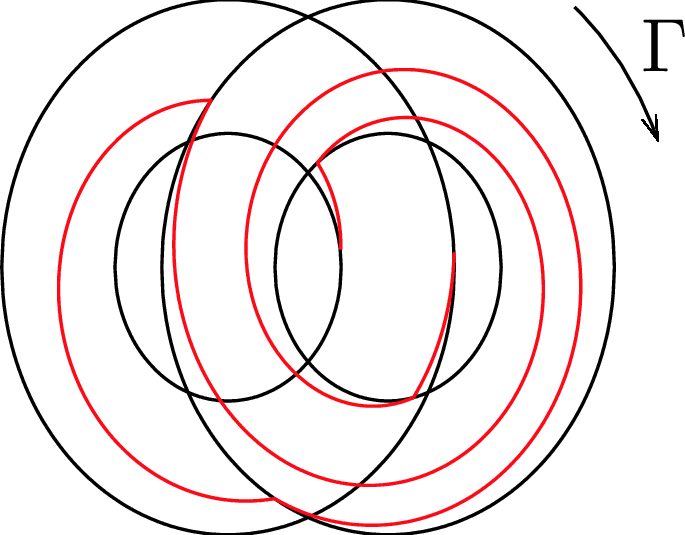}
\caption[A planar linked-twist map]{One iteration of the planar linked-twist map. Part~(a) shows some initial conditions in the form of a red horizontal line across the left-hand annulus $A_+$. Part~(b) shows the image of these points under the twist map $\Phi$ and part~(c) shows the image under the linked-twist map $\Theta=\Gamma\circ\Phi$.}
\label{fig:planar_twist}
\end{figure}

We define twist maps on $A_{\pm}$ as follows: let $\Phi,\Gamma:\mathbb{R}^2\to\mathbb{R}^2$ be given, respectively, by
\beq*
\Phi(u,v)=\left\{
\begin{array}{r@{\quad}l}
M_+\circ\Lambda\circ M_+^{-1}(u,v) &\text{if }(u,v)\in A_+ \\
(u,v) &\text{otherwise.}
\end{array}
\right.
\eeq*
\beq*
\Gamma(u,v)=\left\{
\begin{array}{r@{\quad}l}
M_-\circ\Lambda^{-1}\circ M_-^{-1}(u,v) &\text{if }(u,v)\in A_- \\
(u,v) &\text{otherwise,}
\end{array}
\right.
\eeq*
A planar linked-twist map $\Theta:A\to A$ is given by the composition $\Theta=\Gamma\circ\Phi$. Figure~\ref{fig:planar_twist} illustrates its behaviour.

We make some comments. First, we have defined $\Theta$ as the composition of one twist map $\Phi$ and one twist map $\Gamma$, in contrast to our definition of a toral linked-twist map which was the composition of $j$ `horizontal' and $k$ `vertical' twists. We can of course define a more general planar linked-twist map $\Gamma^k\circ\Phi^j$ for $j,k\in\mathbb{N}$. In fact all of the results we prove will go through with only trivial alterations; the cost however would be more cumbersome notation in a number of places. For this reason alone we take $j=k=1$.

Second, the map $\Theta$ preserves the Lebesgue measure on $A$; see \citet{woj}.

Third, let us consider the planar linked-twist map as an abstract linked-twist map. The linked-twist map $\Phi$ restricted to $A_+$ is given by $M_+\circ\Lambda\circ M_+^{-1}$ where $M_+$ is the smooth embedding of cylinder $L$ into the plane, and where $\Lambda$ denotes the twist map. $\Gamma$ restricted to $A_-$ is given by $M_-\circ\Lambda^{-1}\circ M_-^{-1}$. 

It is convenient to express $\Gamma$ in terms of $\Lambda$ rather than $\Lambda^{-1}$ (the latter not fitting our exacting definition of a twist map because the twist function has negative derivative). To this end we introduce a map $B:L\to L$ given simply by $B(r,\theta)=(r,-\theta)$; it is easy to show that $\Lambda^{-1}=B\circ\Lambda\circ B^{-1}$. Our two embeddings are thus $M_+$ and $M_-\circ B$ and we will compare the signs of the determinants of their Jacobians in order to determine whether $\Theta$ is co- or counter-twisting.
We have
\beq*
DM_{\pm}(r,\theta)=\pm\left( \begin{array}{cc}
\cos\theta & -r\sin\theta \\
\sin\theta & r\cos\theta
\end{array}\right)
\eeq*
and so determinants
\beq*
r\left(\cos^2\theta+\sin^2\theta\right),
\eeq*
which clearly are both positive. It is easy to see that $DB$ will have determinant $-1$. Thus the embedding of $L$ into $A_+$ preserves orientation whereas the embedding of $L$ into $A_-$ reverses it; $\Theta$ is co-twisting.

We illustrate the map's behaviour in Figure~\ref{fig:planar_twist}.

In Chapter~\ref{chapter3} we will prove the following:
\begin{thm}
\label{thm:main_plane}
Let $r_0=2$ and $r_1=\sqrt{7}$. Then the planar linked-twist map $\Theta:A\to A$ has the Bernoulli property, which is to say that it is isomorphic to a Bernoulli shift.
\end{thm}

This verifies a conjecture of \citet{woj}, in the particular case where the annuli are as stated.

We comment that a weakness of our method is the need to be specific about the size of the annuli. We discuss this more in Chapter~\ref{concl} where we are able to isolate which part of our proof would need to be improved upon to obtain a more general result and discuss our ideas for how this might be achieved.

\subsection{Linked-twist maps on the two-sphere}\label{Intro.Theorems.Sphere}

In this section we introduce a linked-twist map on the sphere. One of the main features of interest is the construction of the pair of embeddings for which we make use of Jacobi's elliptic functions. We will state a theorem to be proven in Chapter~\ref{chapter4}.

\begin{figure}[htp]
\centering
\includegraphics[totalheight=0.28\textheight]{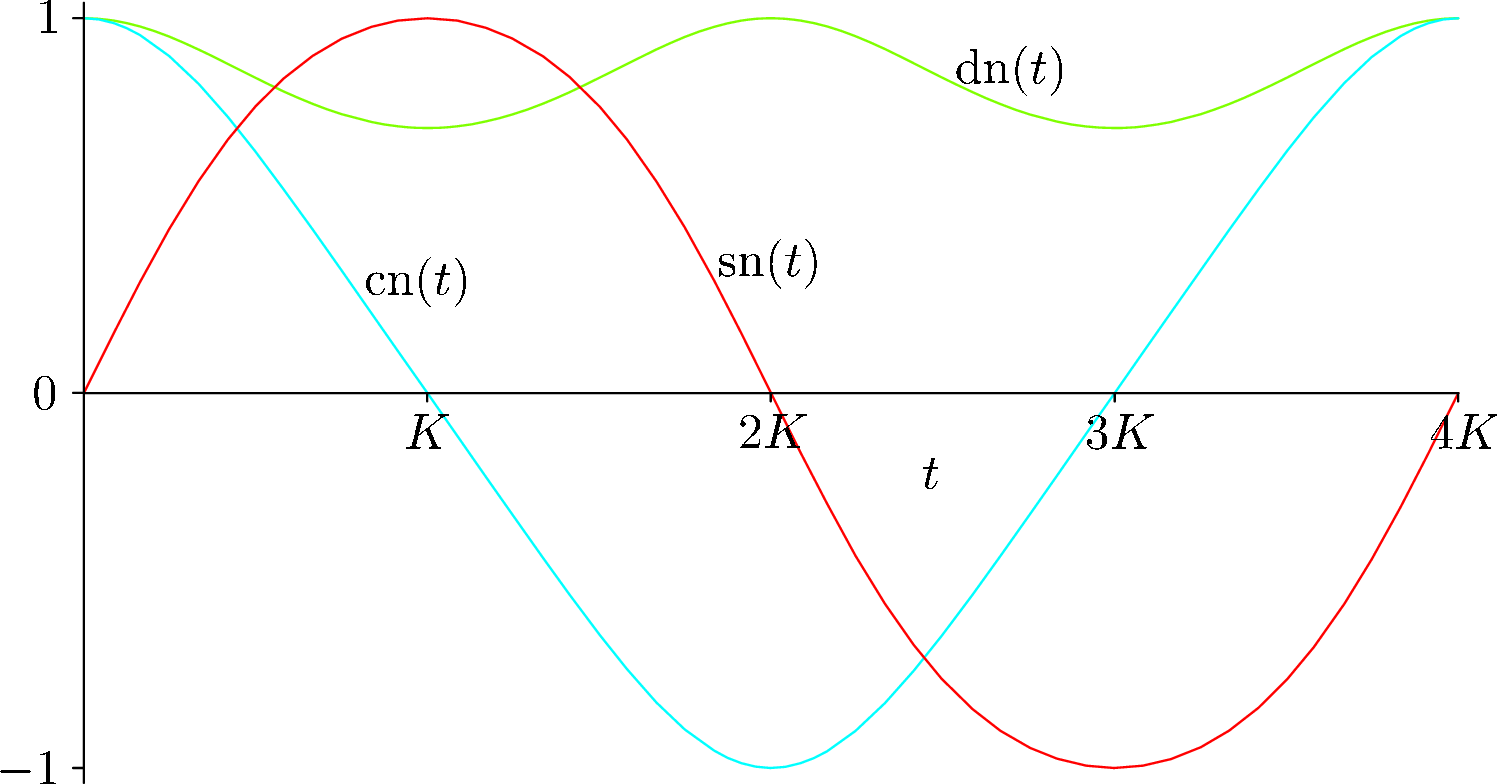}
\caption[The Jacobi elliptic functions]{The Jacobi elliptic functions: $\sn(t)$ shown in red resembles a stretched sine function, whereas $\cn(t)$ in blue resembles a stretched cosine. The function $\dn(t)$ in green has no analogy amongst the standard trigonometric functions.}
\label{fig:elliptic_functions}
\end{figure}

We review a small number of facts about Jacobi's elliptic functions; we review several more in Section~\ref{SLTM.Main.Embedding}. For a comprehensive treatment see \citet{ww} or alternatively see the excellent paper of \citet{meyer}. We plot the functions $\sn,\cn,\dn:\mathbb{R}\to[-1,-1]$ although we do not define them explicitly. Commonly each function depends upon a parameter $k\in(0,1)$ also but we will always take $k=\sqrt{2}/2$ so we omit this dependence. Let
\beq*
K(k)=\int_0^{\pi/2}\left(1-k^2\sin^2t\right)^{-1/2}\d t.
\eeq*
We comment that $K=K(\sqrt{2}/2)\approx 1.85$. Functions $\sn$ and $\cn$ are periodic with period $4K$ whereas $\dn$ is periodic with period $2K$; see Figure~\ref{fig:elliptic_functions}.

\begin{figure}[htp]
\centering
(a)\includegraphics[totalheight=0.27\textheight]{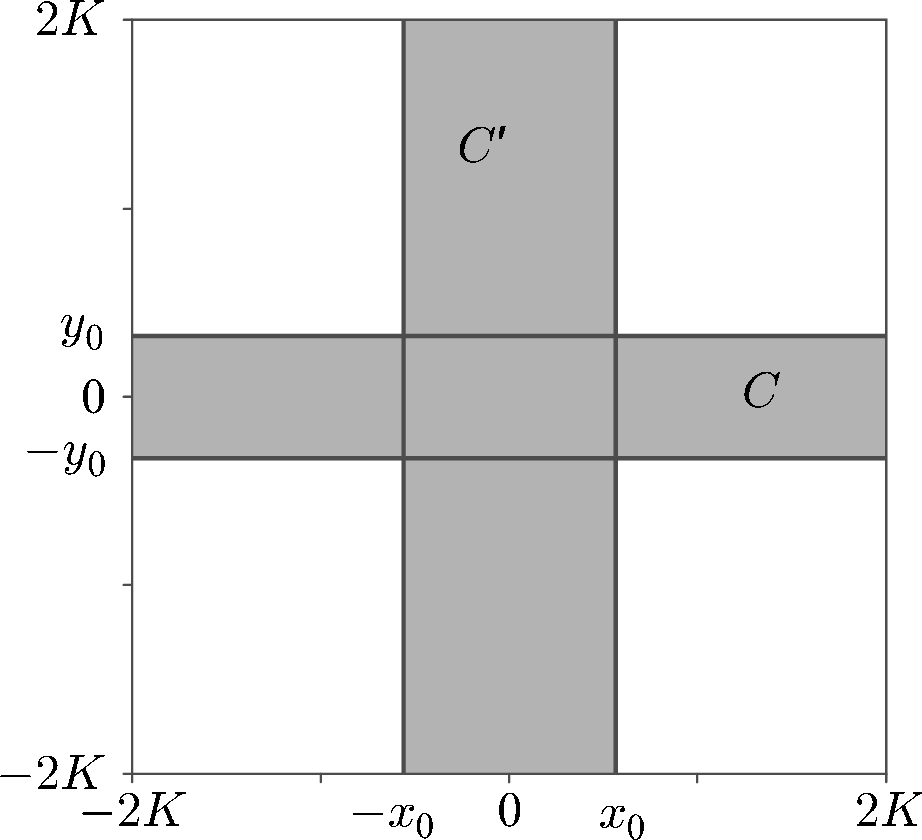}\quad
(b)\includegraphics[totalheight=0.24\textheight]{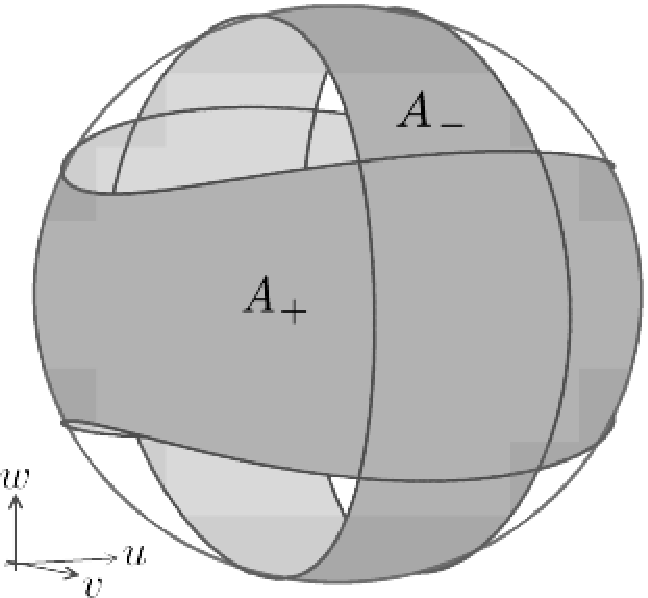}
\caption[The manifold $A\subset\mathbb{S}^2$]{It is simplest to define the manifold $A\subset\mathbb{S}^2$ as the image of $C\cup C'$ with respect to the function $E$ defined above. Notice that our restrictions on the size of $C,C'$ lead to two distinct areas of intersection of the annuli $A_+$ and $A_-$ and four `holes' which are not part of $A$. We have used dark shading for $A$ and light shading for the `reverse' side of $A$, as seen through the holes.}
\label{fig:R_sphere}
\end{figure}

It will be convenient to denote $\mathbb{S}^1=[-2K,2K]$ with the opposite ends identified. Let $C=\mathbb{S}^1\times I$ where $I=[-y_0,y_0]$ for some $0<y_0<K$. We define a second `rotated' cylinder $C'=I'\times\mathbb{S}^1$ where $I'=[-x_0,x_0]$ for some $0<x_0<K$. We can picture $C,C'$ as subsets of the two-torus $\mathbb{T}^2=\mathbb{S}^1\times\mathbb{S}^1$ as in Figure~\ref{fig:R_sphere}(a).

Now let $E:\mathbb{T}^2\to\mathbb{S}^2\subset\mathbb{R}^3$ be given by
\footnote{We are most grateful to Dr.\ Holger Waalkens for suggesting this map to us.}
\beq*
E(x,y)=\left(\sn(x)\dn(y),\cn(x)\cn(y),\dn(x)\sn(y)\right).
\eeq*
The restriction of $E$ to either $C$ or $C'$ (though \emph{not} their union) is a diffeomorphism between that cylinder and its image in $\mathbb{S}^2\subset\mathbb{R}^3$ (we prove this in Section~\ref{SLTM.Main.Embedding}). Let $A_+=E(C)$, $A_-=E(C')$ and $A=A_+\cup A_-$. Our linked-twist map will be defined on $A$ which is illustrated in Figure~\ref{fig:R_sphere}(b). (The fact that $E$ is only injective when restricted to one of $C$ or $C'$ accounts for the fact that $C\cap C'$ has one connected component but $A_-\cap A_+$ has two. We will say much more on this in Chapter~\ref{chapter4}.)

As in our construction of an abstract linked-twist map we define a twist map on $\Phi:A_+\to A_+$ to be the composition
\beq*
\Phi=E\circ T\circ E^{-1},
\eeq*
where $T:C\to C$ is a linear twist map as defined previously. We extend $\Phi$ to all of $A$ by declaring it equal to the identity function on $A_-\backslash A_+$.

One way to define a twist map $\Gamma:A_-\to A_-$ would be to define a twist map on $C'$, say $T':C'\to C'$, and define $\Gamma=E\circ T'\circ E^{-1}$. Instead we introduce a diffeomorphism $N:C\to C'$ given by
\beq*
N(x,y)=\left(\frac{x_0}{y_0}y,x\right).
\eeq*
It is easy to check that $N(C)=C'$. We define $\Gamma:A_-\to A_-$ as the composition
\beq*
\Gamma=E\circ N\circ T\circ N^{-1}\circ E^{-1}
\eeq*
and declare $\Gamma$ equal to the identity function on $A_+\backslash A_-$. The advantage of this definition is that it is trivial to show that the Jacobian matrix for $N$ has negative determinant, and this leads immediately to the conclusion that if $E$ is orientation-preserving then $E\circ N$ must be orientation-reversing, and \emph{vice versa}. Define the linked-twist map
\beq*
\Theta=\Gamma\circ\Phi,
\eeq*
then $\Theta$ is co-twisting. As with the linked-twist map in the plane we could define a more general linked-twist map on the sphere as the composition $\Gamma^k\circ\Phi^j$ for $j,k\in\mathbb{N}$. Our work would again require only trivial alterations but at a cost of more cumbersome notation.

In Chapter~\ref{chapter4} we will prove the following.

\begin{thm}\label{thm:main_sphere}
The linked-twist map $\Theta:A\to A$ has the Bernoulli property, which is to say that it is isomorphic to a Bernoulli shift.
\end{thm}

%% file: chapter1/chapter1.tex
\chapter{Literature Review}
\setcounter{equation}{0}
\setcounter{figure}{0}
\label{chapter1}

The ergodic theory of hyperbolic systems is a significant branch of the dynamical systems theory and the literature is appropriately rich and diverse. In this chapter we provide some key definitions and results on which our work builds, but in doing so we barely scratch the surface of all that is out there.

The chapter is divided into four sections. In Section~\ref{LitReview.ErgTh} we provide some basic definitions from ergodic theory, starting with the relatively weak property of \emph{ergodicity} and building up to the strongest \emph{Bernoulli} property. In Section~\ref{LitReview.Hyp} we review some key concepts and results from the \emph{hyperbolic} theory. The systems we study in later chapters will all display \emph{non-uniform hyperbolicity} and here we describe in detail what this means.

In Section~\ref{LitReview.Ltm} we survey the literature pertaining to the \emph{linked-twist maps} we described in the previous chapter. Finally in Section~\ref{LitReview.App} we survey a few examples of applications for which linked-twist maps provide a natural model. The existence of such applications goes some way to explaining the recent resurgence in interest in linked-twist maps, as is perhaps best evidenced by the book of \citet{sturman}.

\input chapter1/section1

\input chapter1/section2

\input chapter1/section3

\input chapter1/section4

%% file: chapter1/section1.tex
\section{Ergodic theory}\label{LitReview.ErgTh}
\setcounter{equation}{0}

Ergodic theory is concerned with dynamical systems on measure spaces. It is typically highly non-trivial to prove that a given dynamical system has any of the ergodic properties we will present. However, the pay-off for doing so is a substantial amount of information about the behaviour of `most' trajectories.

\subsection{Ergodicity and mixing}\label{LitReview.ErgTh.Erg}

All definitions and results in this section may be found in \citet{bs}. Another standard reference for this material is \citet{kh}.

Let $(M,\mathfrak{U},f,\mu)$ be a measure-preserving dynamical system. Here $M$ is a set and will usually be furnished with some additional structure; typically we might require that $M$ be a compact metric space, or a Riemannian manifold. $\mathfrak{U}$ denotes a $\sigma$-algebra of subsets of $M$, $f$ a transformation of $M$ into itself and $\mu$ a positive measure defined on $\mathfrak{U}$.

Typically $\mu$ will be finite and so without loss of generality we may assume that it is a probability measure, i.e.\ $\mu(M)=1$. We will assume that $f$ preserves $\mu$, in the sense that for each set $A\in\mathfrak{U}$ we have $\mu(f^{-1}(A))=\mu(A)$.

\begin{defn}[Ergodicity]
A dynamical system $(M,\mathfrak{U},f,\mu)$ is said to be ergodic if whenever $A\in\mathfrak{U}$ has the property that $f(A)=A$, then either $\mu(A)=0$ or $\mu(A)=1$.
\end{defn}

Ergodicity may be thought of as \emph{indecomposability}, in the sense that two disjoint, non-trivial (i.e.\ positive measure) invariant sets are not possible.

A stronger condition than ergodicity is the following:

\begin{defn}[Strong mixing]
A dynamical system $(M,\mathfrak{U},f,\mu)$ is said to be strong mixing if for all sets $A,B\in\mathfrak{U}$ one has
\beq*
\lim_{n\to\infty}\mu\left(f^{-n}(A)\cap B\right)=\mu(A)\mu(B).
\eeq*
\end{defn}

The strong mixing (typically called just \emph{mixing}) property implies ergodicity and can be thought of as points `losing memory' of where they started. This is the kind of property we would like to prove for our dynamical systems.

In fact, it will be possible to prove a stronger property known as the \emph{Bernoulli property}. The Bernoulli property is significantly more abstract than the other ergodic properties we have presented but the pay-off is substantial; Bernoulli systems behave, in a rigorous sense, as randomly as possible.

\subsection{The Bernoulli property}\label{LitReview.ErgTh.Bern}

A good reference for the material in this section is \citet{wig1}. Let $S=\{1,2,...,N\}$ be a collection of $N$ symbols, where $N$ is an integer strictly greater than one. A (bi-infinite) symbol sequence has the form $s=...,s_{-1},s_0,s_1,...$ where each $s_i\in S$. The space $\Sigma^N$ of all such symbol sequences is naturally thought of as the bi-infinite Cartesian product $\cdots\times S\times S\times S\times\cdots$. We can define a metric on $\Sigma^N$: if $t=...,t_{-1},t_0,t_1,...$ is another symbol sequence then let
\beq*
d(s,t)=\sum_{i=-\infty}^{+\infty}\frac{\delta_i}{2^{|i|}}\qquad\text{where }\delta_i=\left\{
\begin{array}{r@{\quad}l}
0 &\text{if }s_i=t_i \\
1 &\text{otherwise.}
\end{array}
\right.
\eeq*
See \citet{devaney} for a proof that $d$ is indeed a metric on $\Sigma^N$. Intuitively points in $\Sigma^N$ are close if their sequences agree on a long central block. It is shown in \citet{sturman} that the metric space $(\Sigma^N,d)$ is compact, totally disconnected and perfect (i.e.\ a \emph{Cantor set}) and has the cardinality of the continuum.

We now outline how one may define a measure on $\Sigma^N$, following the approach of \citet{arnoldavez}. Let $A_i^j$ be the set of points in $\Sigma^N$ having $j\in S$ for the $i^{\text{th}}$ element in the symbol sequence. These sets generate a $\sigma$-algebra of subsets of $\Sigma^N$. We also define the \emph{cylinder sets}
\beq*
A^{j_1\cdots j_k}_{i_1\cdots i_k}=\bigcap_{h=1}^k A^{j_h}_{i_h}.
\eeq*
Define a normalised measure $\mu$ on $S$ by insisting that for each $j\in S$ we have $\mu(j)\geqslant 0$ and $\sum_{j=1}^N\mu(j)=1$. The measure of a set $A_i^j$ is defined by $\mu\left(A^j_i\right)=\mu(j)$ and we extend this measure to the cylinder sets via the identity
\beq*
\mu\left(A^{j_1\cdots j_k}_{i_1\cdots i_k}\right)=\prod_{h=1}^k\mu(j_h).
\eeq*
It can be shown that $\mu$ satisfies the axioms of a measure.

The last part of our construction is a map $\sigma$ of $\Sigma^N$ into itself, known as a \emph{shift map}. It is expressed most concisely by the relationship $[\sigma(s)]_i=s_{i+1}$, although perhaps intuitively it is preferable to insert a period at some point in the symbol sequence (written without commas), and look at where that period occurs in the symbol sequence of the image:
\beq*
s=\cdots s_{-2}s_{-1}.s_0s_1s_2\cdots,\qquad\sigma(s)=\cdots s_{-2}s_{-1}s_0.s_1s_2\cdots.
\eeq*
If the domain is all of $\Sigma^N$ then $\sigma$ is often called a \emph{full shift on $N$ symbols}. In this case it can be shown (see \citet{wig1}) that $\sigma$ is a homeomorphism of $\Sigma^N$, that it has a countable infinity of periodic orbits including orbits of all periods, an uncountable infinity of non-periodic orbits, a dense orbit and moreover (see \citet{sturman}) that $\sigma$ preserves the measure $\mu$ constructed previously, with respect to which it is mixing.

\begin{defn}[Bernoulli property]
A dynamical system $(M,\mathfrak{U},f,\nu)$ is said have the Bernoulli property if it is (metrically) isomorphic to a Bernoulli shift. More formally, we require that the following diagram commutes
\begin{center}
\begin{pspicture}(0,.5)(3,3) 
\psset{nodesep=5pt} 
\rput(.5,2.5){\rnode{A}{$\Sigma^N$}} 
\rput(2.5,2.5){\rnode{B}{$\Sigma^N$}} 
\rput(.5,1){\rnode{C}{$M$}} 
\rput(2.5,1){\rnode{D}{$M$}}
\ncline{->}{A}{B}\Aput{$\sigma$}
\ncline{->}{B}{D}\Aput{$\phi$} 
\ncline{->}{A}{C}\Bput{$\phi$} 
\ncline{->}{C}{D}\Bput{$f$}
\end{pspicture} 
\end{center}
where $\phi:(\Sigma^N,\mu)\to(M,\nu)$ is an isomorphism.
\end{defn}

A map having the Bernoulli property automatically has all of the properties of the full shift on $N$ symbols, given above. An example of a dynamical system having the Bernoulli property is the well-known \emph{baker's map} of the unit square. See \citet{sturman} for further details of the map and a proof of this result. In this example the isomorphism can be constructed explicitly; this is aided by discontinuities in the map, in contrast to the maps we will consider.

%% file: chapter1/section2.tex
\section{Hyperbolicity}\label{LitReview.Hyp}
\setcounter{equation}{0}

Hyperbolicity is an important part of the dynamical systems theory and the focus of a great deal of active research. Knowing that a certain dynamical system has hyperbolic structure gives us access to a number of results and techniques for demonstrating ergodic properties. All of the dynamical systems we consider in this thesis display hyperbolicity; here we outline some of the key definitions and results we will repeatedly rely upon.

\subsection{Uniform hyperbolicity}\label{LitReview.Hyp.UH}

Let $(M,\mathfrak{U},f,\mu)$ be a dynamical system. Throughout this section we will assume that $M$ has some structure beyond being just a set, likewise $f$ some smoothness properties. In particular we assume that $M$ is a compact, smooth ($C^{\infty}$) $n$-dimensional Riemannian manifold and $f$ a smooth ($C^{\infty}$) diffeomorphism of $M$ (i.e.\ a differentiable map with differentiable inverse; we have defined these terms in the case $n=2$ in Section~\ref{intro.abs.rev}). In section~\ref{LitReview.Hyp.smws} we will discuss what happens when we relax these conditions somewhat, but for now let us keep the exposition as clean as possible. Denote by $Df_x$ the Jacobian matrix of $f$ evaluated at $x\in M$.

\begin{defn}[Hyperbolic fixed point]
Let $x\in M$ be a fixed point of $f$, i.e.\ $f(x)=x$. Then $x$ is said to be hyperbolic if none of the eigenvalues of $Df_x$ have magnitude one.
\end{defn}

In some neighbourhood $U$ of a hyperbolic fixed point $x\in M$, the dynamics of $f$ will closely resemble the behaviour of the linearised system. To be precise, there is a corresponding neighbourhood $V$ of the origin and a homeomorphism $h:V\to U$ such that $f(h(y))=h(Df_xy)$ for all $y\in V$. This is the well-known Hartman-Grobman theorem; see \citet{robinson}.

More generally we will define hyperbolicity on a set, rather than at a fixed point. The simplest case is where the hyperbolicity is \emph{uniform}. The following definition is taken from \citet{bs}.

\begin{defn}[Uniformly hyperbolic set, Anosov diffeomorphism]
Let $(M,\mathfrak{U},f,\mu)$ be a measure-preserving dynamical system. Let $U\subset M$ be a non-empty open subset and $f:U\to f(U)\subset M$ a smooth diffeomorphism. A compact, $f$-invariant set $\Lambda\subset U$ is said to be uniformly hyperbolic if there exist constants $c>0$ and $0<\lambda<1$, and if there is a continuous splitting of the tangent space $T_xM=E^s(x)\oplus E^u(x)$ at each $x\in\Lambda$ such that
\beql\label{eqn:UH1}
Df_xE^s(x)=E^s(f(x))\and Df_xE^u(x)=E^u(f(x)),
\eeql
\beql\label{eqn:UH2}
\|Df_x^nv_s\|\leqs c\lambda^n\|v_s\|\quad\text{ for }v_s\in E^s(x)
\eeql
\beql\label{eqn:UH3}
\|Df_x^{-n}v_u\|\leqs c\lambda^n\|v_u\|\quad\text{ for }v_u\in E^u(x).
\eeql
If $\Lambda=M$ then $f$ is called an Anosov diffeomorphism.
\end{defn}

Condition (\ref{eqn:UH1}) says that stable and unstable directions should be invariant under the differential map $Df_x$, whereas conditions (\ref{eqn:UH2}) and (\ref{eqn:UH3}) give estimates on the contraction of stable subspaces under forward iteration, and of unstable subspaces under backward iteration respectively. We call this hyperbolicity `uniform' because the constants $c$ and $\lambda$ are independent of the point $x\in M$. Notice that this is precisely the situation we found when analysing the cat map in Section~\ref{intro.motiv.catmap}.

It turns out that few dynamical systems are uniformly hyperbolic. In the next section we discuss a generalisation.

\subsection{Pesin theory}\label{LitReview.Hyp.Pesin}

We begin by defining what we mean by \emph{non-uniformly hyperbolic}, then go on to describe some celebrated results due to \citet{pes} which have hugely influenced the study of such systems over the past three decades or so. \citet{bap} and \citet{sturman} both give good accounts of the material in this section; our definitions are taken from these sources.

We will not use the results from this section as they do not apply to the linked-twist maps we study, which are not diffeomorphisms. However they provide a natural bridge from uniform hyperbolicity to the theorem of \citet{ks} that we will describe next and use extensively thereafter.

\begin{defn}[Non-uniformly hyperbolic]
The measure-preserving dynamical system $(M,\mathfrak{U},f,\mu)$ is said to be non-uniformly (completely) hyperbolic if there exist measurable functions $0<\lambda_-(x)<1<\lambda_+(x)$ and $\eps(x)$ such that $\eps(x)=\eps(f(x))$ (i.e.\ $\eps$ is invariant along trajectories), and if there is a splitting of the tangent space $T_xM=E^s(x)\oplus E^u(x)$ for each $x$, and finally a function $c(x)$ so that, for each $k\in\mathbb{Z}$ and $n>0$ we have
\beql\label{eqn:NUH1}
Df_x^kE^s(x)=E^s(f^k(x))\and Df_x^kE^u(x)=E^u(f^k(x)),
\eeql
\beql\label{eqn:NUH2}
\|Df_x^nv_s\|\leqs c\left(f^k(x)\right)\lambda_-^n(x)\|v_s\|\quad\text{ for }v_s\in E^s(x)
\eeql
\beql\label{eqn:NUH3}
\|Df_x^nv_u\|\geqs c^{-1}\left(f^k(x)\right)\lambda_+^n(x)\|v_u\|\quad\text{ for }v_u\in E^u(x),
\eeql
\beql\label{eqn:NUH4}
\angle\left(E^s(x),E^u(x)\right)\geqs c^{-1}(x),
\eeql
\beql\label{eqn:NUH5}
c\left(f^k(x)\right)\leqs c(x)e^{\eps(x)|k|}.
\eeql
\end{defn}

Conditions (\ref{eqn:NUH1}), (\ref{eqn:NUH2}) and (\ref{eqn:NUH3}) are analogous to the conditions we impose on uniformly hyperbolic systems, though the replacement of independent constants with functions of $x$ means that they are less restrictive. Condition (\ref{eqn:NUH4}) says that the stable and unstable directions are transversal. The final condition (\ref{eqn:NUH5}) is perhaps a little more subtle and deals with the rate at which our contraction or expansion estimates in conditions (\ref{eqn:NUH2}) and (\ref{eqn:NUH3}) deteriorate along a trajectory. It says that this deterioration is sub-exponential and is thus dominated by the exponential contraction or expansion.

The most important tool for analysing non-uniformly hyperbolic systems is the Lyapunov exponent.

\begin{defn}[Lyapunov exponent]
For a dynamical system $(M,\mathfrak{U},f,\mu)$, the Lyapunov exponent $\chi^{\pm}(x,v)$ at the point $x\in M$ and in the direction $v\in T_xM$ is given by
\beq*
\chi^{\pm}(x,v)=\lim_{n\to\pm\infty}\frac{1}{n}\log\left\|Df_x^nv\right\|,
\eeq*
whenever this limit exists, where $\|\cdot\|$ denotes the standard Euclidean norm in $\mathbb{R}^n$.
\end{defn}

The importance of Lyapunov exponents is illustrated by the fact that non-uniformly hyperbolic systems are commonly known as \emph{systems with non-zero Lyapunov exponents}. The result to which this epithet alludes is the following.

\begin{thm}[\citet{pes}]
A dynamical system $(M,\mathfrak{U},f,\mu)$ is non-uniformly (completely) hyperbolic if for almost every $x\in M$ the Lyapunov  exponent $\chi(x,v)$ is non-zero for every non-zero $v\in T_xM$.
\end{thm}

Pesin derived two further results which lay the foundations for our results on linked-twist maps. The first is the famous stable manifold theorem. We need the following definition, taken from \citet{bap}.

\begin{defn}[Local invariant manifolds]
Let $B^s(0,\varepsilon)$ be the open $\varepsilon$-neighbourhood of the origin in $E^s(x)$, and similarly $B^u(0,\varepsilon)$. A local stable manifold of $x$ has the form
\beql
\label{eqn:st}
\gamma^s(x)=\exp_x\{(x,\psi^s(x)):x\in B^s(0,\varepsilon)\}
\eeql
for some $\varepsilon>0$, where $\psi^s:B^s(0,\varepsilon)\to E^u(x)$ is a smooth map satisfying $\psi^s(0)=0$ and $D\psi^s(0)=0$. The trajectories of $x$ and $y\in\gamma^s(x)$ approach each other at exponential rate as $n\to+\infty$.
Transposing $u$ and $s$ in the above and considering $n\to-\infty$ yields the description of local unstable manifolds; we omit further details.
\end{defn}

\begin{thm}[\citet{pes}]
If $f:M\to M$ is of class $C^{1+\alpha}$ and is non-uniformly hyperbolic, then for almost every $x\in M$ there exists a local stable manifold $\gamma^s(x)$ with the properties that $x\in\gamma^s(x)$, $T_x\gamma^s(x)=E^s(x)$ and if $y\in\gamma^s(x)$ and $n\geqs 0$ then
\beq*
\rho\left(f^n(x),f^n(y)\right)\leqs T(x)\lambda^ne^{\eps n}\rho(x,y),
\eeq*
where $\rho$ is the distance in $M$ induced by the Riemannian metric and $T:M\to(0,\infty)$ is a Borel function satisfying
\beq*
T\left(f^m(x)\right)\leqs T(x)e^{10\eps|m|},\quad m\in\mathbb{Z}.
\eeq*
\end{thm}

The other result that will be crucial to our work shows that the manifold $M$ has an \emph{ergodic partition}, the definition of which is given in the following statement.

\begin{thm}[\citet{pes}]
\label{thm:pes_erg_part}
If $f:M\to M$ is of class $C^{1+\alpha}$ and is non-uniformly hyperbolic then $M$ is either a finite or countably infinite union of disjoint measurable sets $M_0,M_1,...$ such that $\mu(M_0)=0$ and $\mu(M_i)>0$ for all other subsets, each $M_i$ is $f$-invariant (i.e.\ $f(M_i)=M_i$) and the restriction of $f$ to any $M_i$ is ergodic.
\end{thm}

This completes our very brief exposition of Pesin theory. As we have mentioned, in order to use Pesin-type results as above we need to appeal to more general work of \citet{ks} which we review in the following section.

\subsection{Smooth maps with singularities}\label{LitReview.Hyp.smws}

This section surveys the results of \citet{ks}, our exposition following closely that found in the appendix of \citet{p1}. In stating the results it is necessary to introduce several nested full-measure sets. We tabulate these (Table~2.1) in the hope that it helps the reader through the construction.

Let $X$ be a complete metric space with a metric $\rho$. Let $M\subset X$ be an open subset which is also a Riemannian manifold, with the Riemannian metric inducing $\rho|_M$.  Assume there is some $r>0$ such that for each $x\in M$ the exponential map $\exp_x$, restricted to the ball
\beq*
B(x)=B(x,\min(r,\text{dist}_{\rho}(x,X\backslash M))),
\eeq*
is injective.

Let $(M,\mathfrak{U},f,\mu)$ be a measure-preserving dynamical system as before, but with the difference that we define $f$ only on an open set $N\subset M$ into $M$. The measure $\mu$ is an $f$-invariant probability measure on $M$ and we require that, on $N$, the function $f$ is $C^2$ and injective. Finally we denote $\text{sing}(f)=M\backslash N$.

We call $f$ a \emph{smooth map with singularities}. This completes the definition of the map itself. We now describe two conditions of a technical nature that place restrictions on the nature of the points at which $Df$ is undefined, and on the growth of $D^2f$ near to these points.

\begin{table}
\centering\label{tab:X}
\begin{tabular}{|c|l|}\hline
Set&Description\\ \hline
$X$& Complete metric space\\
$M$& $n$-dimensional Riemannian manifold\\
$N$& Open set on which $f$ is defined\\
$J$& Intersection of all images and pre-images of $N$\\
$E$& Set on which Lyapunov exponents exist\\ \hline
\end{tabular}
\caption{Full measure subsets of the complete metric space X. As we have listed the sets, each contains those below it.}
\end{table}

In keeping with the notation of \citet{p1} and \citet{sturman} we say that $f$ satisfies the condition (KS1) if and only if there are positive constants $a$ and $c_1$ so that for every $\eps>0$ we have
\beql\label{eqn:KS1}
\mu\left(B(\text{sing}(f),\varepsilon)\right)\leqslant c_1\varepsilon^a,
\eeql
where $B(\text{sing}(f),\varepsilon)$ means the open $\varepsilon$-neighbourhood, with respect to $\rho$, of the set $\text{sing}(f)$.

We say that $f$ satisfies the condition (KS2) if and only if there are positive constants $b$ and $c_2$ so that for every $x\in N$ we have
\beql\label{eqn:KS2}
\|D^2f(x)\|\leqslant c_2\rho(x,\text{sing}(f))^{-b},
\eeql
where $\|D^2f(x)\|$ denotes the supremum of $\|D^2(\exp_z^{-1}\circ f\circ\exp_y)\|$, taken over those $x$ for which $x\in B(y)$ and $f(x)\in B(z)$.

Informally, the two conditions state that `most' points have neighbourhoods free from singularities and that the second derivative does not get large `too quickly' (e.g.\ exponentially) as we approach these
singularities.

If $f$ satisfies condition (KS1) then $\mu(\text{sing}(f))=0$ and so $\mu(N)=1$. Let $J=\bigcap_{n=-\infty}^{\infty}f^n(N)$. The $f$-invariance of $\mu$ implies that $\mu(J)=1$ also (we prove this in Chapter~\ref{chapter4}). The Multiplicative Ergodic Theorem (originally proven by \citet{osel} but see also \citet{cm} and \citet{bap}) holds for smooth maps with singularities satisfying the above conditions. We have

\begin{thm}[\citet{cm}]
\label{thm:met}
Suppose that
\beql
\label{eqn:oseledec}
\int_X\log^+\|Df\|_{\text{op}}d\mu<\infty\qquad
\text{and}\qquad\int_X\log^+\|Df^{-1}\|_{\text{op}}d\mu<\infty,
\eeql
where $\log^+(\cdot)=\max\{\log(\cdot),0\}$ and $\|\cdot\|_{\text{op}}$ denotes the operator norm induced by $\rho$. Then there is an $f$-invariant set $E\subset J$, $\mu(E)=1$, so that for every $x\in E$ and every non-zero $v\in T_xX$, the Lyapunov exponents $\chi^{\pm}(x,v)$ exist.
\end{thm}

The following theorem provides the framework for our work in Chapters~\ref{chapter3} and~\ref{chapter4}; in addition to the main reference, see also \citet{p1} and \citet{sturman}. The theorem contains the definition of an \emph{ergodic partition} to which we will refer several times, and includes a more complete description than was given in Theorem~\ref{thm:pes_erg_part}.

\begin{thm}[\citet{ks}]
\label{thm:ks}
Let $f$ be a smooth ($C^2$ at least) map with singularities, defined a.e.\ on smooth manifold $M$ as above.
\begin{enumerate}
\item[\emph{(a)}] Suppose $f$ satisfies the conditions (KS1) and (KS2) and the hypothesis of Theorem~\ref{thm:met} above. Then for a.e.\ $x\in M$ and for all non-zero tangent vectors $v\in T_xM$, the Lyapunov exponents $\chi^{\pm}(x,v)$ exist. Corresponding to any positive (respectively, negative) Lyapunov exponents, there exist local unstable (stable) manifolds
$\gamma^{u(s)}(x)$ of the form we have described.
\item[\emph{(b)}] If additionally for a.e.\ $x\in M$ and for all non-zero $v\in T_xM$ we have $\chi(x,v)\neq 0$ then $M$ decomposes into an (at most) countable family of positive measure, $f$-invariant, pairwise-disjoint sets $M_i$ on which the restriction of $f$ is ergodic. Furthermore each set $M_i$ has the form $M_i=\bigcup_{j=1}^{n(i)}M_i^j$ where, for each $j$, $f^{n(i)}|_{M_i^j}$ is Bernoulli. Such a system will be said to have an ergodic partition.
\item[\emph{(c)}] If additionally for a.e.\ $x,y\in M$ there exist integers $m,n$ such that
\beql
\label{eqn:mip1}
f^m(\gamma^u(x))\cap f^{-n}(\gamma^s(y))\neq\emptyset
\eeql
(we say that $f$ satisfies the manifold intersection property), then in the decomposition of $M$ there is just one positive-measure set, i.e.\ $M_i=M$ for all $i$, and so $f$ is ergodic.
\item[\emph{(d)}] If additionally, for a.e.\ $x,y\in M$, the condition (\ref{eqn:mip1}) is satisfied for each pair of sufficiently large integers $m,n$ (we say that $f$ satisfies the repeated manifold intersection property), then $f$ has the Bernoulli property.
\end{enumerate}
\end{thm}

\subsection{The Sinai-Liverani-Wojtkowski approach to proving ergodicity}\label{LitReview.Hyp.SLW}

We discuss some work of \citet{liverani95ergodicity} based on work of \citet{sinai_70} in which the authors establish criteria for certain maps to be ergodic. The class of systems to which their results apply is broad and includes some (according to \citet{nicol2}) of the linked-twist maps we will study, if not all of them. There are a large number of technical hypotheses in the statement of their theorem and we do not intend to use their results, so we limit our exposition to an informal discussion.

In part we do not appeal to their result because the aforementioned technical considerations mean that this would be far from trivial (although this alone is not justification, as the same can be leveled at our approach!). However, when we come to make our closing remarks and discuss our ideas for future works we will have good cause to conjecture that our approach offers benefits that theirs does not.

The significance of the works of \citet{pes} and later of \citet{ks} is that they extended the class of systems for which an ergodic partition (a \emph{local} property) can be established. Conversely \citet{sinai_70} and later \citet{liverani95ergodicity} extend the class of systems for which ergodicity (a \emph{global} property) can be established.

Consider the condition (\ref{eqn:mip1}) given in Theorem~\ref{thm:ks}, which is the bridge between local and global properties in that theorem. The condition is sufficient because it allows one to conclude that any integrable observable on $M$ (that is, any function $g\in L^1(M,\mathbb{R})$) that is $f$-invariant, is constant $\mu$-almost everywhere. It can be shown that this condition is equivalent to ergodicity; see for example \citet{bs}. This idea dates back to \citet{hopf} and such a construction is sometimes called a \emph{Hopf chain}; for more details we recommend the introductory sections of \citet{liverani95ergodicity}.

Proving that (\ref{eqn:mip1}) is satisfied, in general, will require some specific knowledge of the nature of local invariant manifolds. If we are able to conclude that the length of $f^m(\gamma^u(x))$ diverges with $m$ and that the length of $f^{-n}(\gamma^s(y))$ diverges with $n$, and moreover if we can give a `useful' characterisation of their respective \emph{orientations}, then we might hope to satisfy the condition. Either, or both, of these may be far from trivial within a generic non-uniformly hyperbolic system.

As we have seen, such systems do not have uniform growth rates for local invariant manifolds, nor indeed uniform lower bounds on the original sizes of those manifolds. Thus in general one needs a more sophisticated approach, and this is precisely where the Sinai-Liverani-Wojtkowski approach comes in.

The Sinai-Liverani-Wojtkowski method (as we call it) works by constructing a connected `chain' of local invariant manifolds with one end at $x$ and the other at $y$. So far this is just Hopf's method, but rather than relying on growth and orientation to deduce this connection, their method is to very carefully partition $M$ using small overlapping `squares' whose sides are parallel to the stable and unstable directions respectively.

Their arguments relate the width of a chosen partition to the conditional measure of those points within a given square whose local invariant manifolds completely cross that square. In this manner they are able to conclude that if a point $x\in M$ satisfies a certain local condition, then there is an open set containing $x$ that is itself contained within a single ergodic component. It is an easy corrolary that if $\mu$-a.e.\ $x\in M$ satisfies this condition then the map is ergodic.

As we have briefly argued, the real achievement of these methods is to deduce ergodicity in systems where either the growth or orientation (or both) of local invariant manifolds are not well-behaved, in some sense. The linked-twist maps for which we prove strong ergodic properties do not fall into this category; we think in particular of the linked-twist map in the plane: that local invariant manifolds grow arbitrarily long is known and mentioned elsewhere; that the orientations of these manifolds can be characterised in a useful manner is the cornerstone of our proof.

%% file: chapter1/section3.tex
\section{Linked-twist maps}\label{LitReview.Ltm}
\setcounter{equation}{0}

The linked-twist map literature spans almost three decades and includes results describing certain ergodic properties of the toral and the planar linked-twist maps we have mentioned. In Section~\ref{LitReview.Ltm.Torus} we review the relevant results for toral maps and in Section~\ref{LitReview.Ltm.Plane} we do the same for the planar maps. For a comprehensive overview of this literature the reader is directed to \citet{sturman}. In Section~\ref{LitReview.Hyp.Other} we list some other explicitly defined maps for which strong ergodic properties have been established, in the hope that this will help to put the results for linked-twist maps into context.

\subsection{Linked-twist maps on the two-torus}\label{LitReview.Ltm.Torus}

Let $H=G^k\circ F^j:R\to R$, with $R\subset\mathbb{T}^2$, be a toral linked-twist map as defined in Section~\ref{Intro.Theorems.Torus}. Recall that the product $jk$ is positive if and only if a toral linked-twist map $H_{j,k}$ is co-twisting. The following theorems describe the ergodic properties of these maps. 

\begin{thm}[\citet{be}]
If $jk>0$ and $H$ is composed of smooth twists, then $H$ has an ergodic partition.
\end{thm}

The smooth twists are as depicted in Figure~\ref{fig:twist_function}(b). Furthermore the authors sketch a geometrical argument with which this may be extended to the Bernoulli property.

\begin{thm}[\citet{d1}]
If $jk>0$ then periodic and homoclinic points of $H$ are dense, and $H$ is topologically mixing.
\end{thm}

Devaney's result is more topological in nature; in fact this theorem was motivated by the similarities between toral \ltms and the cat map. We observe also that Devaney does not need to put restrictions on the nature of the twist functions (aside, of course, from those conditions mentioned in Section~\ref{Intro.Background.Abstract}, which all twist functions must satisfy).

\begin{thm}[\citet{woj}]
If $jk>0$ and $H$ is composed of linear twist maps, then $H$ is Bernoulli. Alternatively, if $jk<-4$ then $H$ has an ergodic partition.
\end{thm}

Wojtkowski's result on the torus is the main source of inspiration for our proof of the Bernoulli property in the plane. In fact he proved only the $K$-property; it follows from \citet{ch} that the system is necessarily Bernoulli. Wojtkowski also considers planar linked-twist maps in the same paper; we mention this in the next section.

In order to state the next result we briefly describe what is meant by the \emph{strength} of a twist. Consider a `horizontal' twist map $F^j:P\to P$ as defined in Section~\ref{Intro.Theorems.Torus}, i.e.\ $F^j(x,y)=(x+jf(y),y)$ where $f:I\to\mathbb{S}^1$ is a twist function. The strength of the twist map $F^j$ is defined to be 
\beq*
s_F=\sgn(j)\inf_{y\in I}|jf'(y)|.
\eeq*
Denote $I=[i_0,i_1]$. In the common case where $f$ is affine then $s_F=j/(i_1-i_0)$. The strength of a `vertical' twist map is defined analogously.

\begin{thm}[\citet{p1}]
If $s_Fs_G<-C_0\approx -17.24445$ and $|j|,|k|\geqs 2$ then $H$ is Bernoulli.
\end{thm}

Przytycki constructs an intricate argument allowing him to prove this result for counter-twisting maps. We have mentioned before that these are more difficult to study than their co-twisting counterparts; a comparison of the criteria in this theorem and the previous one exemplifies the situation.

Two other results for toral \ltms are known to us, though they take us a little further afield so we do not state them precisely. Both are due to Nicol. In the first paper (\citeyear{nicol2}) he constructs a \ltm which has the Bernoulli property, despite having local invariant manifolds and positive Lyapunov exponents only on a null set. In the second paper (\citeyear{nicol1}) he considers a Bernoulli \ltm of infinite entropy having smooth local invariant manifolds and positive Lyapunov exponents a.e.\ with some discontinuities. He shows that the map is stochastically stable.

\subsection{Linked-twist maps in the plane}\label{LitReview.Ltm.Plane}

Let $\Theta=\Gamma^k\circ\Phi^j:A\to A$, with $A\subset\mathbb{R}^2$, be a planar linked-twist map as defined in Section~\ref{Intro.Theorems.Plane}. The study of these maps was motivated by a number of authors. \citet{thurston} encountered linked-twist maps such as these in his study of diffeomorphisms of surfaces and \citet{braun} showed that similar maps arise as an approximate model of the global flow for the St\"{o}rmer problem. \citet{bowen} showed that certain linked-twist maps on such a manifold have positive topological entropy, and asked whether they possessed any ergodic properties. The following results describe what is known. As before $jk>0$ is the co-twisting case.

\begin{thm}[\citet{d2}]
For non-zero $j,k$ there is an invariant zero-measure Cantor set on which the map $\Theta$ is topologically conjugate to a subshift of finite type.
\end{thm}

Devaney's paper motivates our construction of a similar invariant set for a toral linked-twist map.

\begin{thm}[\citet{woj}]
If $jk>0$ and the twists are sufficiently strong, then $\Theta$ has an ergodic partition. Alternatively if $jk<0$ and a different (stronger) twist condition is satisfied, then also $\Theta$ has an ergodic partition.
\end{thm}

We give further details in Section~\ref{PLTM.Wojtkowski}. As we have mentioned, Wojtkowski's work is the main source of inspiration for our proof. We mention also some unpublished notes of \citet{p_preprint}; the planar linked-twist maps are amongst those that he discusses. In \citet{p2} the author considers again this large class of maps and shows that under certain conditions periodic saddles and homoclinic points are dense, and that the maps are topologically transitive.

\subsection{Other maps with strong ergodic properties}\label{LitReview.Hyp.Other}

Finally we describe some other systems for which strong ergodic properties have been established. The list, although not exhaustive, contains the majority of examples of which we are aware. As such it demonstrates the relative scarcity of such results and serves to underline the significance of our constructions.

The main classes of examples and specific examples of which we are aware are: geodesic flows on manifolds with negative curvature (\citet{anosov_sinai}, \citet{bb}, \citet{burns}); gases of hard spheres (\citet{kss}); symplectic Anosov and pseudo-Anosov diffeomorphisms (\citet{anosov_sinai}, \citet{gerber}, \citet{mackay06}); geodesic flows on surfaces with special metrics and potentials (\citet{donnay}, \citet{bg}, \citet{knauf}); systems like Wojtkowski's (\citeyear{woj2}); certain rational maps of the sphere (\citet{bk}); and the Belykh map (\citet{sataev}).

Moreover we mention the important work of \citet{katok} who showed that Bernoulli diffeomorphisms may occur on any surface. Beginning with a hyperbolic toral automorphism similar to that which we have encountered (and which is uniformly hyperbolic), Katok `slows down' trajectories in a neighbourhood of the origin. The manner in which this slowing down is accomplished is a little sophisticated and we do not intend to provide the details here; an excellent account can be found in \citet{bap}.

One of the consequences is that a Bernoulli map, \emph{derived from} an Anosov diffeomorphism, can be shown to exist on $\mathbb{S}^2$. This is particularly interesting given the well-known results of \citet{hirsch} and \citet{shiraiwa} that $\mathbb{S}^2$ cannot support an Anosov diffeomorphism.

%% file: chapter1/section4.tex
\section{Applications of the linked-twist map theory}\label{LitReview.App}
\setcounter{equation}{0}

In recent years the study of linked-twist maps has taken on a new significance owing to developments in our understanding of the mechanisms underlying good mixing of fluids. \citet{ottbook} has shown that the single most important feature to incorporate in the design of any fluid mixing device is the `crossing of streamlines', by which we mean that flow occurs periodically in two transversal directions. That linked-twist maps provide a suitable paradigm for this design process was highlighted in \citet{ow_science} and has been discussed at much greater length in \citet{ow2} and \citet{sturman}.

\subsection{DNA microarrays}\label{LitReview.App.Torus}

One important example of physical systems that may be analysed within the \ltm framework are certain models of DNA microarrays. In this section we summarise some ideas contained in \citet{hsw}; see that paper and the references therein for further details. DNA microarrays have been used widely in biochemical analysis for a number of years. Amongst their uses are gene discovery and mapping, gene regulation studies, disease diagnosis and drug discovery and toxicology.

A DNA microarray consists of DNA strands (`probes') fixed to a surface such as glass or silicon. This array is placed in a hybridization chamber containing a solution of DNA or mRNA (the `target'). Hybridization occurs when the target strand combines with a complementary probe strand, as governed by base-pairing rules.

Hybridization is most efficient when each target strand can move throughout the solution and encounter every probe. Two processes lead to this: diffusion and advection. The former cannot be relied upon to produce the desired result in reasonable time because the typical situation involves low Reynolds number and thus no turbulence. Advection in such devices has consequently been the focus of a great deal of research into how one might induce good mixing.

It is now well established that chaotic advection provides a source of efficient mixing in many fluid problems, in particular those on the `microfluidic' scale of the DNA microarrays. Two designs for such devices, both relying upon cyclic removal and reinjection of fluid into the mixing chamber, are detailed in \citet{mspflsh} and \citet{rpbcsmcc}. Typically two different source-sink pairs are used.

Two factors which have a great impact on the efficiency of such mixers are the locations at which fluid should be removed and reinjected, and the time for which such a source-sink pair should be active. Here linked-twist map theory can help to inform the design. The motion of fluid in such a device can bear striking resemblances to the motion of a point in the domain of a toral linked-twist map. Analysing these mixing devices in this manner has lead to the proposal of new mixing protocols.

\subsection{Channel-type micromixers}\label{LitReview.App.Plane}

There are many areas of applications where the \ltms which most naturally act as models are defined on surfaces other than the two-torus. One such example are channel-type micromixers. In this section we briefly summarise some ideas presented in \citet{sdamsw}; see that paper and the references therein for further details. 

Mixing of the fluid flowing through a microchannel is highly desirable in a number of situations, including the homogenisation of solutions of reagents in chemical reactions and in the control of dispersion of material along the direction of Poiseuille flows. \citet{sdamsw} argue that, at low Reynolds number and using a `simple' channel (i.e.\ one that is straight and has smooth walls) any mixing is a consequence of diffusion only. Moreover they conclude that the rate at which this happens, even in a microchannel, is slow compared with convection along the channel. To reduce the length of channel required for mixing to occur one needs to introduce transversal components of the flow which stretch and fold volumes of fluid over the cross section of the channel, thus reducing the distance over which diffusion must act.

Such transversal flows may be generated by placing ridges on the floor of the channel, at an oblique angle to the flow. The ridges present anisotropic resistance to viscous flows resulting in transversal flow, which then circulates back across the top of the channel. Consequently flow along the channel becomes helical. The helical motion of the flow is a motion that can be approximated by certain planar linked-twist maps.

\subsection{Other examples}\label{LitReview.App.Other}

The existence of the above examples illustrates that the linked-twist map framework can be a useful tool in the development of models for certain mixing devices. In the case of planar linked-twist maps there are numerous other examples we could have mentioned, the blinking vortex flow of \citet{aref1} being a prime case in point. Here a pair of point vortices in an unbounded inviscid fluid are alternately switched on. Further work on this system was conducted by \citet{kro}.

The work has applications to the study of tidal flow close to a headland jutting out into the sea. For details see \citet{sg}, \citet{sb} and \citet{samunpub}. In the latter reference a kinematic model to study mixing and transport by eddies is developed. For details on how \ltms may be used as a paradigm for studying such systems we direct the reader to \citet{w3}.

Yet more examples are given by the electroosmotic stirrers of \citet{qb}; the cavity flows introduced by \citet{cro} and developed further by \citet{lo} and \citet{jmo}; and the egg-beater flows of \citet{fo2}.

\Ltms on the two-sphere do not lend themselves to applications as readily as their counterparts in the plane, but we are able to extract an example from the field of quantum ergodicity. \citet{marklof_o'keefe} have shown that the quantum eigenstates of linked-twist maps defined on the two-torus are equidistributed. This result uses the fact that the corresponding classical linked-twist maps are ergodic. \citet{o'keefe} demonstrates that an analogous result can be shown to hold on the two-sphere, if one is able to show that the corresponding classical maps are ergodic. We mention another possible application for which linked-twist maps on the sphere might be a useful analytical tool in Chapter~\ref{concl}.

Finally we mention an application to granular mixing. There are numerous situations in the pharmaceutical, food, chemical, ceramic, metallurgical and construction industries where an understanding of the behaviour of granular media is crucial. However the literature dedicated to the mixing of these materials is not nearly as developed as its counterpart for fluids. Recent work has shown that this is yet another example of a physical process which may be modelled using linked-twist maps. For further details see \citet{sturman_granular} and the references therein.

%% file: chapter2/chapter2.tex
\chapter{A horseshoe in a toral linked-twist map}
\setcounter{equation}{0}
\setcounter{figure}{0}
\label{chapter2}

This chapter is motivated by work of \citet{d2} in which the author establishes the existence of a topological horseshoe in a planar linked-twist map. We construct, using similar ideas, the counterpart for a toral linked-twist map. To our knowledge this is not to be found in the literature.

Let $j,k$ be positive integers such that $N=(j-1)(k-1)\geqs 2$ (i.e.\ each of $j,k$ is at least 2 and one is at least 3), and let $G^k\circ F^j:R\to R$ be a toral linked-twist map, as in Section~\ref{Intro.Theorems.Torus}. We show that there exists a zero-measure, compact, invariant set within $R\subset\mathbb{T}^2$ (a `horseshoe') on which the dynamics are topologically conjugate to a full shift on $N$ symbols. (By contrast, Devaney's construction yields a conjugacy with a \emph{subshift} instead; we discuss this further in Chapter~\ref{concl}.)

In Section~\ref{Horseshoe.ConMos} we develop some notation with which to state a theorem due to \citet{moser}, which provides sufficient critera for the existence of the horseshoe. In Section~\ref{Horseshoe.Strips} we show that toral linked-twist maps as above satisfy these criteria; this entails a detailed geometrical construction. Finally in Section~\ref{Horseshoe.UH} we provide some extra analysis to show that the map restricted to the horseshoe is uniformly hyperbolic.

Before we begin we remark on our notation. It will be natural in this chapter to reserve the letter $H$ for \emph{horizontal strips}, which we define below. To avoid confusion, throughout this chapter, the linked-twist map will always be denoted by $G^k\circ F^j$ (and \emph{not} by $H_{j,k}$), whereas $H,H_i$ etc.\ denote horizontal strips.

\input chapter2/section1

\input chapter2/section2

\input chapter2/section3

%% file: chapter2/section1.tex
\section{The Conley-Moser conditions}\label{Horseshoe.ConMos}
\setcounter{equation}{0}

Here we describe sufficient criteria for a two-dimensional invertible map to possess an invariant Cantor set on which the aforementioned conjugacy exists. These are commonly known as the \emph{Conley-Moser conditions}, having first been introduced in \citet{moser}. A lucid and comprehensive account of this material may be found in \citet{wig1}.

\subsection{Horizontal and vertical curves and strips}

We begin with some definitions. Recall that a real-valued function $f$ defined on connected domain $D\subset\mathbb{R}$ is \emph{Lipschitz continuous} if and only if there exists a constant $c>0$ and for every pair $a,b\in D$ we have
\beq*
|f(a)-f(b)|\leqs c|a-b|.
\eeq*
We will say that such an $f$ is $c$-\emph{Lipschitz}. We use this notation to define \emph{curves} in $S=P\cap Q\subset R$. Recall that $S=[x_0,x_1]\times[y_0,y_1]$.

\begin{defn}[$m_h$-horizontal and $m_v$-vertical curves]
An $m_h$-horizontal curve is the graph of an $m_h$-Lipschitz function $h:[x_0,x_1]\to[y_0,y_1]$. An $m_v$-vertical curve is the graph of an $m_v$-Lipschitz function $v:[y_0,y_1]\to[x_0,x_1]$.
\end{defn}

We use such curves to form the boundaries of \emph{strips} as follows.

\begin{defn}[$m_h$-horizontal and $m_v$-vertical strips]
Given two non-intersecting $m_h$-horizontal curves of functions $h_1$ and $h_2$, with $h_1(x)<h_2(x)$ for each $x\in[x_0,x_1]$, an $m_h$-horizontal strip is the set
\beq*
H=\left\{(x,y)\in S: x\in[x_0,x_1], y\in[h_1(x),h_2(x)]\right\}.
\eeq*
The $m_h$-horizontal curves are then referred to as the horizontal boundaries of $H$. The vertical boundaries of $H$ are contained within the lines $x=x_0$ and $x=x_1$.

Similarly given two non-intersecting $m_v$-vertical curves of functions $v_1$ and $v_2$, with $v_1(y)<v_2(y)$ for each $y\in[y_0,y_1]$, an $m_v$-vertical strip is the set
\beq*
V=\left\{(x,y)\in S: y\in[y_0,y_1], x\in[v_1(y),v_2(y)]\right\}.
\eeq*
The $m_v$-vertical curves are then referred to as the vertical boundaries of $V$. The horizontal boundaries of $V$ are contained within the lines $y=y_0$ and $y=y_1$.
\end{defn}

\begin{figure}[htp]
\centering
\includegraphics[totalheight=0.35\textheight]{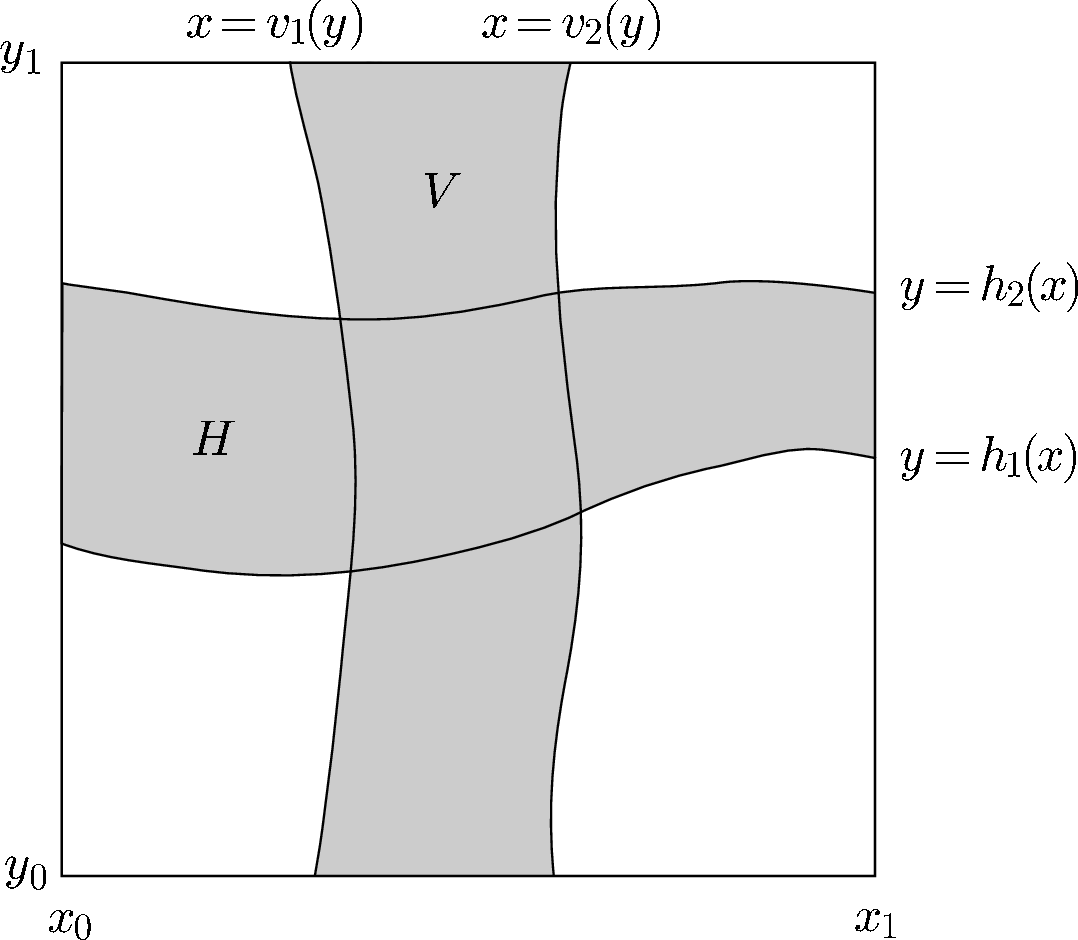}
\caption[Horizontal and vertical curves]{The region $S\subset R$ showing horizontal curves of $y=h_1(x)$ and $y=h_2(x)$, vertical curves of $x=v_1(y)$ and $x=v_2(y)$, a horizontal strip $H$ and a vertical strip $V$.}
\label{fig:curves_and_strips}
\end{figure}

We illustrate some curves and strips in Figure~\ref{fig:curves_and_strips}. Define the \emph{width} of a strip as follows:

\begin{defn}[Width of $m_h$-horizontal and $m_v$-vertical strips]
Let $H$ be an $m_h$-horizontal strip as above. Its width is given by
\beq*
d(H) = \max_{x\in[x_0,x_1]}|h_2(x)-h_1(x)|.
\eeq*
Similarly let $V$ be an $m_v$-vertical strip as above. Its width is given by
\beq*
d(V) = \max_{y\in[y_0,y_1]}|v_2(y)-v_1(y)|.
\eeq*
\end{defn}

\subsection{The Conley-Moser conditions}

Let $\psi:S\to R$ be a map and let $I=\{1,2,...,N\}$ be an index set for some $N\in\mathbb{N}$. Let $\{H_i\}_{i\in I}$ be a set of disjoint $m_h$-horizontal strips and let $\{V_i\}_{i\in I}$ be a set of disjoint $m_v$-vertical strips. The \emph{Conley-Moser conditions} on $\psi$ are as follows:

\begin{cond}
\label{cond:CM1}
$0\leqs m_vm_h<1$.
\end{cond}

\begin{cond}
\label{cond:CM2}
$\psi$ maps $H_i$ homeomorphically onto $V_i$ (i.e.\ $\psi(H_i)=V_i$) for each $i\in I$. Moreover, the horizontal boundaries of $H_i$ map to the horizontal boundaries of $V_i$ and the vertical boundaries of $H_i$ map to the vertical boundaries of $V_i$.
\end{cond}

\begin{cond}
\label{cond:CM3}
Suppose $H\subset\bigcup_{i\in I}H_i$ is an $m_h$-horizontal strip and let
\beq*
\tilde{H}_i=\psi^{-1}(H)\cap H_i.
\eeq*
Then $\tilde{H}_i$ is an $m_h$-horizontal strip for each $i\in I$ and $d(\tilde{H}_i)\leqslant n_hd(H)$ for some $0<n_h<1$.

Similarly suppose $V\subset\bigcup_{i\in I}V_i$ is an $m_v$-vertical strip and let
\beq*
\tilde{V}_i=\psi(V)\cap V_i.
\eeq*
Then $\tilde{V}_i$ is an $m_v$-vertical strip for each $i\in I$ and $d(\tilde{V}_i)\leqslant n_vd(V)$ for some $0<n_v<1$.
\end{cond}

The result we will use is the following.

\begin{thm}[\citet{moser}]
\label{thm:ConleyMoser}
Suppose $\psi:S\to R$ satisfies Conditions~\ref{cond:CM1}, \ref{cond:CM2} and~\ref{cond:CM3}. Then $\psi$ has an invariant Cantor set $\Lambda$ on which it is topologically conjugate to a full shift on $N$ symbols, i.e.\ the following diagram commutes:
\smallskip
\begin{center}
\begin{pspicture}(0,.5)(3,3)
\psset{nodesep=5pt} 
\rput(.5,2.5){\rnode{A}{$\Lambda$}} 
\rput(2.5,2.5){\rnode{B}{$\Lambda$}} 
\rput(.5,1){\rnode{C}{$\Sigma^N$}} 
\rput(2.5,1){\rnode{D}{$\Sigma^N$}}
\ncline{->}{A}{B}\Aput{$\psi$}
\ncline{->}{B}{D}\Aput{$\phi$} 
\ncline{->}{A}{C}\Bput{$\phi$} 
\ncline{->}{C}{D}\Bput{$\sigma$}
\end{pspicture} 
\end{center}
\smallskip
where $\phi:\Lambda\to\Sigma^N$ is a homeomorphism and $\sigma:\Sigma^N\to\Sigma^N$ is the shift map on the space of symbol sequences, defined in section~\ref{LitReview.ErgTh.Bern}.
\end{thm}

We remark that the restriction of $\psi$ to $\Lambda$ has the Bernoulli property.

%% file: chapter2/section2.tex
\section{Construction of the strips}\label{Horseshoe.Strips}
\setcounter{equation}{0}

We now construct strips satisfying the Conley-Moser conditions of the previous section. Our first task will be to define a certain quadrilateral $M\subset S$. It will transpire that the images and pre-images of $M$ (with respect to the linked-twist map $G^k\circ F^j$) have certain convenient properties; in particular they contain the horizontal and vertical strips we require.

\subsection{Construction of the quadrilateral $M\subset S$}

\begin{figure}[htp]
\centering
(a)\includegraphics[totalheight=0.17\textheight]{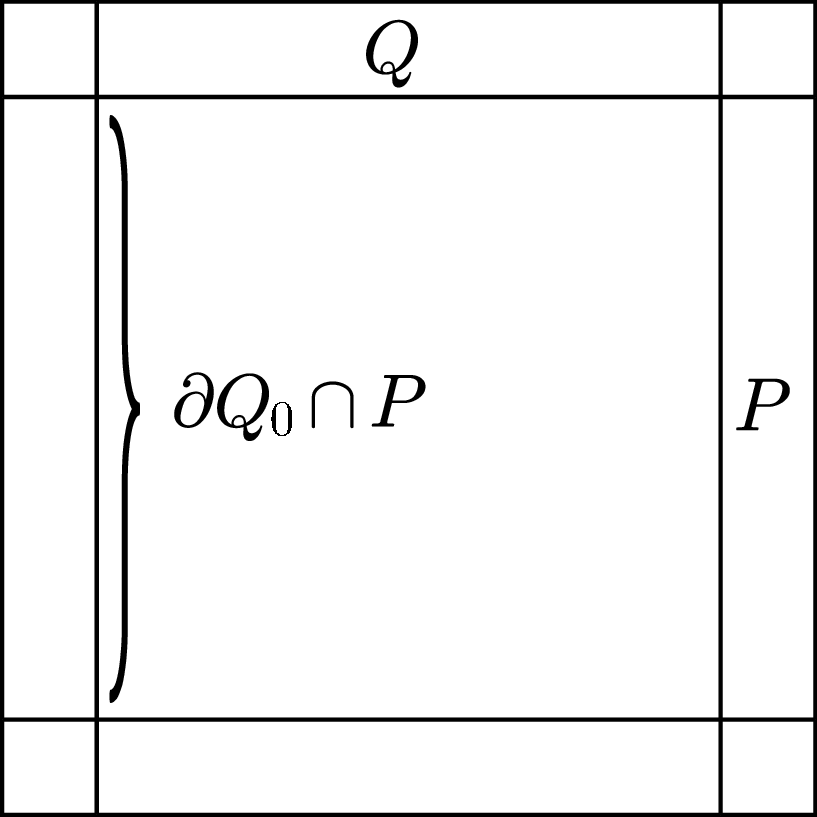}\quad
(b)\includegraphics[totalheight=0.17\textheight]{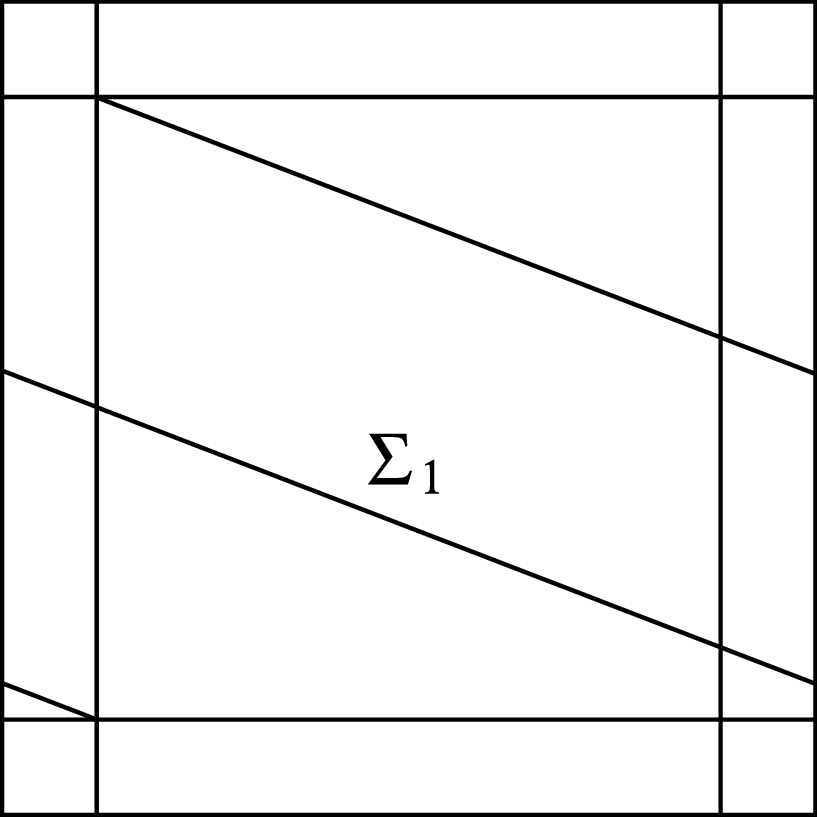}\quad
(c)\includegraphics[totalheight=0.17\textheight]{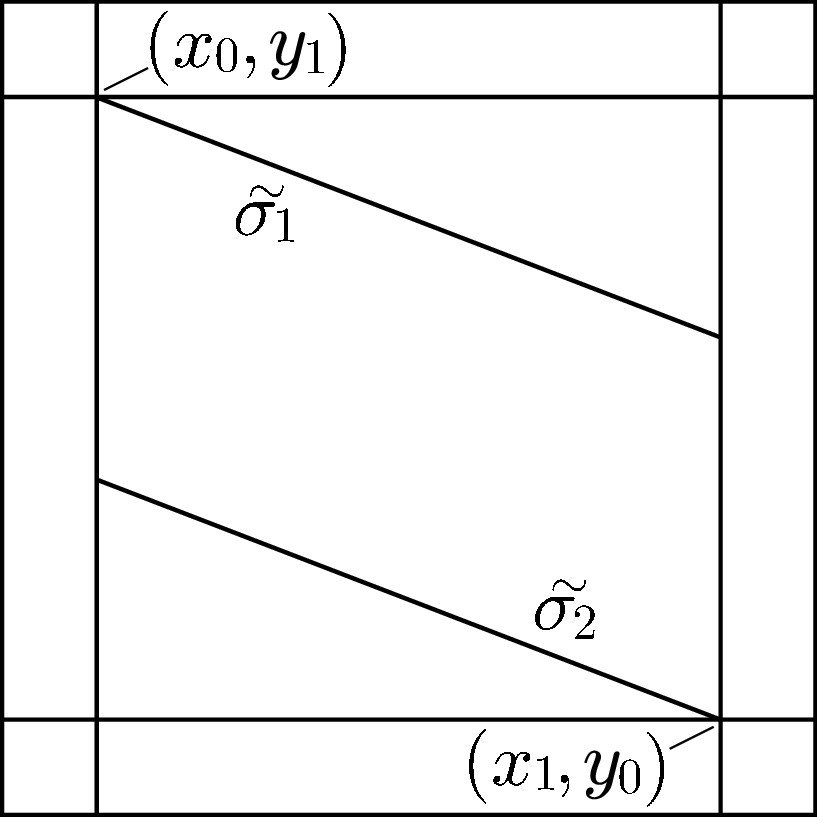}
\caption[Construction of $M$: Part I]{Construction of the horizontal boundaries of $M$.}
\label{fig:M1}
\end{figure}

\begin{figure}[htp]
\centering
(a)\includegraphics[totalheight=0.17\textheight]{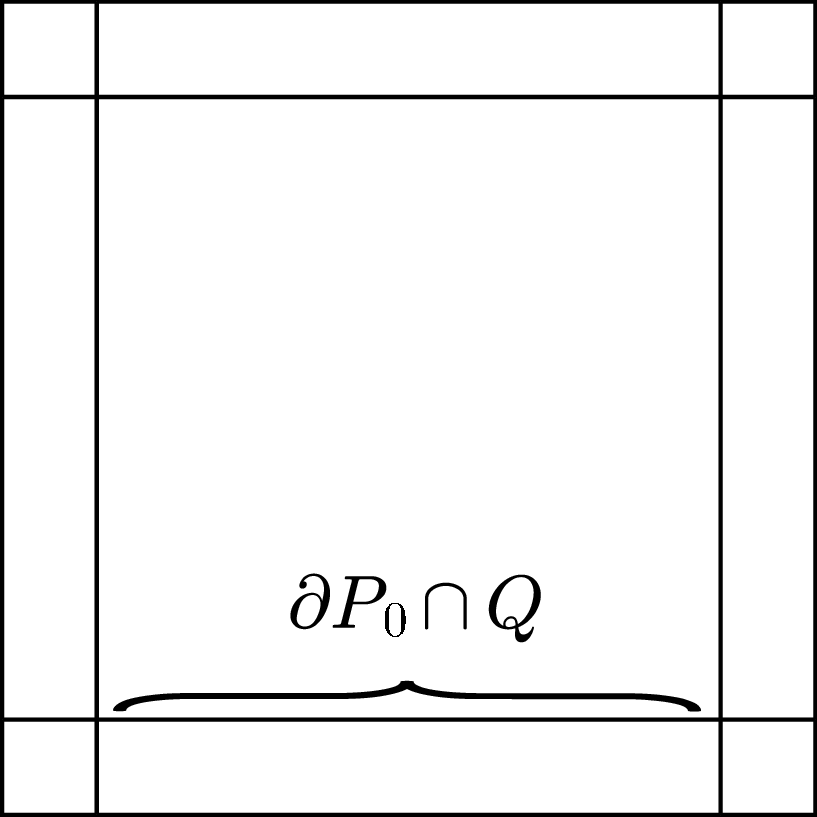}\quad
(b)\includegraphics[totalheight=0.17\textheight]{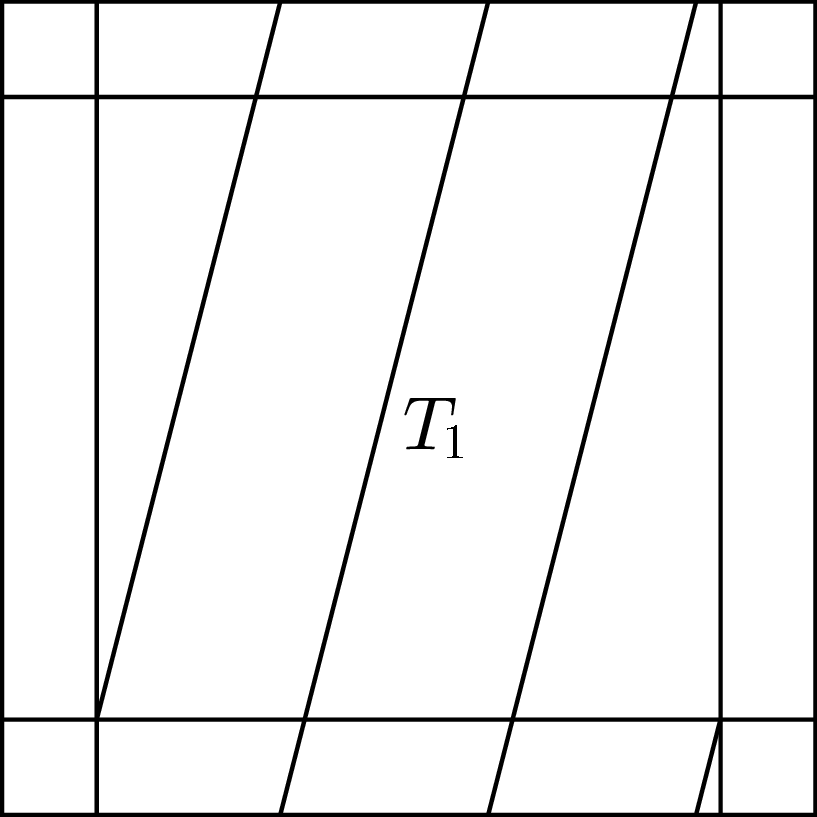}\quad
(c)\includegraphics[totalheight=0.17\textheight]{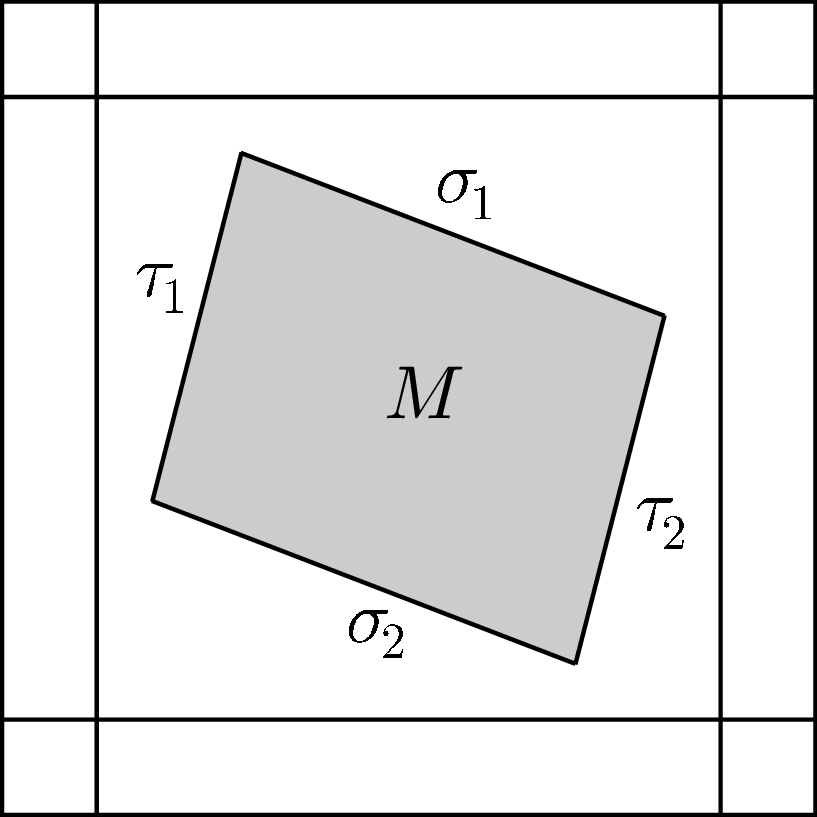}
\caption[Construction of $M$: Part II]{Construction of the vertical boundaries of $M$.}
\label{fig:M2}
\end{figure}

Figures~\ref{fig:M1} and~\ref{fig:M2} illustrate the construction of $M$, which we now describe. Recall our notation for the manifold $R\subset\mathbb{T}^2$, established in Section~\ref{Intro.Theorems.Torus}; in particular we have a `horizontal' annulus $P$ with boundaries $\pd P_0$ (on which $y=y_0$) and $\pd P_1$ (on which $y=y_1$), and a `vertical' annulus $Q$ with boundaries $\pd Q_0$ (on which $x=x_0$) and $\pd Q_1$ (on which $x=x_1$).

Consider the portion of the boundary of $Q$ given by $\partial Q_0\cap P$. This is shown in Fig~\ref{fig:M1}(a). Let $\Sigma_1=F^{-j}(\partial Q_0\cap P)$, as shown in \ref{fig:M1}(b). For illustrative purposes we have taken $j=2$. We observe that $\Sigma_1\cap S$ consists of $j$ disjoint pieces, each of which stretches across $S$ from $\partial Q_0$ to $\partial Q_1$.

One of these pieces has as an end-point the point $(x_0,y_1)$. Let $\tilde{\sigma_1}\subset\Sigma_1\cap S$ denote this piece, as shown in \ref{fig:M1}(c). Similarly we define $\tilde{\sigma_2}\subset\Sigma_2\cap S$, which is derived in a similar manner from $\Sigma_2=F^{-j}(\partial Q_1\cap P)$ and shown in the same figure; it has as an end-point the point $(x_1,y_0)$.

Analogously, part~(a) of Figure~\ref{fig:M2} shows $\partial P_0\cap Q$ and \ref{fig:M2}(b) shows its image with respect to the map $G^k$. We have illustrated this using $k=3$. Define $T_1=G^k(\partial P_0\cap Q)$, then $T_1\cap S$ consists of $k$ disjoint pieces, each of which stretches across $S$ from bottom to top.

Let $\tilde{\tau_1}\subset T_1\cap S$ be that piece which has as an end-point the point $(x_0,y_0)$. Similarly define $\tilde{\tau_2}$ to be that piece of $T_2\cap S$ which has as an end-point the point $(x_1,y_1)$; here $T_2=G^k(\partial P_1\cap Q)$.

Part~(c) of Figure~\ref{fig:M2} shows the quadrilateral $M\subset S$, which is bounded by the four lines $\tilde{\sigma_1},\tilde{\sigma_2},\tilde{\tau_1}$ and $\tilde{\tau_2}$. Notice that the boundary consists of two pairs of parallel lines and so $M$ is a parallelogram. Finally, we denote by $\sigma_1\subset\tilde{\sigma_1}$ that part of $\tilde{\sigma_1}$ which is a part
of the boundary of $M$, and similarly $\sigma_2\subset\tilde{\sigma_2}$, $\tau_1\subset\tilde{\tau_1}$ and $\tau_2\subset\tilde{\tau_2}$.

\subsection{Intersection of $M$ with its images and pre-images}

In this section we prove results concerning the intersection of $M$ with its image and pre-image (with respect to the linked-twist map $G^k\circ F^j$) respectively. The resulting sets \emph{resemble} but technically are \emph{not} collections of $m_v$-vertical and $m_h$-horizontal strips, because they do not stretch completely across $S$. It turns out that they \emph{are} the intersection of such
strips with $M$. We state the required result as a proposition.

We recall that the inverse of $G^k\circ F^j$ is given by $F^{-j}\circ G^{-k}$.

\begin{prop}
\label{prop:properties_of_M}
$M$ has the following properties:
\begin{enumerate}
\item $M\cap (G^k\circ F^j(M))$ is the intersection of $M$ with $(j-1)(k-1)$ disjoint $m_v$-vertical strips. These strips intersect only the boundaries $\sigma_1$ and $\sigma_2$ of $M$ (i.e.\ they do not intersect the boundaries $\tau_1$ or $\tau_2$).
\item Similarly, $(F^{-j}\circ G^{-k}(M))\cap M$ is the intersection of $M$ with $(j-1)(k-1)$ disjoint $m_h$-horizontal strips. These strips intersect only the boundaries $\tau_1$ and $\tau_2$ of $M$ (i.e.\ they do not intersect the boundaries $\sigma_1$ or $\sigma_2$).
\item The above holds with $0<m_vm_h<1$.
\end{enumerate}
\end{prop}

\begin{figure}[htp]
\centering
(a)\includegraphics[totalheight=0.25\textheight]{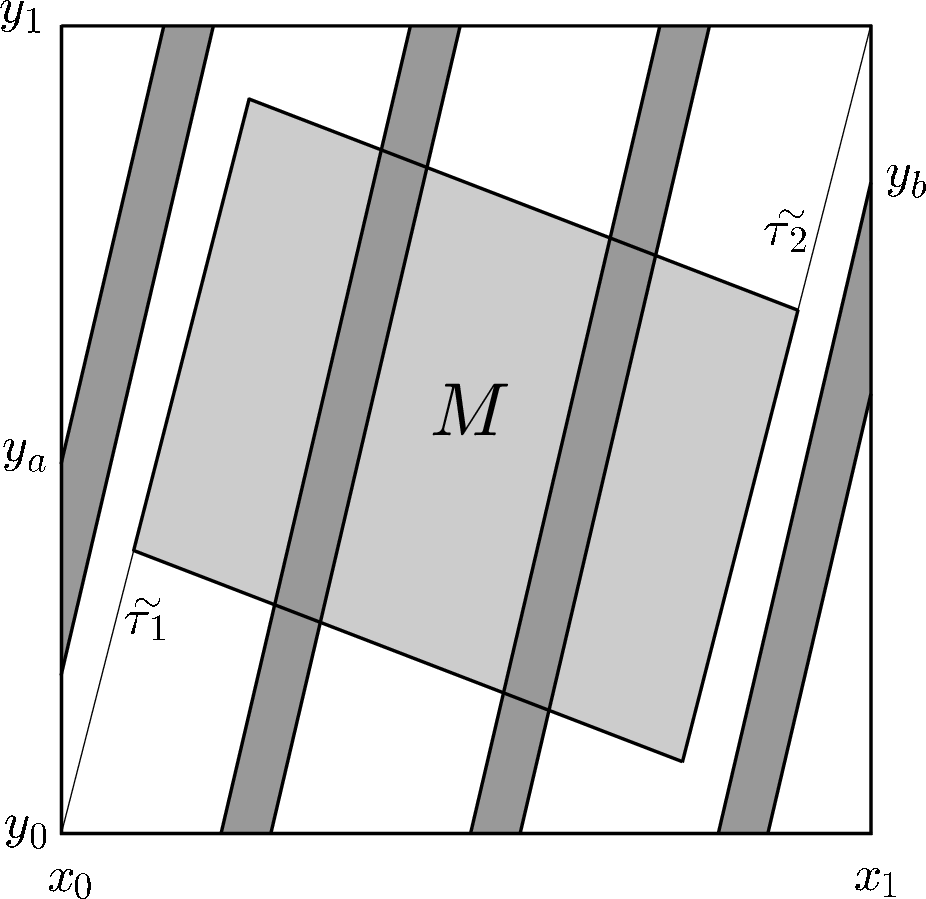}\quad
(b)\includegraphics[totalheight=0.257\textheight]{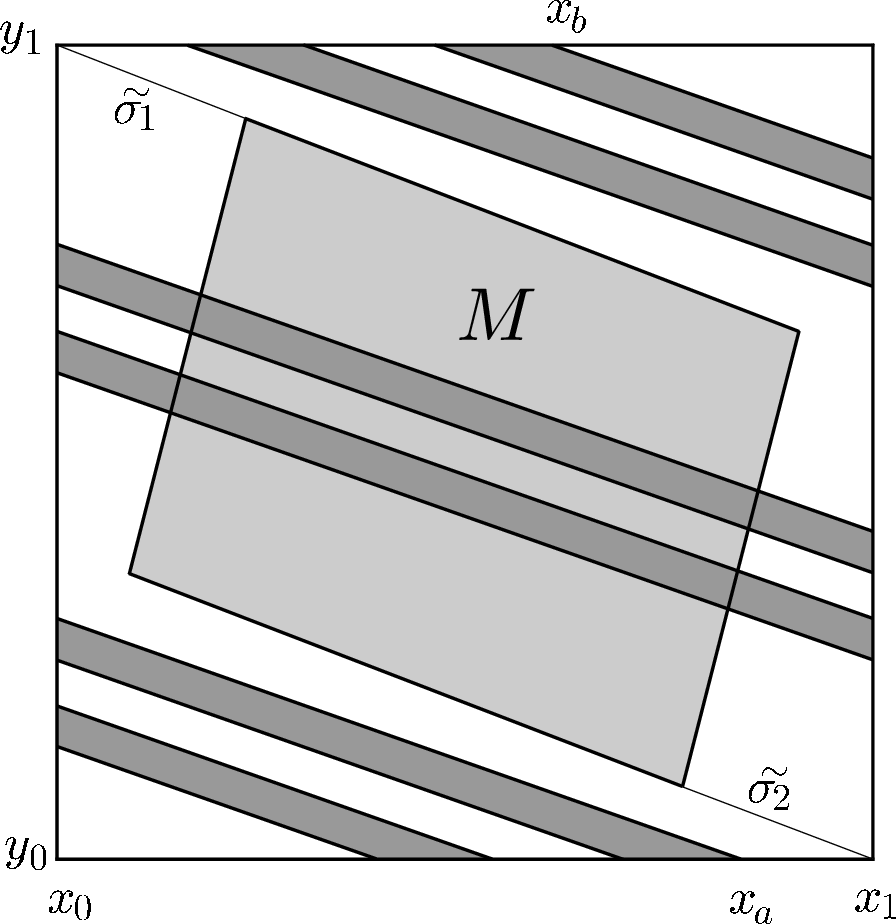}
\caption{Intersection of $M$ with the horizontal and vertical strips.}
\label{fig:strips3}
\end{figure}

\begin{proof}
We prove each part in turn.

\begin{figure}[htp]
\centering
(a)\includegraphics[totalheight=0.2\textheight]{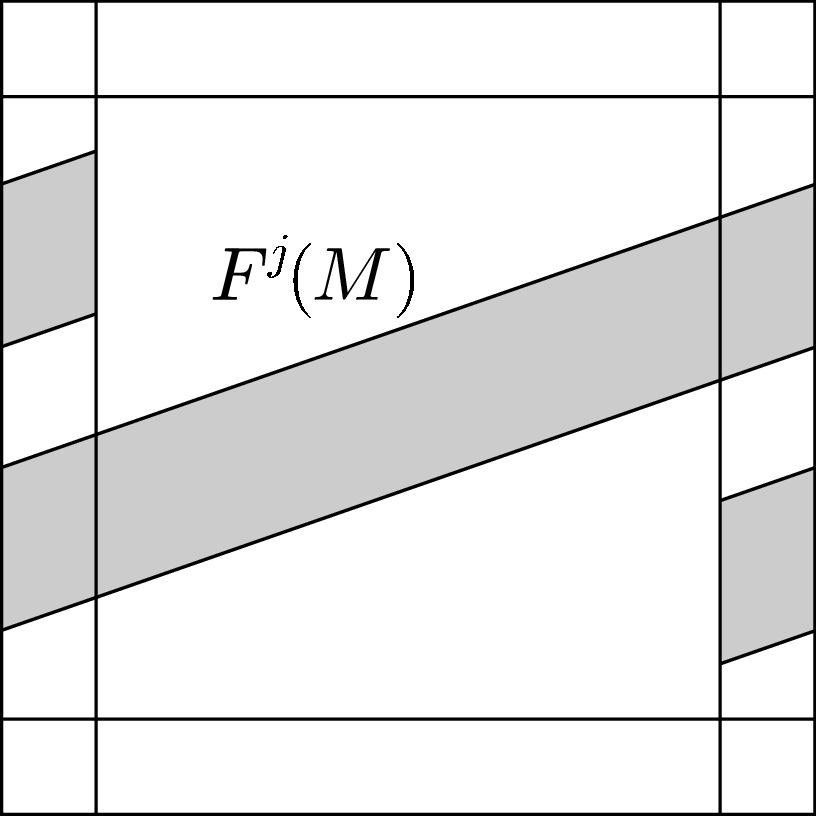}\quad
(b)\includegraphics[totalheight=0.2\textheight]{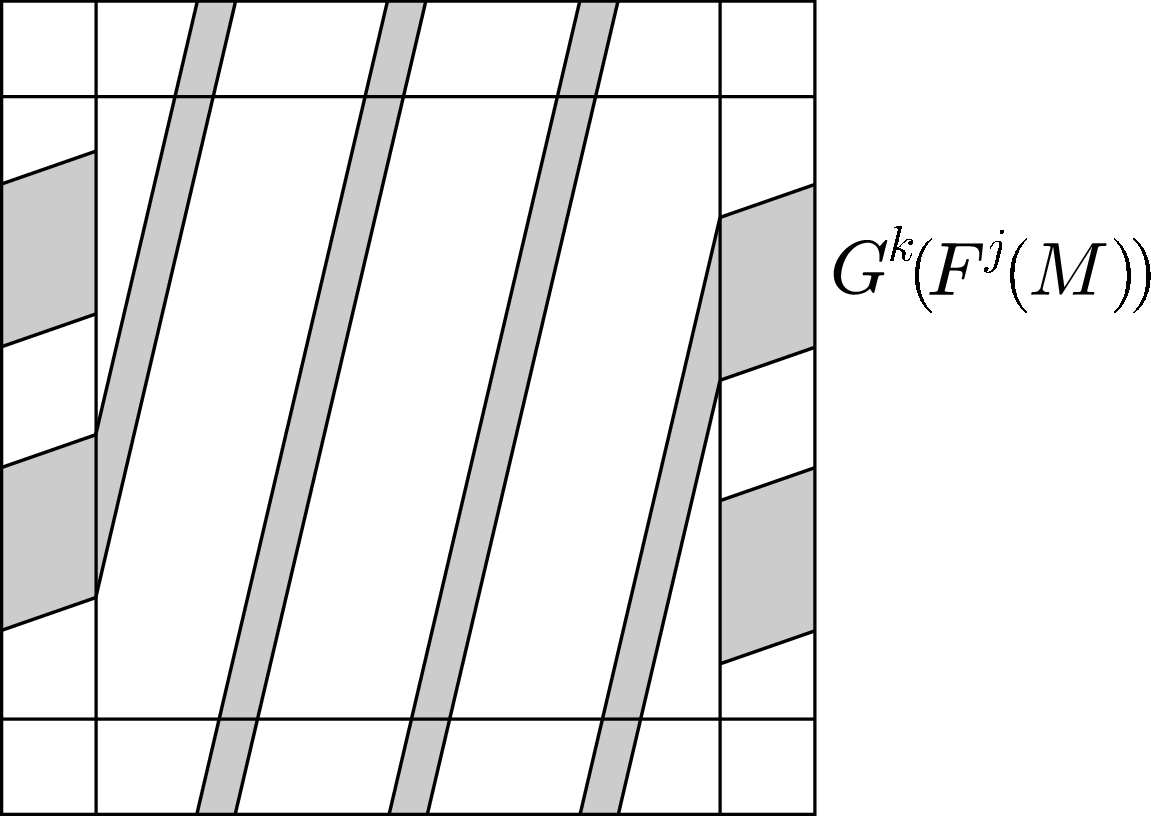}
\caption{Construction of vertical strips crossing $S$.}
\label{fig:strips1}
\end{figure}

\begin{enumerate}
\item Consider $F^j(M)$, as illustrated in Figure~\ref{fig:strips1}(a) with $j=2$. By the construction of $M$ this strip crosses $S$ horizontally $(j-1)$ times and, because the linked-twist map is a homeomorphism, no two of these crossings intersect each other. Part~(b) of the figure shows $G^k\circ F^j(M)$, illustrated with $k=3$. Observe that $(G^k\circ F^j(M))\cap S$ consists of $(j-1)(k+1)=4$ disjoint pieces.

The boundaries of $M$ are straight lines and $f$ and $g$ are affine, so the boundaries of these pieces are straight lines. By definition, the $(j-1)(k-1)=2$ such pieces that cross $S$ completely (i.e.\ intersect both $\partial P_0$ and $\partial P_1$) are disjoint $m_v$-vertical strips.

In Figure~\ref{fig:strips3}(a) we show the same situation in closer detail and with the original set $M$ overlaid. To
complete the proof of the first part of the proposition, it suffices to show that none of the $(j-1)(k+1)$ pieces of $(G^k\circ F^j(M))\cap S$ intersect $(\tilde{\tau_1}\cup\tilde{\tau_2})\subset G^k\circ F^j(\partial P)$. Assume for a contradiction that this is \emph{not} the case. Consequently we can find some $(x,y)\in R$ for which $(x,y)\in (G^k\circ F^j(M))\cap (G^k\circ F^j(\partial P))$. Then $(F^{-j}\circ G^{-k}(x,y))\in M\cap\partial P=\emptyset$, an obvious contradiction.

\begin{figure}[htp]
\centering
(a)\includegraphics[totalheight=0.2\textheight]{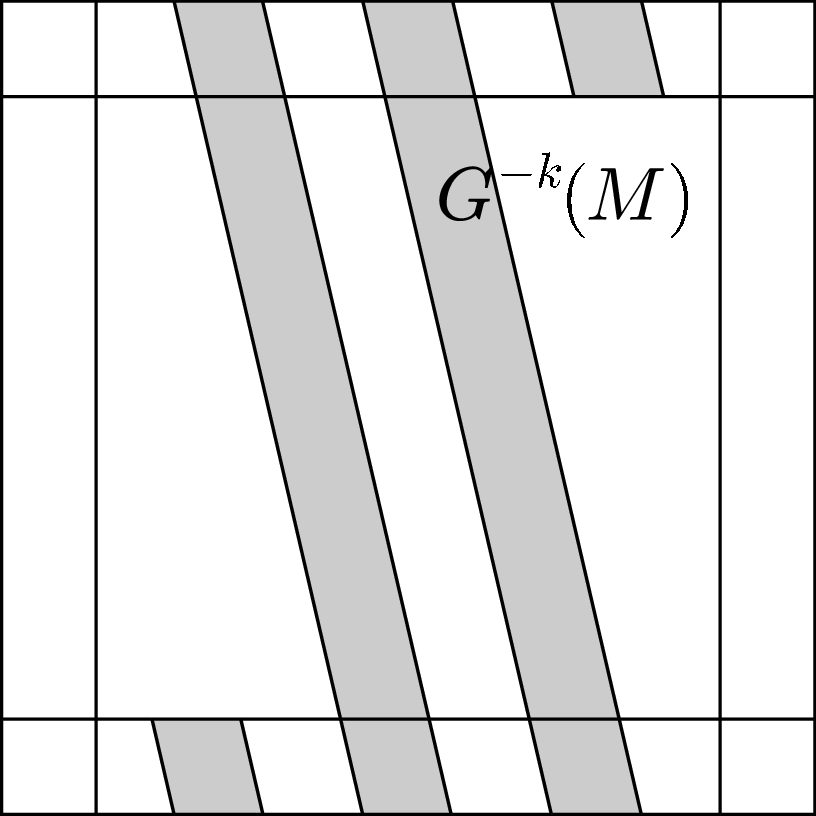}\quad
(b)\includegraphics[totalheight=0.2\textheight]{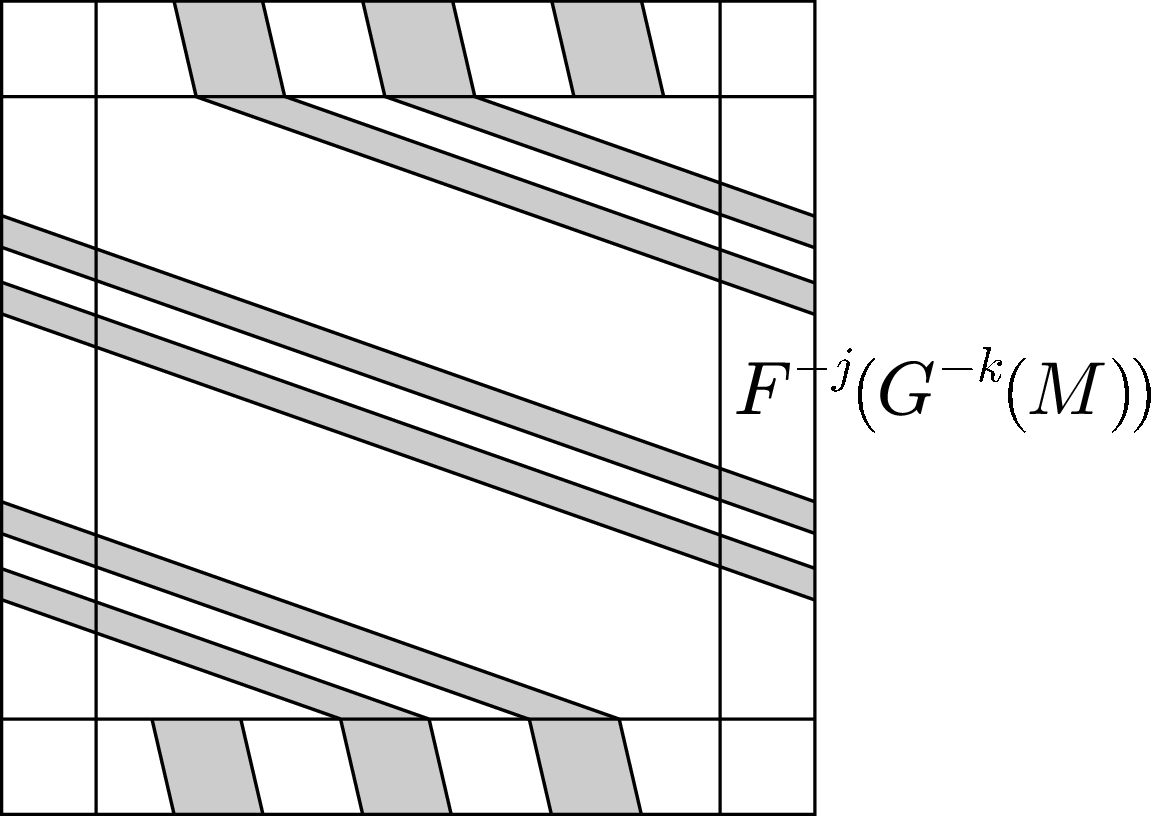}
\caption{Construction of horizontal strips crossing $S$.}
\label{fig:strips2}
\end{figure}

\item The proof of the second part is similar. $G^{-k}(M)$ crosses $S$ vertically $(k-1)$ times and these strips are disjoint; see Figure~\ref{fig:strips2}(a). $(F^{-j}\circ G^{-k}(M))\cap S$ consists of $(j+1)(k-1)=6$ disjoint pieces each bounded by straight lines as shown in part~(b) of the figure. Those pieces which cross $S$ completely are disjoint $m_h$-horizontal strips.

Now consider Figure~\ref{fig:strips3}(b). In order to show the second part of the lemma it suffices to show that none of these $(j+1)(k-1)$ pieces intersect $\sigma_1$ or $\sigma_2$. This follows (as in part~1 of the proof) from the fact that any such point would have an image in $M\cap\partial Q=\emptyset$, a clear contradiction.

\item Last of all we consider the size of $m_vm_h$. It is clear from Figure~\ref{fig:strips3} part~(b) that each boundary of each $m_h$-horizontal strip is an $m_h$-horizontal curve, and moreover that it is a straight line of constant gradient
\beq*
\frac{y_1-y_0}{j+x_b-x_a}
\eeq*
with respect to $x$. Here $x_a$ and $x_b$ are as shown in the figure. It will be enough for present purposes to observe that $0<x_b-x_a<x_1-x_0<1$, thus this gradient (which is a suitable value for $m_h$) satisfies
\beq*
\frac{y_1-y_0}{j+1}<m_h<\frac{y_1-y_0}{j}.
\eeq*

Similarly from Figure~\ref{fig:strips3}(a) we conclude that $m_v$ may be taken to be the gradient
\beq*
\frac{x_1-x_0}{k+y_b-y_a}
\eeq*
with respect to $y$. The values $y_a$ and $y_b$ are as shown in the figure and satisfy $0<y_b-y_a<1$. Consequently
\beq*
\frac{x_1-x_0}{k+1}<m_v<\frac{x_1-x_0}{k}
\eeq*
and thus
\beq*
0<\frac{(x_1-x_0)(y_1-y_0)}{(j+1)(k+1)}<m_vm_h< \frac{(x_1-x_0)(y_1-y_0)}{jk}<1.
\eeq*
\end{enumerate}
\end{proof}

Recall our notation that $N=(j-1)(k-1)$ (and our assumption that $N\geqs 2$) and $I=\{1,2,...,N\}$ is an index set. Let $\{H_i\}_{i\in I}$ be the connected, pairwise-disjoint subsets of $(F^{-j}\circ G^{-k}(M))\cap M$. Proposition~\ref{prop:properties_of_M} says that for each $i\in I$ we have
\beql\label{eqn:H_intersect_M}
H_i=M\cap\bar{H_i},
\eeql
where $\bar{H_i}$ is an $m_h$-horizontal strip. Similarly, let $\{V_i\}_{i\in I}$ be the connected, pairwise-disjoint subsets of $M\cap (G^k\circ F^j(M))$. Then for each $i\in I$ we have
\beql\label{eqn:V_intersect_M}
V_i=M\cap\bar{V_i},
\eeql
where $\bar{V_i}$ is an $m_v$-vertical strip.

\subsection{Existence of the horseshoe}

In this section we show that the strips we have constructed satisfy the Conley-Moser conditions.

\begin{prop}
\label{prop:mapping_of_h_strips}
For each $i\in I$ the linked-twist map $G^k\circ F^j$ maps $H_i$ homeomorphically onto $V_i$ (up to permutation of the $V_i$). Moreover the horizontal boundaries of $H_i$ are mapped to the horizontal boundaries of $V_i$, and the vertical boundaries of $H_i$ are mapped to the vertical boundaries of $V_i$.
\end{prop}

\begin{proof}
$G^k\circ F^j$ is a homeomorphism of $R=P\cup Q$ and so it is immediate that any one of the $N$ disjoint pieces $H_i$ (a connected component of $(F^{-j}\circ G^{-k}(M))\cap M$) must map to one and only one of the $N$ disjoint pieces $V_i$ (the connected components of $M\cap (G^k\circ F^j(M))$), and must do so homeomorphically. Furthermore, elementary topology (see, for example, \citet{armstrong}) tells us that boundaries must map to boundaries.

The horizontal boundaries of the $H_i$ are shown in Figure~\ref{fig:strips2}(b) and their images under $F^j$, shown in part~(a) of that figure, are contained within the vertical boundaries of $G^{-k}(M)$. Applying the map $G^k$ in turn, these vertical boundaries become the horizontal boundaries of $M$, which contain the horizontal boundaries of $V_i$.

Similarly, the vertical boundaries of the $V_i$ are shown in Figure~\ref{fig:strips1}(b). The map $G^{-k}$ takes these into the horizontal boundaries of $F^j(M)$, as in part~(a) of that figure, and these in turn are mapped by $F^{-j}$ into the vertical boundaries of $M$. The vertical boundaries of the $H_i$ are contained within these vertical boundaries of $M$ and thus the result.
\end{proof}

At this point we discuss how the Conley-Moser conditions (as given) are not perfectly suited to the present task, how we get around this and how else we might have gotten around this.

The Conley-Moser conditions describe quite specifically how horizontal and vertical curves and strips are mapped onto each other. We have defined curves as stretching completely across $S$ and consequently the boundaries do not map in the required manner; to remedy this we consider the intersections of such curves with $M$, as in~(\ref{eqn:H_intersect_M}) and~(\ref{eqn:V_intersect_M}).

Defining curves and strips as we do allows us to determine our estimates on widths using the orthogonal $(x,y)$ coordinates on $\mathbb{T}^2$ at the expense of then having to intersect these curves and strips with $M$ in order that the Conley-Moser conditions are satisfied.

A different approach would be, as we have suggested, to define the curves and strips only on $M$; the pay-off here is obvious, but we are then forced to adopt new coordinates on $M$ and this clearly complicates matters in its own way. We are of the opinion that the former is the `lesser of two evils'.

Consequently we are forced to adopt the (cumbersome) notation that strips denoted \emph{with} an overbar stretch across $S$, whereas strips denoted \emph{without} an overbar stretch only across $M$. Thus $\bar{H}$ represents an $m_h$-horizontal strip whereas $H$ represents the corresponding intersection $\bar{H}\cap M$.

\begin{prop}
\label{prop:width_of_strips}
The proposition has two parts:
\begin{enumerate}
\item Let $\bar{H}$ be an $m_h$-horizontal strip such that $H=\bar{H}\cap M$ is contained in $\bigcup_{i\in I}H_i$ and let
\beq*
\tilde{H_i}=(F^{-j}\circ G^{-k}(H))\cap H_i.
\eeq*
Then there exists an $m_h$-horizontal strip $\bar{\tilde{H_i}}$ such that $\tilde{H_i}=\bar{\tilde{H_i}}\cap M$ and
\beq*
d(\tilde{H_i})\leqslant n_hd(H)
\eeq*
for some $0<n_h<1$.
\item Similarly, let $\bar{V}$ be an $m_v$-vertical strip such that $V=\bar{V}\cap M$ is contained in $\bigcup_{i\in I}V_i$ and let
\beq*
\tilde{V_i}=(G^k\circ F^j(V))\cap V_i.
\eeq*
Then there exists an $m_v$-vertical strip $\bar{\tilde{V_i}}$ such that $\tilde{V_i}=\bar{\tilde{V_i}}\cap M$ and
\beq*
d(\tilde{V_i})\leqslant n_vd(V)
\eeq*
for some $0<n_v<1$.
\end{enumerate}
\end{prop}

\begin{figure}[htp]
\centering
(a)\includegraphics[totalheight=0.2\textheight]{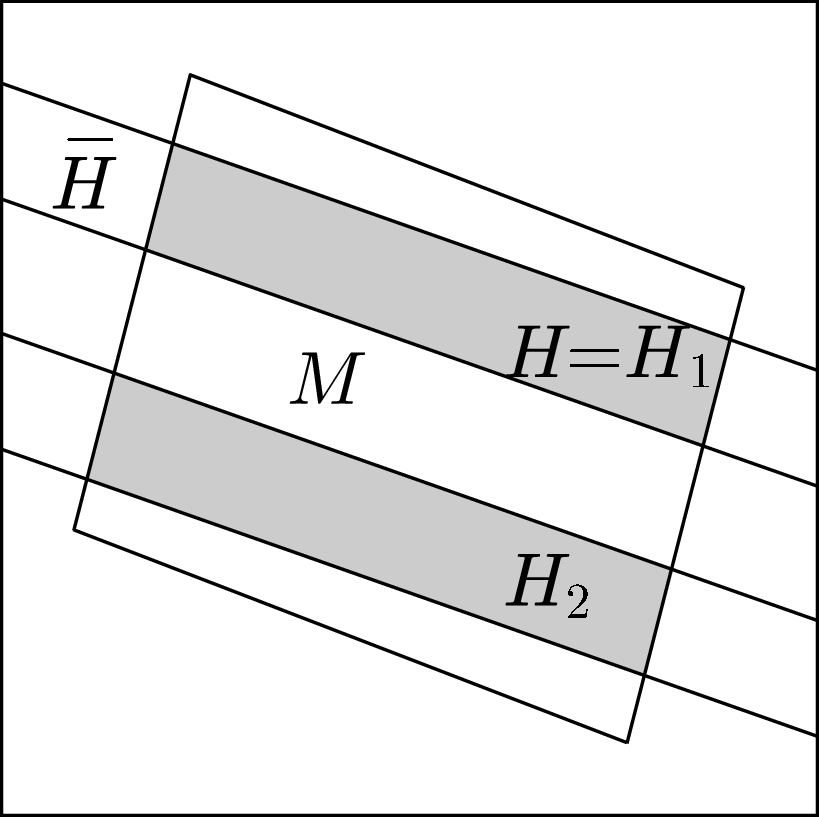}\quad
(b)\includegraphics[totalheight=0.2\textheight]{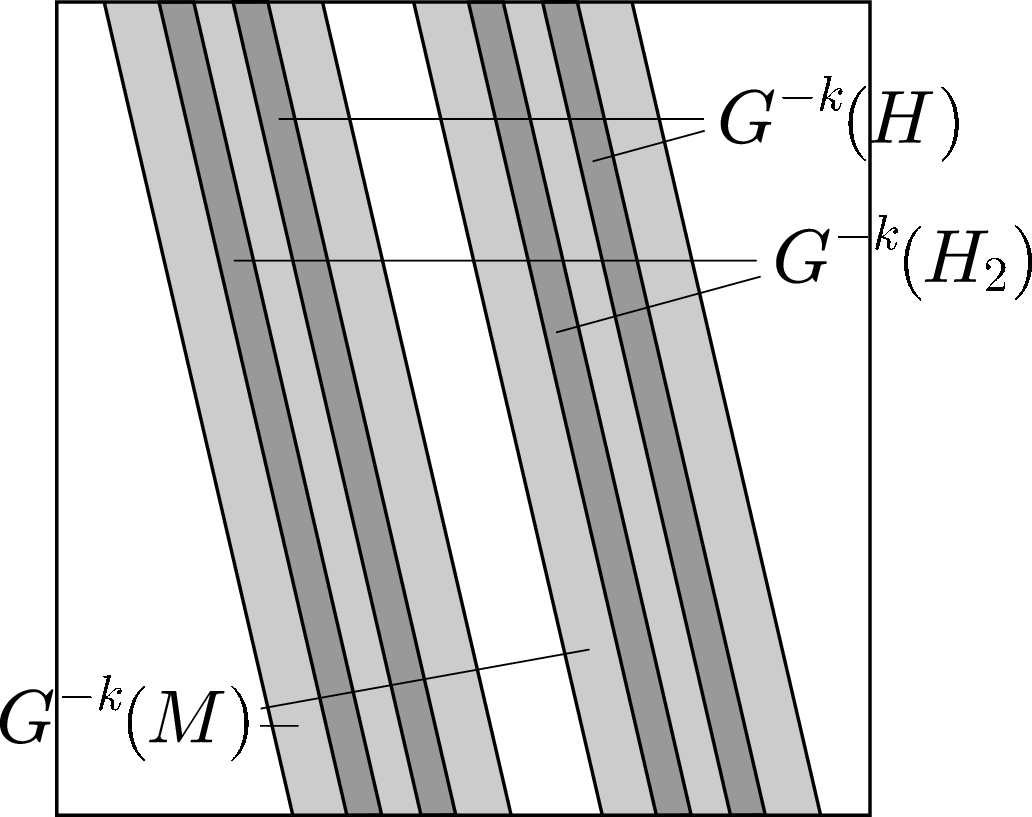}
\caption{Pre-image of horizontal strips: Part I}
\label{fig:final1}
\end{figure}

\begin{figure}[htp]
\centering
(a)\includegraphics[totalheight=0.2\textheight]{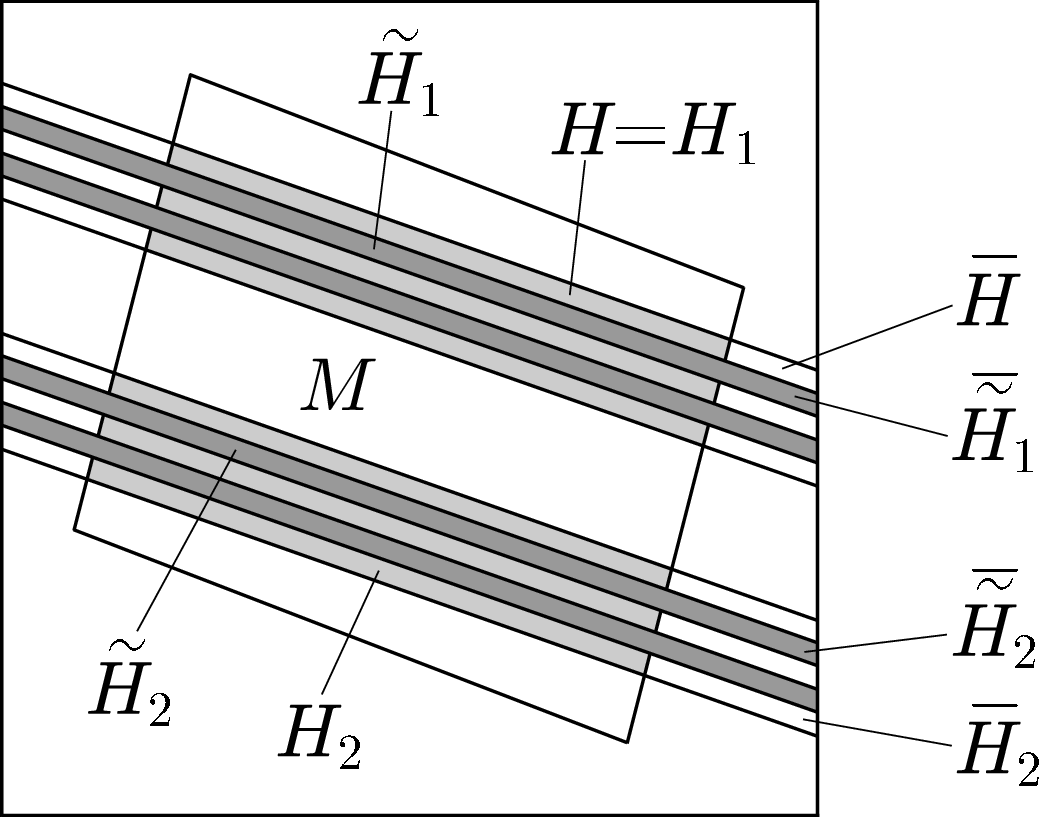}\quad
(b)\includegraphics[totalheight=0.2\textheight]{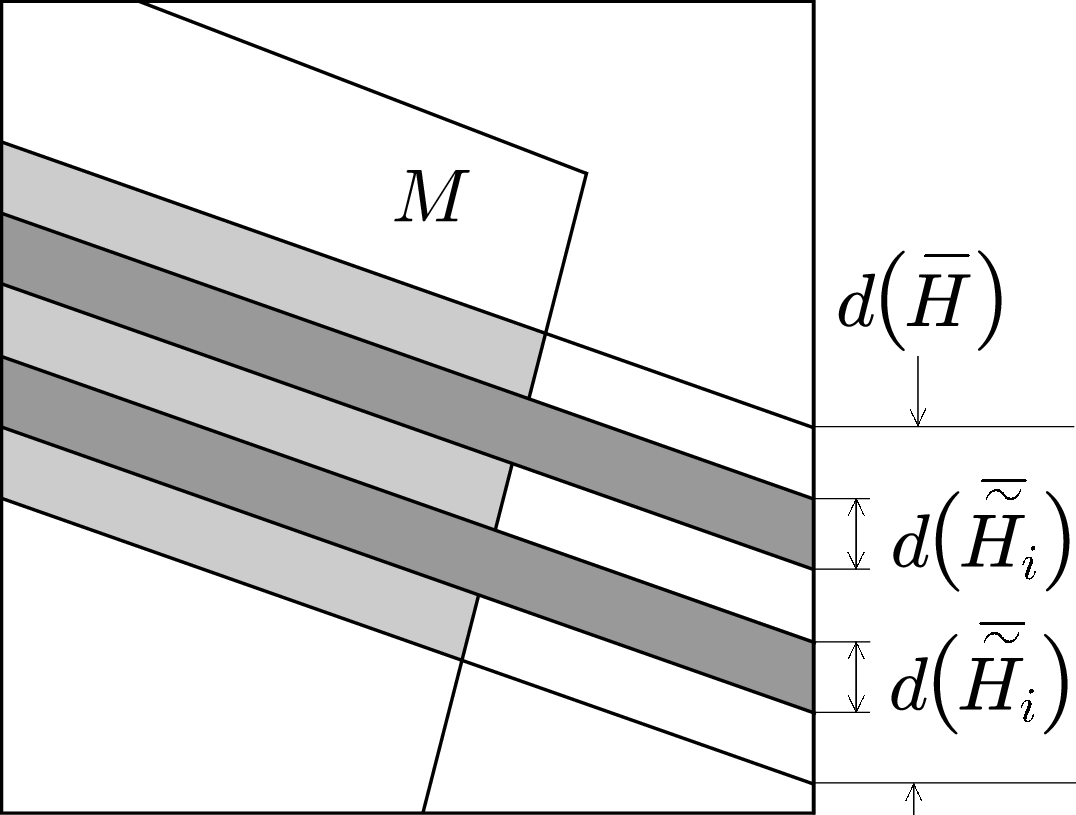}
\caption{Pre-image of horizontal strips: Part II}
\label{fig:final2}
\end{figure}

\begin{proof}
We prove just the first part, the second being similar.

Figure~\ref{fig:final1}(a) shows an $m_h$-horizontal strip $\bar{H}$ and two strips across $M$, $H_1=H=\bar{H}\cap M$ and $H_2$ (both shaded). Figure~\ref{fig:final1}(b) shows the pre-images of $M$, $H=H_1$ and $H_2$ with respect to $G^k$, illustrated with $k=3$. We observe that the boundaries of each are straight lines and that each of $G^{-k}(M)\cap S$, $G^{-k}(H)\cap S$ and $G^{-k}(H_2)\cap S$ consists of
$(k-1)$ $m_v$-vertical curves. Moreover $\left(G^{-k}(H\cup H_2)\cap S\right)\subset G^{-k}(M)\cap S$.

Recall from Figure~\ref{fig:strips2}(b) that $(F^{-j}\circ G^{-k}(M))\cap S$ consists of $(j+1)(k-1)$ disjoint pieces and that $(j-1)(k-1)$ of these are $m_h$-horizontal strips. Figure~\ref{fig:final2}(a) shows only these pieces. Notice that $F^{-j}\circ G^{-k}(H)$ (and similarly
$F^{-j}\circ G^{-k}(H_2)$) stretches completely across each piece and has straight-line boundaries. In other words, $(F^{-j}\circ G^{-k}(H))\cap S$ consists of $(j-1)(k-1)$ $m_h$-horizontal strips. We denote these $\left\{\bar{\tilde{H_i}}\right\}_{i\in I}$, so that
$\tilde{H_i}=\bar{\tilde{H_i}}\cap M$.

It remains to show that $d(\tilde{H_i})\leqslant n_hd(H)$ for some $0<n_h<1$. Each $H_i$ is the intersection of $M$ with an $m_h$-horizontal strip of width $d(\bar{H})$. Similarly each intersection $H_i\cap (F^{-j}\circ G^{-k}(H_{i'}))$ (for $i,i'\in I$), is the
intersection of $M$ with another $m_h$-horizontal strip of width $d(\bar{\tilde{H_i}})$. To obtain the result, consider Figure~\ref{fig:final2}(b); it is clear that
\beq*
(j-1)(k-1)d(\bar{\tilde{H_i}})\leqslant d(\bar{H}).
\eeq*
Simply observe that $d(\tilde{H_i})=d(\bar{\tilde{H_i}})$ and that $d(H)=d(\bar{H})$, and one obtains the required result with
\beq*
0<n_h=\frac{1}{(j-1)(k-1)}<1.
\eeq*
\end{proof}

We are now in a position to prove the existence of the horseshoe.

\begin{thm}
\label{thm:horseshoe}
Let each of $j,k$ and the product $N=(j-1)(k-1)$ be at least 2. Then the toral linked-twist map $G^k\circ F^j$ has an invariant Cantor set
$\Lambda$ on which it is topologically conjugate to a full shift on $N=(j-1)(k-1)$ symbols.
\end{thm}

\begin{proof}
We show that Propositions~\ref{prop:properties_of_M}, \ref{prop:mapping_of_h_strips} and~\ref{prop:width_of_strips} together imply that Conditions~\ref{cond:CM1}, \ref{cond:CM2} and~\ref{cond:CM3} hold. The result then follows from Theorem~\ref{thm:ConleyMoser}.

Proposition~\ref{prop:properties_of_M} says that $G^k\circ F^j(M)$ contains $N$ disjoint $m_v$-vertical strips, which do not intersect the boundaries $\tau_1$ or $\tau_2$ of $M$. Denote by $\bar{H}_i$ these strips and by $H_i$ their respective intersections with $M$. Similarly $F^{-j}\circ G^{-k}(M)$ contains $N$ disjoint $m_h$-horizontal strips which do not intersect $\sigma_1$ or $\sigma_2$. We denote these by $\bar{V}_i$ and denote by $V_i$ their respective intersections with $M$. We have $0<m_vm_h<1$, satisfying Condition~\ref{cond:CM1}.

It is important to notice that $\bar{H_i}\cap\bar{V_j}=H_i\cap V_j$ for all $i,j\in I$, i.e.\ no intersections between horizontal and vertical strips occur \emph{outside} of $M$.

Proposition~\ref{prop:mapping_of_h_strips} shows that $G^k\circ F^j$ acts as dictated by Condition~\ref{cond:CM2}, i.e.\ the $H_i$ are mapped homeomorphically onto the $V_i$ and horizontal (respectively, vertical) boundaries are mapped to horizontal (vertical) boundaries. Thus $G^k\circ F^j$ satisfies Condition~\ref{cond:CM2}.

Finally, Proposition~\ref{prop:width_of_strips} shows that $G^k\circ F^j$ acts as specified by Condition~\ref{cond:CM3}. In particular the
pre-image $F^{-j}\circ G^{-k}(H_i)$ forms another $N$ horizontal strips when it intersects with the $N$ strips $H_i$, and each of the new strips has width strictly less than the original. Analogous behaviour occurs for the images $G^k\circ F^j(V_i)$. Thus Proposition~\ref{prop:width_of_strips} shows that $G^k\circ F^j$ satisfies Condition~\ref{cond:CM3}.
\end{proof}

%% file: chapter2/section3.tex
\section{Uniform hyperbolicity of the horseshoe}\label{Horseshoe.UH}
\setcounter{equation}{0}

We conclude the chapter with a proof that the restriction of $G^k\circ F^j$ to the invariant set $\Lambda$ satisfies the uniform hyperbolicity conditions given in section~\ref{LitReview.Hyp}. This result should not be surprising to us. The toral linked-twist map we consider is seen to be a generalisation of the uniformly hyperbolic toral automorphism known as the cat map. The non-uniformity of the present system derives from the fact that stretching and contraction, the very essence of hyperbolicity, are a consequence of points entering and returning to $S$ (this feature is highlighted in Wojtkowski's (\citeyear{woj}) proof of the Bernoulli property). That we have no lower bounds on this return-time means that any uniform hyperbolicity constants we propose can be violated by some trajectory.

With $\Lambda$ the situation is different. By construction the return-time to $S$ is one iteration for each point. Let $z\in\Lambda$; we have the Jacobian
\beq*
D(G^k\circ F^j)_z= \begin{pmatrix} 1 & 0 \\ k\beta & 1 \end{pmatrix}
 \begin{pmatrix} 1 & j\alpha \\ 0 & 1 \end{pmatrix}
 =\begin{pmatrix} 1 & j\alpha \\ k\beta & 1+jk\alpha\beta \end{pmatrix},
\eeq*
which is independent of $z$ itself. Denote this Jacobian by $A$ for convenience, then the eigenvalues of $A$ are given by
\beq*
\lambda_{\pm}=\frac{1}{2}\left(jk\alpha\beta\pm\sqrt{(jk\alpha\beta)^2-4}\right).
\eeq*
It is easily checked that these are real and distinct, and moreover that
\beq*
0<\lambda_-<1<\lambda_-^{-1}=\lambda_+.
\eeq*
Eigenvectors of $A$ are given by
\beq*
v_{\pm}=\begin{pmatrix} j\alpha \\ \lambda_{\pm}-1 \end{pmatrix}.
\eeq*

Define subspaces $E^s(z)$ and $E^u(z)$ in the tangent space $T_z\mathbb{T}^2$ to be the spans of $v_-$ and $v_+$ respectively. Because these vectors are linearly independent they form a basis for the tangent space, i.e.\ $T_z\mathbb{T}^2=E^s(z)\oplus E^u(z)$. The properties required of a uniformly hyperbolic system are easily satisfied with $C=1$ and $\lambda=\lambda_-$; indeed
\beq*
AE^s(z)=\left\{Apv_-:p\in\mathbb{R}\right\}=\left\{p\lambda_-v_-:p\in\mathbb{R}\right\}=E^s\left(G^k\circ F^j(z)\right),
\eeq*
and similarly $AE^u(z)=E^u\left(G^k\circ F^j(z)\right)$. Together these satisfy (\ref{eqn:UH1}). Finally let $v_s\in E^s(z)$ and $v_u\in E^u(z)$, then
\beq*
\|A^nv_s\|=\|\lambda_-^nv_s\|=|\lambda_-^n|\|v_s\|\and\|A^{-n}v_u\|=\|\lambda_+^{-n}v_u\|=|\lambda_-^n|\|v_u\|,
\eeq*
where $\|\cdot\|$ denotes the standard Euclidean norm; these satisfy (\ref{eqn:UH2}) and (\ref{eqn:UH3}) respectively.

%% file: chapter3/chapter3.tex
\chapter{The Bernoulli property for a planar linked-twist map}
\setcounter{equation}{0}
\setcounter{figure}{0}
\label{chapter3}

In this chapter we establish the Bernoulli property for the planar linked-twist map defined in Section~\ref{Intro.Theorems.Plane}. Our starting point is the work of \citet{woj} who proved that such systems can have an ergodic partition. We outline his criteria for this in Section~\ref{PLTM.Wojtkowski}.

A crucial element of our proof is the introduction of new coordinates on the manifold $A$ which will enable us to improve Wojtkowski's estimates on the orientation of local invariant manifolds. We introduce these in Section~\ref{PLTM.New_Coords}. They simplify the proof that in the case $2\leqs r_0<r_1\leqs\sqrt{7}$, Wojtkowski's condition for an ergodic partition is satisfied and we give this proof in the same section.

In Section~\ref{PLTM.Manifolds} we describe the local unstable manifold $\gamma^u(w)$ of a point $w=(u,v)\in A$. Following Wojtkowski's lead it will be convenient to consider the unstable manifolds $\gamma^u(\omega)$, where $\omega=(r,\theta)=M_+^{-1}(w)$. We define its length and discuss how this definition might be extended to $\Theta_{\Sigma}^n(\gamma^u(\omega))$. There are some technical considerations but we show that for $\mu$-a.e.\ such $w$ our definition does indeed hold. We conclude by showing that the length of $\Theta_{\Sigma}^n(\gamma^u(\omega))$ grows exponentially with $n$.

In Section~\ref{PLTM.Cone_C} we express the planar linked-twist map $\Theta$ in our new coordinates. We then show that a certain tangent cone is preserved by the differential $D\Theta$. Finally, in Section~\ref{PLTM.Bernoulli} we are able to give the estimate we have mentioned on the orientation of the unstable manifolds and, following this, a largely geometrical proof that the Bernoulli property is satisfied.

\input chapter3/section1

\input chapter3/section2

\input chapter3/section3

\input chapter3/section4

\input chapter3/section5

%% file: chapter3/section1.tex
\section{Wojtkowski's results}\label{PLTM.Wojtkowski}
\setcounter{equation}{0}

In this section we describe Wojtkowski's (\citeyear{woj}) criteria for the planar linked-twist map to have an ergodic partition.

Let $w=(u,v)\in\Sigma$ and denote by $\alpha(w)\in(0,\pi)$ the angle at which the segment connecting $w$ to $(-1,0)$ meets the segment connecting $w$ to $(1,0)$. Let
\beql
\label{eqn:alpha_condition}
\eta=\sup_{w\in\Sigma}\frac{\cot\alpha(w)}{r(w)},
\eeql
where $r(w)$ denotes the Euclidean distance from $w$ to $(-1,0)$. We also denote by $c$ the infimum of the derivative of the twist function, which in this case is just given by $2\pi/(r_1-r_0)$. Wojtkowski proved the following:
\begin{thm}[\citet{woj}] \label{thm:wojt_erg_part}
If
\beql
\label{eqn:condn_W}
c>2\eta
\eeql
then the linked-twist map $\Theta:A\to A$ is the union of (at most) countably many ergodic components.
\end{thm}
Wojtkowski's \emph{conjecture}, made in the same paper, is that under the assumptions of Theorem~\ref{thm:wojt_erg_part} then $\Theta:A\to A$ has the $K$-property. This would, by the work of \citet{ch}, imply that it has the Bernoulli property.

We discuss the proof of Theorem~\ref{thm:wojt_erg_part} briefly. It is easily argued that $\mu$-a.e.\ $w\in A$ lands in $\Sigma$ under iteration of $\Theta$ and, furthermore, returns to $\Sigma$ infinitely many times.\footnote{One simply notices that those points \emph{not} satisfying this condition must be rigid rotations around one of the annuli, and must have rational angle of rotation, else their orbit would be dense and hit $\Sigma$. From the nature of the twist function (in particular the strict monotonicity) one easily infers that such points are contained within a set of measure zero.} Thus for a full-measure set of points we may talk of the \emph{return map to} $\Sigma$, or just the \emph{return map} as we shall usually abbreviate it. Following \citet{woj} we shall actually define the return map on $M_+^{-1}(\Sigma)\subset L$ rather than on $\Sigma$ itself, as follows: 

\begin{defn}[First-return map to $\Sigma$]
Let $(r,\theta)\in M_+^{-1}(\Sigma)\subset L$. The first-return map to $\Sigma$ is the map $\Theta_{\Sigma}:M_+^{-1}(\Sigma)\to M_+^{-1}(\Sigma)$ given by
\beq*
\Theta_{\Sigma}=M_+^{-1}\circ\Theta^i\circ M_+,
\eeq*
where $i$ is the smallest (strictly) positive integer for which $\Theta^i(M_+(r,\theta))\in\Sigma$.
\end{defn}

\begin{figure}[htp]
\centering
\includegraphics[totalheight=0.16\textheight]{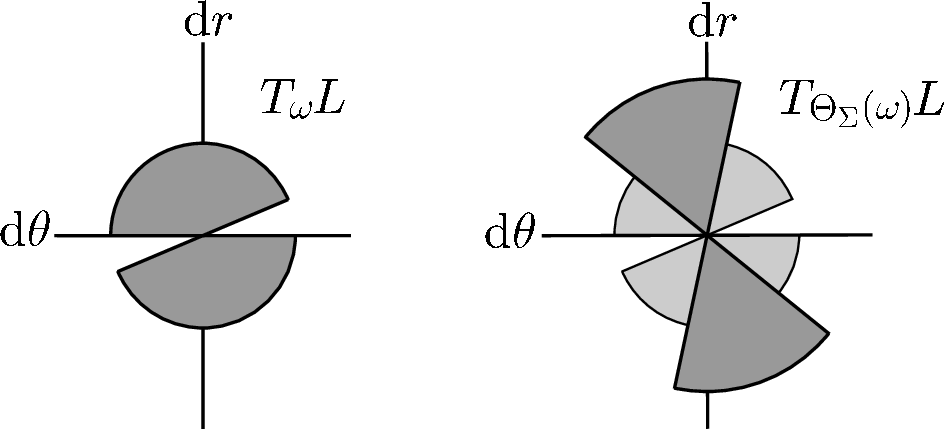}
\caption[The cone $U$]{The invariant expansive cone $U\subset T_{\omega}L$ is shown in the left-hand figure. In the right-hand figure is the image of the cone under the differential map $D\Theta_{\Sigma}$ (dark-shaded) with the original cone (light-shaded) included for comparison. Observe how the cone is mapped into itself and vectors within it are expanded.}
\label{fig:U}
\end{figure}

For $(r,\theta)\in M_+^{-1}(\Sigma)$ let $\beta_1=\text{d}r$, $\beta_2=\text{d}\theta$ give coordinates in the tangent space $T_{(r,\theta)}L$ and define the cone
\beq*
U(r,\theta)=\left\{(\beta_1,\beta_2):\frac{\beta_2}{\beta_1}\geqslant\frac{-c}{2}\right\}.
\eeq*
Wojtkowski establishes that $U$ is invariant under, and expanded by, the derivative $D\Theta_{\Sigma}$. We illustrate the situation in Figure~\ref{fig:U}. More precisely, define the cone field
\beq*
U_+=\bigcup_{(r,\theta)\in M_+^{-1}(\Sigma)}U(r,\theta)
\eeq*
and let $\|\cdot\|$ be the norm in $T_{(r,\theta)}L$ induced by the Riemannian metric, i.e.\ $\|(\beta_1,\beta_2)\|=\sqrt{\beta_1^2+r^2\beta_2^2}$. We have the following:

\begin{prop}[\citet{woj}]
\label{prop:wojtkowski_planar_cones}
$D\Theta_{\Sigma}(U_+)\subset U_+$. Furthermore there is a constant $\lambda>1$, independent of $(r,\theta)$ or $\beta$, and for vectors $\beta\in U_+$ we have $\|D\Theta_{\Sigma}\beta\|\geqslant\lambda\|\beta\|$.
\end{prop}
For a proof see the original paper or \citet{sturman}.

%% file: chapter3/section2.tex
\section{New coordinates for the manifold $A$}\label{PLTM.New_Coords}
\setcounter{equation}{0}

At the heart of our proof of the Bernoulli property for $\Theta$ is a new coordinate system. In this section we introduce these coordinates and with them prove that if the condition
\beql
\label{eqn:size_of_annuli}
2\leqs r_0<r_1\leqs\sqrt{7}
\eeql
is satisfied, then the linked-twist map $\Theta$ has an ergodic partition. Theorem~\ref{thm:wojt_erg_part} says that this amounts to proving that the condition~(\ref{eqn:condn_W}) is satisfied.

We use the notation $\mathbb{S}^1=[-\pi,\pi]$ with opposite ends identified, $\mathcal{R}=[r_0,r_1]$ and $-\mathcal{R}=[-r_1,-r_0]$.

\subsection{Construction of the new coordinates on $A_+$}

We introduce new coordinates $(x(u,v),y(u,v))$ for the manifold $A$. The key feature is that if $(u,v)\in\Sigma_+$ then $x(u,v)$ is given by the Euclidean distance from $(u,v)$ to the centre of annulus $A_+$ and $y(u,v)$ is given by the distance from $(u,v)$ to the centre of annulus $A_-$.\footnote{Coordinates such as these are commonly called \emph{two-centre bi-polar coordinates}, however we do not adopt this name because we will not extend the coordinates to all of $A$ in this manner.} Thus for $(u,v)\in\Sigma_+$ we have
\beql
\label{eqn:xy_condition}
(x,y)=\left(\sqrt{(1+u)^2+v^2},\sqrt{(1-u)^2+v^2}\right).
\eeql
Similarly on $\Sigma_-$ the \emph{magnitudes} of $x$ and $y$ are determined in this way, but one of the coordinates assumes a negative value so that points are uniquely defined. Unfortunately it is not useful to extend $x,y$ to all of $A$ in the same manner so we will take an alternative approach. We begin by defining the coordinates on $A_+$.

Let $(u,v)\in A_+$ and let $(r,\theta)=M_+^{-1}(u,v)\in L=\mathcal{R}\times\mathbb{S}^1$. Setting $x=r$ satisfies the first part of (\ref{eqn:xy_condition}) because $r$ is the Euclidean distance from $(u,v)$ to $(-1,0)$. However the $\theta$ coordinate will \emph{not} in general give the Euclidean distance to $(1,0)$.

\begin{figure}[htp]
\centering
\includegraphics[totalheight=0.4\textheight]{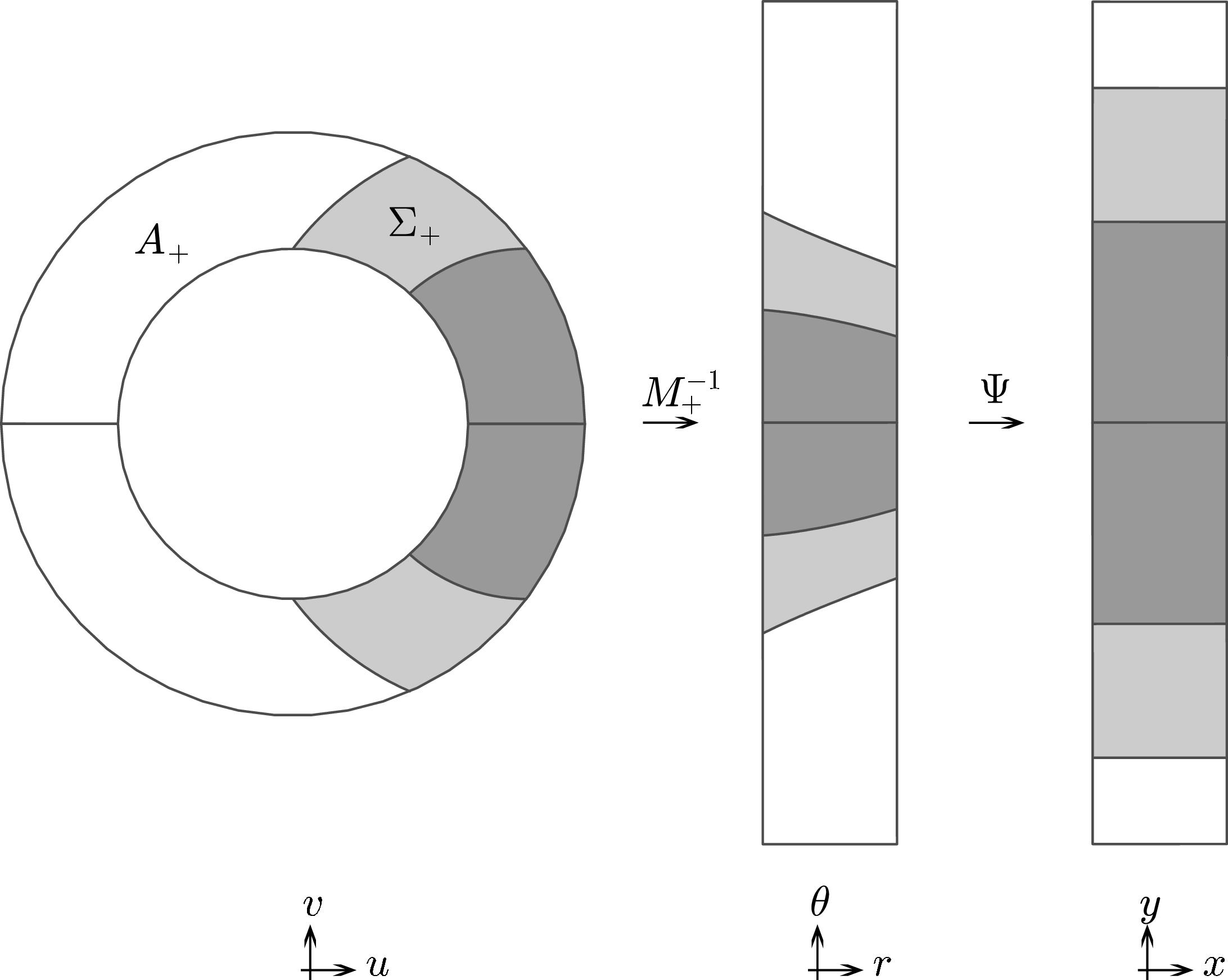}
\caption[Coordinate systems on $A+$]{The region $A_+\subset\mathbb{R}^2$, illustrated in the three coordinate systems. Left-to-right: Cartesians $(u,v)$ in the plane; polars $(r,\theta)\in\mathbb{R}_0^+\times\mathbb{S}^1$; and new coordinates $(x,y)\in\mathbb{S}^1\times\mathbb{S}^1$. Shading indicates the three regions for which $\Psi$ takes different forms, as explained in the text.}
\label{fig:new_coords}
\end{figure}

We now wish to define a homeomorphism $\Psi:\mathcal{R}\times\mathbb{S}^1\to \mathcal{R}\times\mathbb{S}^1$ so that
\beq*
(x,y)=\Psi\circ M_+^{-1}(u,v)
\eeq*
is as required. It is clear that $\Psi$ will need to have the form
\beq*
\Psi(r,\theta)=(r,\psi(r,\theta)),
\eeq*
for some function $\psi:\mathcal{R}\times\mathbb{S}^1\to\mathbb{S}^1$ which is to be determined. It suffices to define a homeomorphism $\psi:\mathcal{R}\times[0,\pi]\to[0,\pi]$, for which $\psi(r,0)=0$, $\psi(r,\pi)=\pi$ and with the extra condition that
\beql\label{eqn:psi_symmetry}
\psi(r,-\theta)=-\psi(r,\theta),
\eeql
so that (\ref{eqn:xy_condition}) holds. Once we have done this, we will extend our coordinates to $A_-$. Figure~\ref{fig:new_coords} should help the reader to keep track of our construction.

Let $(u,v)\in\Sigma_+$, which is one of the light-shaded regions in Figure~\ref{fig:new_coords}. We require that $\psi(r,\theta)=\sqrt{(1-u)^2+v^2}$. Expressing $u,v$ in terms of $r,\theta$ gives
\beql
\label{eqn:psi1}
\psi(r,\theta) = \sgn(\theta)\sqrt{r^2-4r\cos\theta+4}.
\eeql
This defines $\psi$ (and thus the new coordinates $(x,y)$) on the domain $\Sigma=\Sigma_+\cup\Sigma_-$. It is clear from the construction (of two-centre bi-polar coordinates) that this is a homeomorphism.

Before we continue it will be convenient to determine the inverse to (\ref{eqn:psi1}) on $\mathcal{R}\times \mathcal{R}$ , i.e.\ the function $\psi^{-1}:\mathcal{R}\times \mathcal{R}\to M_+^{-1}(\Sigma)$ such that
\beq*
\psi(x,\psi^{-1}(x,y))=y.
\eeq*
We claim that the function $\cos^{-1}\left(\tau(x,y)\right)$, where $\tau(x,y)=\frac{1}{4x}\left(x^2-y^2+4\right)$, has this property. Indeed:
\beqa*
\psi\left(x,\cos^{-1}\left(\tau(x,y)\right)\right) &= \sgn(y)\sqrt{x^2-4x\tau(x,y)+4} \\
																									 &= \sgn(y)\sqrt{x^2-x^2+y^2-4+4} = y.
\eeqa*

We now define new coordinates on $A_+\backslash\Sigma$, which consists of two disjoint pieces (recall that we are considering only the `top half' as illustrated; we will deal with the bottom half using the symmetry given by (\ref{eqn:psi_symmetry})). We deal with each separately. Our approach to this is very simple though perhaps this simplicity can be lost amongst the equations we provide. To counter this let us describe heuristically what it is that we are doing first.

Consider Figure~\ref{fig:new_coords} again and look at the middle figure, in $(r,\theta)$ coordinates. We have constructed a homeomorphism from $M_+^{-1}(\Sigma_+)$, which is the upper, light-shaded area, to the corresponding light-shaded area in the right-hand figure. Now we want to do the same for the dark-shaded area immediately below $M_+^{-1}(\Sigma_+)$, that is to create a homeomorphism from it to the corresponding dark-shaded area in the right-hand figure. The simplest way to do this is to leave the points such that $\theta=0$ invariant and `stretch' the others, `upwards' as illustrated, whilst leaving the $r$-coordinate unchanged.

More formally: let $(u,v)$ be a point in that part of $A_+\backslash\Sigma$ which lies \emph{inside} the annulus $A_-$. This is the dark-shaded region in Figure~\ref{fig:new_coords}. $M_+^{-1}(u,v)$ has coordinates $(r,\theta)$ such that $r\in \mathcal{R}$ and $\theta\in[0,\cos^{-1}(\tau(r,r_0))]$. For each $r\in \mathcal{R}$ we define $\psi$ on $\{r\}\times [0,\cos^{-1}(\tau(r,r_0))]$ by
\beql
\label{eqn:psi2}
\psi(r,\theta)=r_0\theta/\cos^{-1}(\tau(r,r_0)).
\eeql

Finally let $(u,v)$ denote a point in that part of $A_+\backslash\Sigma$ which lies \emph{outside} of $A_-$. This is the unshaded region in Figure~\ref{fig:new_coords}. Here $M_+^{-1}(u,v)$ is given by coordinates $(r,\theta)\in \mathcal{R}\times[\cos^{-1}(\tau(r,r_1)),\pi]$. For each $r\in \mathcal{R}$ we define $\psi$ on $\{r\}\times[\cos^{-1}(\tau(r,r_1)),\pi]$ by
\beql
\label{eqn:psi3}
\psi(r,\theta)=r_1+\frac{\theta-\cos^{-1}(\tau(r,r_1))}{\pi-\cos^{-1}(\tau(r,r_1))}(\pi-r_1).
\eeql
This, together with the symmetry requirement (\ref{eqn:psi_symmetry}) completes our definition of $x(u,v),y(u,v)$ on $A_+$. To summarise:
\begin{defn}[Coordinates $(x,y)$ on $A_+$]
Let $(u,v)\in A_+$, then $(x,y)=\Psi\circ M_+^{-1}(u,v)$, where $\Psi(r,\theta)=(r,\psi(r,\theta))$ and $\psi$ takes one of the forms (\ref{eqn:psi1}), (\ref{eqn:psi2}) or (\ref{eqn:psi3}) as described.
\end{defn}

\subsection{Construction of the new coordinates on $A_-$}

\begin{figure}[htp]
\centering
\includegraphics[totalheight=0.25\textheight]{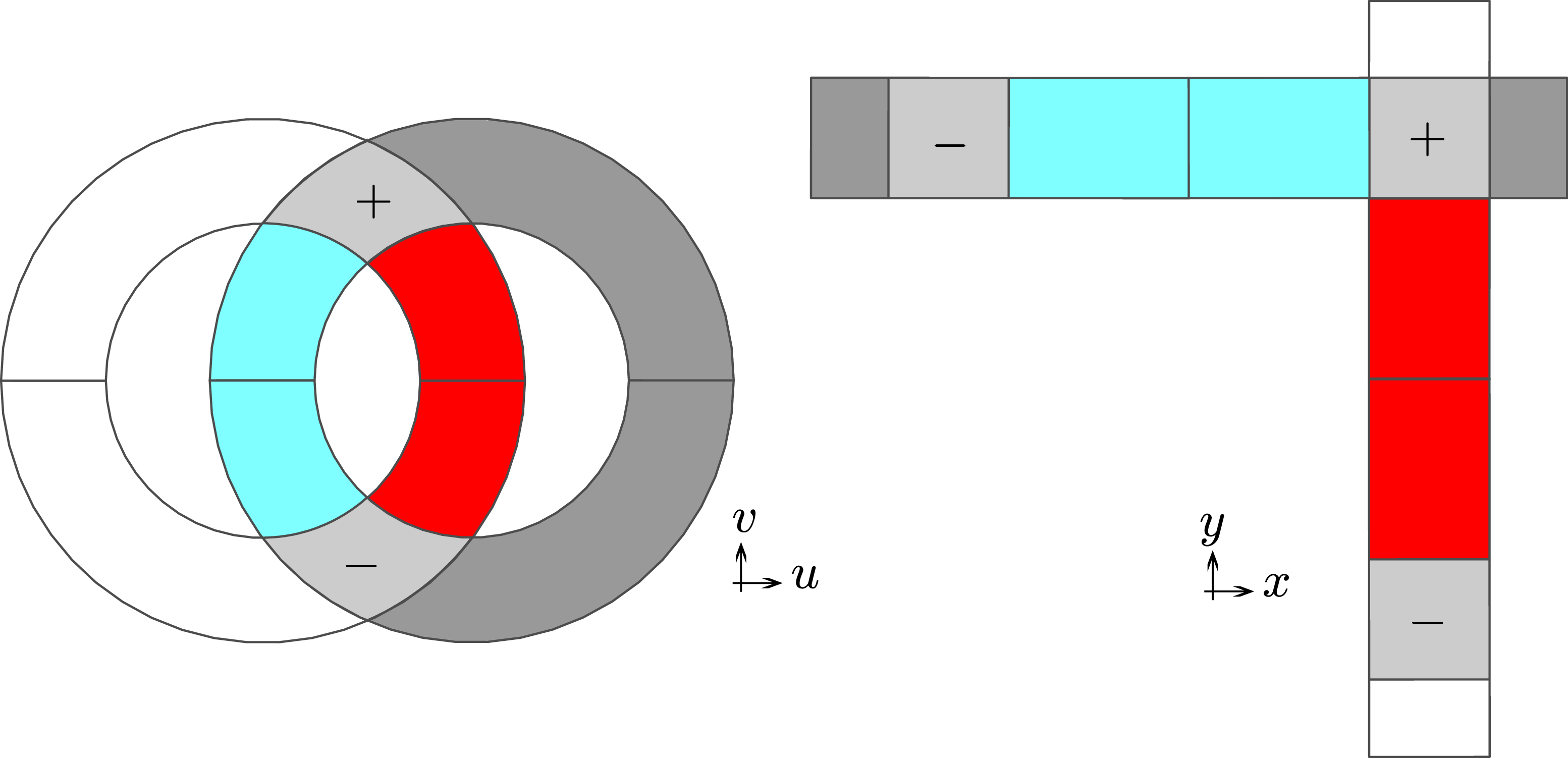}
\caption[Coordinate transformation for $A$]{The manifold $A$, illustrated in its native Cartesian coordinates and in the new coordinates $(x,y)\in\mathbb{S}^1\times\mathbb{S}^1$. Notice that there are two distinct representations of $\Sigma_-$ deriving from the coordinate transformations. In $(x,y)$-coordinates this image is given by $\mathcal{R}\times-\mathcal{R}$ (toward the bottom-right of the right-hand figure), whereas in $(x',y')$-coordinates it is given by $-\mathcal{R}\times\mathcal{R}$ (toward the top-left of the right-hand figure).}
\label{fig:new_coords_all}
\end{figure}

There is a slight problem to overcome in extending the coordinates to $A_-$. Because of the geometry of $A$, the most natural approach (i.e.\ straight-forward symmetry) will lead to us having two different expressions for the new coordinates on $\Sigma_-$. Our strategy is to pursue this naive approach anyway and deal with the problem when it arises; in acknowledgment of this we call the coordinates on $A_-$ `temporary' for the moment. We hope that Figure~\ref{fig:new_coords_all} will help the reader.

We will define temporary new coordinates $(x',y')$ on $A_-$ before extending $(x,y)$ to that domain. We introduce the functions $\iota,\iota^{-1}:\mathbb{R}^2\to\mathbb{R}^2$ given by
\beq*
\iota(x,y)=(y,-x)\and\iota^{-1}(x,y)=(-y,x).
\eeq*
With $\Psi$ as already defined, on $A_-$ we define $(x',y')=\iota^{-1}\circ\Psi\circ M_-^{-1}(u,v)$. Notice that $M_-^{-1}$ expresses $A_-$ in polars (albeit centred at $(1,0)$ and with polar angles shifted by $\pi$), so that the restriction of $\Psi$ to $M_-^{-1}(A_-)$ is a homeomorphism, just as the restriction of $\Psi$ to $M_+^{-1}(A_+)$ was. For $(u,v)\in\Sigma_{\pm}$ we have
\beqa*
(x',y')&=\iota^{-1}\circ\Psi\circ M_-^{-1}(u,v) \\
       &=\iota^{-1}\left(\sqrt{(1-u)^2+v^2},\mp\sqrt{(1+u)^2+v^2}\right) \\
       &=\left(\pm\sqrt{(1+u)^2+v^2},\sqrt{(1-u)^2+v^2}\right).
\eeqa*
Comparing the coordinates $(x,y)$ and $(x',y')$ where their domains overlap we see that $(x',y')=\pm(x,y)$ on $\Sigma_{\pm}$. See Figure~\ref{fig:new_coords_all}. We thus extend the coordinates $(x,y)$ to $A_-$ as follows
\begin{defn}[Coordinates $(x,y)$ on $A_-$]
For $(u,v)\in A_-\backslash\Sigma_-$ define new coordinates $(x,y)=(x',y')=\iota^{-1}\circ\Psi\circ M_-^{-1}(u,v)$. For $(u,v)\in\Sigma_-$ define $(x,y)=-(x',y')$.
\end{defn}

\subsection{The ergodic partition}

The new coordinates simplify our proof that the planar linked-twist map $\Theta:A\to A$ defined with radii $2\leqs r_0<r_1\leqs\sqrt{7}$ has an ergodic partition, as defined within Theorem~\ref{thm:ks}. This amounts to the following:
\begin{lem}
Let $2\leqs r_0<r_1\leqs\sqrt{7}$. Then condition (\ref{eqn:condn_W}) is satisfied.
\end{lem}
\begin{proof}
By symmetry it is enough to show that the condition holds on $\Sigma_+$. Moreover the condition is implied by
\beql
\label{eqn:new_alpha_condition}
\sup_{w\in\Sigma_+}\cot\alpha(w)<\frac{\pi}{\sqrt{7}-2}.
\eeql
Let $(x,y)=\Psi\circ M_+^{-1}(w)$ give $w$ in the new coordinates. The angle $\alpha$ appears in a triangle in which its adjacent sides have lengths $x$ and $y$, and the opposite side has length $2$. The law of cosines says that
\beql
\label{eqn:cos_alpha}
\cos\alpha=\frac{x^2+y^2-4}{2xy}.
\eeql
The partial derivative of (\ref{eqn:cos_alpha}) with respect to $x$ is given by
\beq*
\frac{\pd}{\pd x}(\cos\alpha)=\frac{1}{2y}-\frac{y^2-4}{2x^2y},
\eeq*
and, using $x,y\in[2,\sqrt{7}]$, we calculate that $1/2y\in[1/2\sqrt{7},1/4]$ and $(y^2-4)/2x^2y\in[0,3/16]$. It is easily checked that $1/2\sqrt{7}>3/16$ and so the derivative is always positive. Consequently $\cos\alpha$ is an increasing function of $x$ and thus $\alpha$ is a decreasing function of $x$. By symmetry $\alpha$ is also a decreasing function of $y$.

Combining these facts with (\ref{eqn:cos_alpha}) we find that $\alpha\in[\cos^{-1}(5/7),\pi/3]$. Recall that $\cot$ is positive and decreasing on $(0,\pi/2)$, then
\beq*
\sup_{w\in\Sigma_+}\cot\alpha
=\cot\inf_{w\in\Sigma_+}\alpha
=\cot\cos^{-1}\frac{5}{7}
=\frac{5\sqrt{6}}{12},
\eeq*
and (\ref{eqn:new_alpha_condition}) is seen to be satisfied.
\end{proof}

%% file: chapter3/section3.tex
\section{Growth of local invariant manifolds}\label{PLTM.Manifolds}
\setcounter{equation}{0}

In this section we describe the nature of local unstable manifolds and their images under $\Theta$, and define what we mean by their length. We show that this length diverges on iteration of the map. We remark that \citet{woj} stated that these results are true but did not give a proof.

\subsection{Nature of local unstable manifolds}

At $\mu$-a.e.\ point $(u,v)\in\Sigma$ there is a positive Lyapunov exponent. This follows easily from Wojtkowski's (\citeyear{woj}) construction of the invariant expansive cone $U$, although he does not prove it explicitly (his work pre-dated \citet{ks} so he was unable to use their extension of Pesin theory, and thus had no reason to discuss Lyapunov exponents). The details may be found in \citet{sturman}, however.

Theorem~\ref{thm:ks} (due to \citet{ks}) tells us that as a consequence, associated to each such point is a local unstable manifold. Let $\omega=(r,\theta)=M_+^{-1}(u,v)$. The local unstable manifold is denoted $\gamma^u(\omega)$ and has the form
\beql
\label{eqn:form_of_um}
\gamma^u(\omega)=\exp_{\omega}\{(\omega,\phi^u(\omega)):\omega\in B^u(0,\eps)\}
\eeql
for some $\eps>0$. Here, $B^u(0,\eps)$ is the open $\eps$-neighbourhood of the origin in the unstable subspace $E^u(\omega)\subset T_{\omega}L$ and $\phi^u:B^u(0,\eps)\to E^s(\omega)\subset T_{\omega}L$ is a smooth map satisfying $\phi^u(0)=0$ and $D\phi^u(0)=0$.

We prove that the planar linked-twist map $\Theta:A\to A$ has the Bernoulli property by demonstrating that the `strong' form of the condition~(\ref{eqn:mip1}) is satisfied. We recall that $\mu$-a.e.\ point returns infinitely many times to $\Sigma$ and that we are defining local invariant manifolds in $M_+^{-1}(\Sigma)$ rather than in $\Sigma$ itself. The condition is that for $\mu$-a.e.\ $w,w'\in\Sigma$ and for all sufficiently large integers $m$ and $n$ (depending on $w,w'$) we have
\beql
\label{eqn:mip_pltm}
\Theta^n\left(M_+\left(\gamma^u(M_+^{-1}(w))\right)\right) \cap \Theta^{-m}\left(M_+\left(\gamma^s(M_+^{-1}(w'))\right)\right) \neq\emptyset.
\eeql

There are two key components to our proof that the condition (\ref{eqn:mip_pltm}) is satisfied. The first, which we deal with in this section, is that local unstable manifolds grow under iteration of $\Theta$. We make this notion precise once we have defined what exactly we mean by their \emph{length}. The second component is a useful characterisation of their \emph{direction} or \emph{orientation}. The latter is our principle motivation for having introduced the new coordinates and we return to it in the following section.

\subsection{Length of local unstable manifolds}

We recall some terminology from differential geometry, taken again from \citet{docarmo}: if $M$ is a smooth manifold and $I\subset\mathbb{R}$ an open interval, then a \emph{smooth curve} is a smooth map $\alpha:I\to M$. The image $\gamma=\alpha(I)\subset M$ is referred to as the \emph{trace} of the curve.

\begin{defn}[Length of a smooth curve; length of its trace]
Let $M$ be a smooth manifold, $I\subset\mathbb{R}$ an open interval and $\alpha:I\to M$ a smooth curve. We define the length of the curve $\alpha$ to be
\beql
\label{eqn:curve_length}
\text{length}(\alpha)=\int_I\left\|D\alpha_t\right\|\d t,
\eeql
where $D\alpha_t$ denotes the derivative (i.e.\ the Jacobian) of $\alpha$ evaluated at $t\in I$ and $\|\cdot\|$ denotes the Euclidean norm in $T_{\alpha(t)}M$. We define the length of the trace $\gamma=\alpha(I)$ to coincide with the length of the curve $\alpha$.
\end{defn}

We remark that it is perhaps more common in the differential geometry literature to see the expression $\d\alpha/\d t$ in place of $D\alpha_t$ but it will be more convenient for us to use the latter expression. We also comment that this definition of length can be shown to agree with our geometrical intuition. For more details see \citet{docarmo}.

Now let $\omega=(r,\theta)\in M_+^{-1}(\Sigma)\subset L$ be such that $\gamma^u(\omega)$ exists and is of the form (\ref{eqn:form_of_um}). The unstable subspace and thus the unstable manifold are one-dimensional, so there is an open interval $I\subset\mathbb{R}$ and a diffeomorphism $\alpha:I\to M_+^{-1}(A)$ (i.e.\ the map is diffeomorphic between $I$ and its image in $M_+^{-1}(A)$) such that $\alpha(I)=\gamma^u(\omega)$. The diffeomorphism $\alpha$ is a smooth curve of which $\gamma^u(\omega)$ is the trace, so the length of $\gamma^u(\omega)$ is defined as above.

\begin{figure}[htp]
\centering
\includegraphics[totalheight=0.3\textheight]{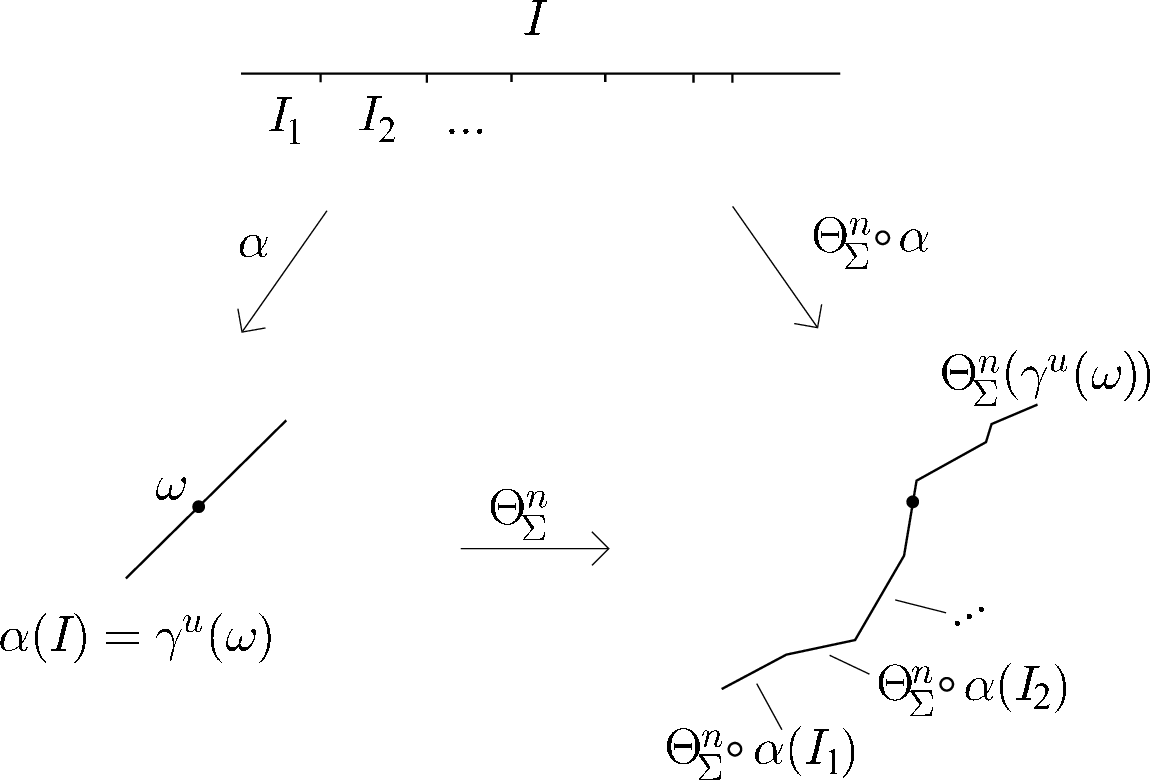}
\caption[Length of $\gamma^u(\omega)$ and of $\Theta_{\Sigma}^n(\gamma^u(\omega))$]{The length of local unstable manifold $\gamma^u(\omega)$ is well-defined because there is a diffeomorphism $\alpha$ of an open interval $I$ and $\alpha(I)=\gamma^u(\omega)$. If $I$ decomposes into a countable (at most) family of intervals on which $\Theta_{\Sigma}^n\circ\alpha$ is differentiable then the length of $\Theta_{\Sigma}^n(\gamma^u(\omega))$ is well-defined also.}
\label{fig:manifolds}
\end{figure}

Now let $n\in\mathbb{N}$ and consider the image of $\gamma^u(\omega)$ with respect to $\Theta^n_{\Sigma}$. In general $\Theta^n_{\Sigma}\circ\alpha:I\to M_+^{-1}(A)$ will \emph{not} be a smooth curve so it is not immediately clear that the length of $\Theta^n_{\Sigma}(\gamma^u(\omega))$ is defined. Suppose however that we are able to determine that $\Theta^n_{\Sigma}\circ\alpha$ is \emph{piecewise smooth}, in the sense that $I$ decomposes into a countable (at most) union of intervals
\beql
\label{eqn:partition_of_I}
I=(i_0,i_1)\cup\bigcup_{h=1}^{\infty}[i_h,i_{h+1})
\eeql
and that the restriction of $\Theta^n_{\Sigma}\circ\alpha$ to any interval $(i_h,i_{h+1})$ is smooth. In this case our definition of length applies to each restriction and we naturally determine the length of $\Theta^n_{\Sigma}\circ\alpha$ by summing over them. We illustrate the situation in Figure~\ref{fig:manifolds}.

Determining that such a partition of $I$ exists for a `typical' point $\omega\in M_+^{-1}(\Sigma)$ is the main technical difficulty in proving that the length of its local unstable manifold grows under iteration with $\Theta$. Our ambition for the remainder of this section is to prove the following theorem:

\begin{thm}
\label{thm:growth_of_unstable_manifolds}
Let $n\in\mathbb{N}$ and let $\lambda>1$ be the constant given by Proposition~\ref{prop:wojtkowski_planar_cones}. For $\mu$-a.e.\ $w\in\Sigma$, let $\omega=M_+^{-1}(w)$, then the lengths of $\gamma^u(\omega)$ and $\Theta_{\Sigma}^n(\gamma^u(\omega))$ are defined and
\beql
\label{eqn:planar_manifold_growth}
\text{length}\left(\Theta_{\Sigma}^n(\gamma^u(\omega))\right)\geqs\lambda^n\text{length}\left(\gamma^u(\omega)\right).
\eeql
\end{thm}

\subsection{Proof of Theorem~\ref{thm:growth_of_unstable_manifolds}}

Our discussion of the length of piecewise smooth curves motivates us to formulate the following proposition:

\begin{prop}
\label{prop:manifolds_intersect_bad_set}
For any $n\in\mathbb{N}$, for $\mu$-a.e.\ $w\in\Sigma$ and for $\omega=M_+^{-1}(w)$, the set of $\omega'\in\gamma^u(\omega)$ at which $D\Theta^n_{\Sigma}$ does not exist is at most countable.
\end{prop}

In proving Proposition~\ref{prop:manifolds_intersect_bad_set} we will use two lemmas. We state and prove these lemmas first, then prove the proposition and last of all prove the theorem. We recommend that the reader skip forward to the proof of the theorem and refer back to the proposition and lemmas in that order. 

\begin{figure}[htp]
\centering
\includegraphics[totalheight=0.2\textheight]{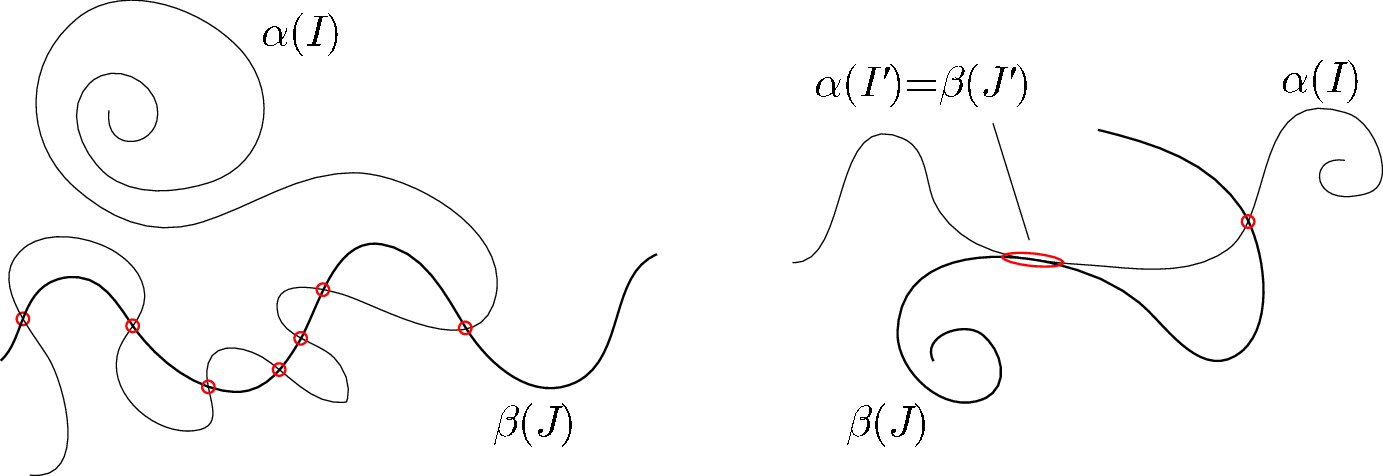}
\caption[Intersection of continuous curves]{Illustration of Lemma~\ref{lem:intersection_of_curves}: when the traces of continuous, injective curves $\alpha:I\to M$ and $\beta:J\to M$ intersect, the set of intersections is either at most countable, as in the left-hand figure, or it contains the images of open intervals $I'\subset I$ and $J'\subset J$ such that $\alpha(I')=\beta(J')$, as shown in the right-hand figure.}
\label{fig:curves}
\end{figure}

Our first lemma is easily stated and proven but we do not know of it anywhere in the literature. Loosely speaking, it says that when the traces of two continuous, injective curves intersect there are only the two possibilities shown in Figure~\ref{fig:curves}.

\begin{lem}
\label{lem:intersection_of_curves}
Let $I,J\subset\mathbb{R}$ be connected open intervals, let $M$ be a metric space and let $\alpha:I\to M$ and $\beta:J\to M$ each be continuous curves of finite length. Either there are open intervals $I'\subset I$ and $J'\subset J$ for which $\alpha(I')=\beta(J')$, or the set $\alpha(I)\cap\beta(J)$ is at most countable.
\end{lem}

\begin{proof}
Assume that there are no intervals $I'$ and $J'$ having the property stated, and let $t_0,t_1\in I$ be any two distinct points such that $\alpha(t_0)\cap\beta(J)\neq\emptyset$ and $\alpha(t_1)\cap\beta(J)\neq\emptyset$. (If there do not exist such points $t_0,t_1$ then, of course, we are done.) By the `no intervals' assumption, there must be some $t\in(t_0,t_1)$ for which $\alpha(t)\cap\beta(J)=\emptyset$. Set
\beq*
d(\alpha(t),\beta(J))=\eps>0,
\eeq*
where $d$ denotes the metric on $M$, and where the distance from the point $\alpha(t)$ to the set $\beta(J)$ is defined, in the usual way, as the infimum over distances $d(\alpha(t),p)$ for $p\in\beta(J)$. 

Now, $d(\alpha(\cdot),\beta(J))$ is a continuous function of $s\in I$, so there is some $\delta>0$ and an open neighbourhood $(t-\delta,t+\delta)\subset(t_0,t_1)\subset I$ such that $(t-\delta,t+\delta)\cap\beta(J)=\emptyset$. By this argument, the number of intersections $\alpha(I)\cap\beta(J)$ can exceed the number of such open neighbourhoods by at most one and there can be at most countably many disjoint open intervals in the interval $I$.
\end{proof}

Our second lemma concerns the nature of the set of points at which $D\Theta^n_{\Sigma}$ is not differentiable.

\begin{lem}
\label{lem:nature_sing_set}
Let $n\in\mathbb{N}$. The set of points in $M_+^{-1}(A)$ at which $D\Theta^n_{\Sigma}$ is not differentiable is given by $M_+^{-1}(s_i\cup\{(-1,0)\})$, where
\beql
\label{eqn:non_diff_pts_pltm}
s_i=\bigcup_{h=0}^{i-1}\Theta^{-h}\left(\pd A_+\cup\Phi^{-1}(\pd A_-)\right),
\eeql
for some $i\in\mathbb{N}$. Moreover $M_+^{-1}(s_i)$ is the union of a single point with a finite union $\bigcup_{h=0}^N\beta_h(J_h)$, where $N\in\mathbb{N}$ and for each $h=1,2,...,N$ the set $J_h\subset\mathbb{R}$ is an open interval and $\beta_h:J_h\to M_+^{-1}(A)$ is a continuous, injective curve of finite length.
\end{lem}

\begin{figure}[htp]
\centering
\includegraphics[totalheight=0.3\textheight]{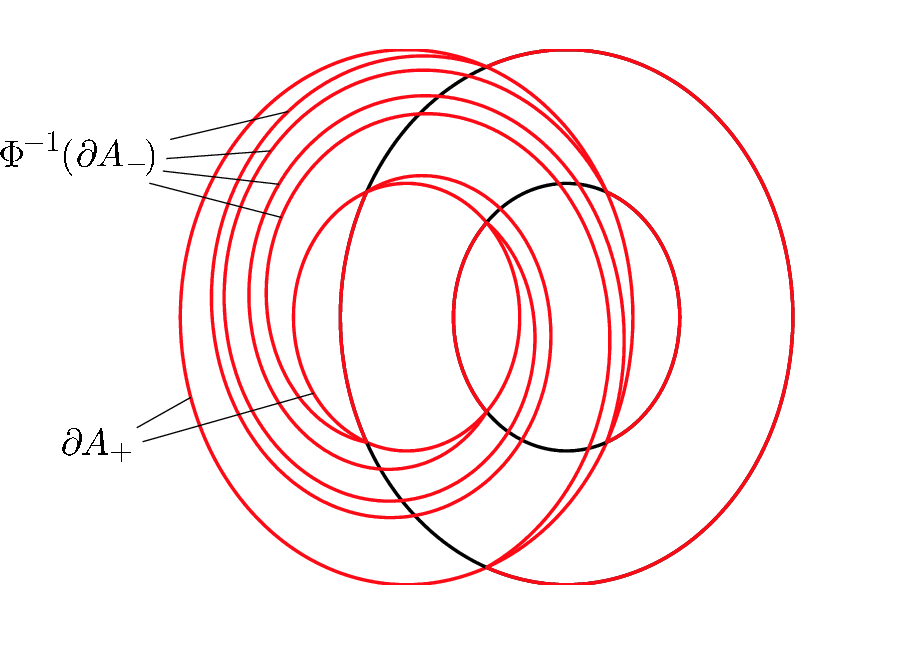}
\caption[Non-differentiable points for $\Theta$]{The manifold $A$ showing those points at which $\Theta$ is non-differentiable. Red lines indicate the set $s_0$ of non-differentiable points, black lines indicate other boundary points.}
\label{fig:singularities}
\end{figure}

In Figure~\ref{fig:singularities} we illustrate the set $s_0$ of non-differentiable points.

\begin{proof}
We prove the statement in the case $n=1$, the proof of the general case is no more difficult. Recall that $D\Theta_{\Sigma}:M_+^{-1}(\Sigma)\to M_+^{-1}(\Sigma)$ is given by $M_+^{-1}\circ\Theta^{i}\circ M_+$ for some positive integer $i$. We investigate whether the derivative
\beql
\label{eqn:deriv_return_map}
D(\Theta_{\Sigma})_{\omega}=D(M_+^{-1})_{\Theta^{i}(M_+(\omega))}D\Theta^{i}_{M_+(\omega)}D(M_+)_{\omega}
\eeql
exists for a given $\omega\in M_+^{-1}(A)$. For an affirmative answer we require that each of the three derivatives on the right-hand side of (\ref{eqn:deriv_return_map}) exist. First note that $M_+$ is differentiable on $A$.

$\Theta=\Gamma\circ\Phi$ where $\Phi$ is differentiable except on $\pd A_+$ and $\Gamma$ is differentiable except on $\pd A_-$. Thus we require that $M_+(\omega)\notin\pd A_+$ and that $\Phi(M_+(\omega))\notin\pd A_-$, i.e.\ that $M_+(\omega)\notin s_0=\left(\pd A_+\cup\Phi^{-1}(\pd A_-)\right)$.

$\Theta^2$ is differentiable if $M_+(\omega)\notin s_0$ and $\Theta(M_+(\omega))\notin s_0$, i.e.\ if $M_+(\omega)\notin s_0\cup\Theta^{-1}(s_0)$. Continuing inductively we see that $\Theta^{i}$ is differentiable except at those $M_+(\omega)\in\bigcup_{h=0}^{i-1}\Theta^{-h}(s_0)$.

Lastly we observe that $M_+^{-1}:A\to\mathbb{R}_0^+\times\mathbb{S}^1$ is differentiable except at the point $(-1,0)$. Combining these conditions we conclude that the set of points in $M_+^{-1}(A)$ at which $D\Theta^n_{\Sigma}$ is not differentiable is given by $M_+^{-1}(s_i\cup\{(-1,0)\})$ for some $i\in\mathbb{N}$, with $s_i$ as defined by (\ref{eqn:non_diff_pts_pltm}).

To see that the second statement holds, observe that $\pd A_-$ (a pair of circles) is the trace of two smooth curves. Applying $\Phi^{-1}$ will leave parts of these circles invariant and skew those parts which cross $A_+$ around that annulus. See Figure~\ref{fig:singularities}. This process stretches each circle (by a finite amount) and means that the curve whose trace is $\Phi^{-1}(A_+)$ is not differentiable at the points which map to $\pd A_+$. There are only finitely many such points however, so this curve is a finite union of smooth pieces. We proceed inductively, observing that each of the maps $\Phi^{-1},\Theta^{-h}$ and $M_+^{-1}$ will introduce only a finite number of non-differentiable points at each step, to obtain the result.
\end{proof}

We are now in a position to prove the proposition.

\begin{proof}[Proof of Proposition~\ref{prop:manifolds_intersect_bad_set}]
Let $i\in\mathbb{N}$ and let $X_i\subset M_+^{-1}(\Sigma)$ contain precisely those $\omega$ for which $\gamma^u(\omega)\cap M_+^{-1}(s_i)$ is uncountable, where $s_i$ is defined by (\ref{eqn:non_diff_pts_pltm}). The set $M_+^{-1}(s_i)$ consists of those points at which $D\Theta^n_{\Sigma}$ is non-differentiable, so $X_i$ consists of those points $\omega\in M_+^{-1}(\Sigma)$ for which we might not be able to determine the length of $\Theta^n_{\Sigma}(\gamma(\omega))$. We will show that $\mu(M_+(X_i))=0$.

Lemma~\ref{lem:nature_sing_set} says that if $\gamma^u(\omega)$ exists for some $\omega\in M_+^{-1}(\Sigma)$, then the set $\gamma^u(\omega)\cap M_+^{-1}(s_i)$ is a finite union
\beql
\label{eqn:non_differentiable_set}
\bigcup_{h=1}^N\alpha(I)\cap\beta_h(J_h),
\eeql
where $I$ and each $J_h$ are connected open intervals of $\mathbb{R}$ and $\alpha:I\to M_+^{-1}(A)$ and each $\beta_h:J_h\to M_+^{-1}(A)$ are continuous curves of finite length. Now assume that the set (\ref{eqn:non_differentiable_set}) is uncountable, i.e.\ $\omega\in X_i$. Then there is at least one $h\in{1,2,...,N}$ for which $\alpha(I)\cap\beta_h(J_h)$ is uncountable. Fix such an $h$.

Lemma~\ref{lem:intersection_of_curves} says that there are intervals $I'\subset I$ and $J_h'\subset J_h$ and that $\alpha(I')=\beta_h(J_h')$. This is significant because it says that $\gamma^u(\omega)$ coincides with $M_+^{-1}(s_i)$ on a subset of positive length. The total length of $M_+^{-1}(s_i)$ is well-defined and finite, so if $\{\omega_l\}_{l\in L}$ is any collection of points $\omega_l\in X_i$ such that $\gamma^u(\omega_{l_1})\cap\gamma^u(\omega_{l_2})=\emptyset$ for $l_1\neq l_2$, then $L$ is at most a countable set.

Our strategy is to find such a countable set that covers $X_i$. To that end we recall the familiar global unstable manifold of $\omega\in X_i$ given by
\beq*
W^u(\omega)=\bigcup_{l=1}^{\infty}\Theta^l\left(\gamma^u(\Theta^{-l}(\omega))\right).
\eeq*
(See for example \citet{bap}.) It is easily shown that $\omega\in W^u(\omega')$ if and only if $\omega'\in W^u(\omega)$ and that otherwise we have $W^u(\omega)\cap W^u(\omega')=\emptyset$. Thus for $\omega\in X_i$ there is a well-defined equivalence class
\beq*
[\omega]=\left\{\omega'\in X_i:\omega'\in W^u(\omega)\right\}.
\eeq*
We can clearly cover $X_i$ for \emph{some} set of distinct elements $\{[\omega_l]\}_{l\in L}$ and by our previous reasoning any such $L$ is at most countable. This is sufficient to give the result:
\beq*
\mu(M_+(X_i))\leqs\mu\left(M_+\left(\bigcup_{l=1}^{\infty}[\omega_l]\right)\right)\leqs\sum_{l=1}^{\infty}\mu(M_+([\omega_l]))=0.
\eeq*
\end{proof}

For any given $n\in\mathbb{N}$ we can use Proposition~\ref{prop:manifolds_intersect_bad_set} and our definition of the length of the trace of a smooth curve to define the length of $\Theta_{\Sigma}^n(\gamma^u(\omega))$ for a set of $\omega$ whose $M_+$-image has full $\mu$-measure in $\Sigma$. This enables us to show that this length grows exponentially with $n$, as in the following theorem.

\begin{proof}[Proof of Theorem~\ref{thm:growth_of_unstable_manifolds}]
Proposition~\ref{prop:manifolds_intersect_bad_set} says that for $\mu$-a.e.\ $w\in\Sigma$ and for $\omega= M_+^{-1}(w)$ we have $\Theta^n_{\Sigma}(\gamma^u(\omega))=\Theta^n_{\Sigma}\circ\alpha(I)$ where $I\subset\mathbb{R}$ is an open interval and $\alpha:I\to M_+^{-1}(A)$ is a smooth curve. Moreover there is a countable (at most) partition of $I$ as given by (\ref{eqn:partition_of_I}) and the restriction of $\Theta^n_{\Sigma}\circ\alpha$ to each open interval $(i_h,i_{h+1})$ is a smooth curve. The length of each smooth curve, by our earlier definition, is given by
\beq*
\text{length}(\Theta^n_{\Sigma}\circ\alpha|_{(i_h,i_{h+1})})=\int_{i_h}^{i_{h+1}}\left\|D\Theta^n_{\Sigma}D\alpha_t\right\|\d t,
\eeq*
and so the length of $\Theta^n_{\Sigma}\circ\alpha$ is given by
\beqa*
\text{length}(\Theta^n_{\Sigma}\circ\alpha) &=\sum_{h=0}^N\int_{i_h}^{i_{h+1}}\left\|D\Theta^n_{\Sigma}D\alpha_t\right\|\d t \\
&\geqs \sum_{h=0}^N\int_{i_h}^{i_{h+1}}\lambda^n\left\|D\alpha_t\right\|\d t \\
&=\lambda^n\sum_{h=0}^N\int_{i_h}^{i_{h+1}}\left\|D\alpha_t\right\|\d t \\
&=\lambda^n\int_I\left\|D\alpha_t\right\|\d t \\
&=\lambda^n\text{length}(\alpha),
\eeqa*
where $N\in\mathbb{N}\cup\{\infty\}$. Here the second line holds because of Proposition~\ref{prop:wojtkowski_planar_cones} and the fact that $D\alpha_t\in U\subset T_{\alpha(t)}L$ for each $t\in I$.
\end{proof}

%% file: chapter3/section4.tex
\section{A new invariant cone for $\Theta$}\label{PLTM.Cone_C}
\setcounter{equation}{0}

In this section we express the map $\Theta:A\to A$ in the new coordinates developed in Section~\ref{PLTM.New_Coords}. Once we have done this we will introduce a new cone in the tangent space, which is invariant under the map. The purpose of this construction is to give an improved estimate on the orientation of local unstable manifolds. This, combined with the growth already established will enable us to show that condition (\ref{eqn:mip_pltm}) is satisfied. We will do that in the following, final section.

\subsection{The map $\Theta$ expressed in the new coordinates}

\begin{figure}[htp]
\centering
(a)\includegraphics[totalheight=0.27\textheight]{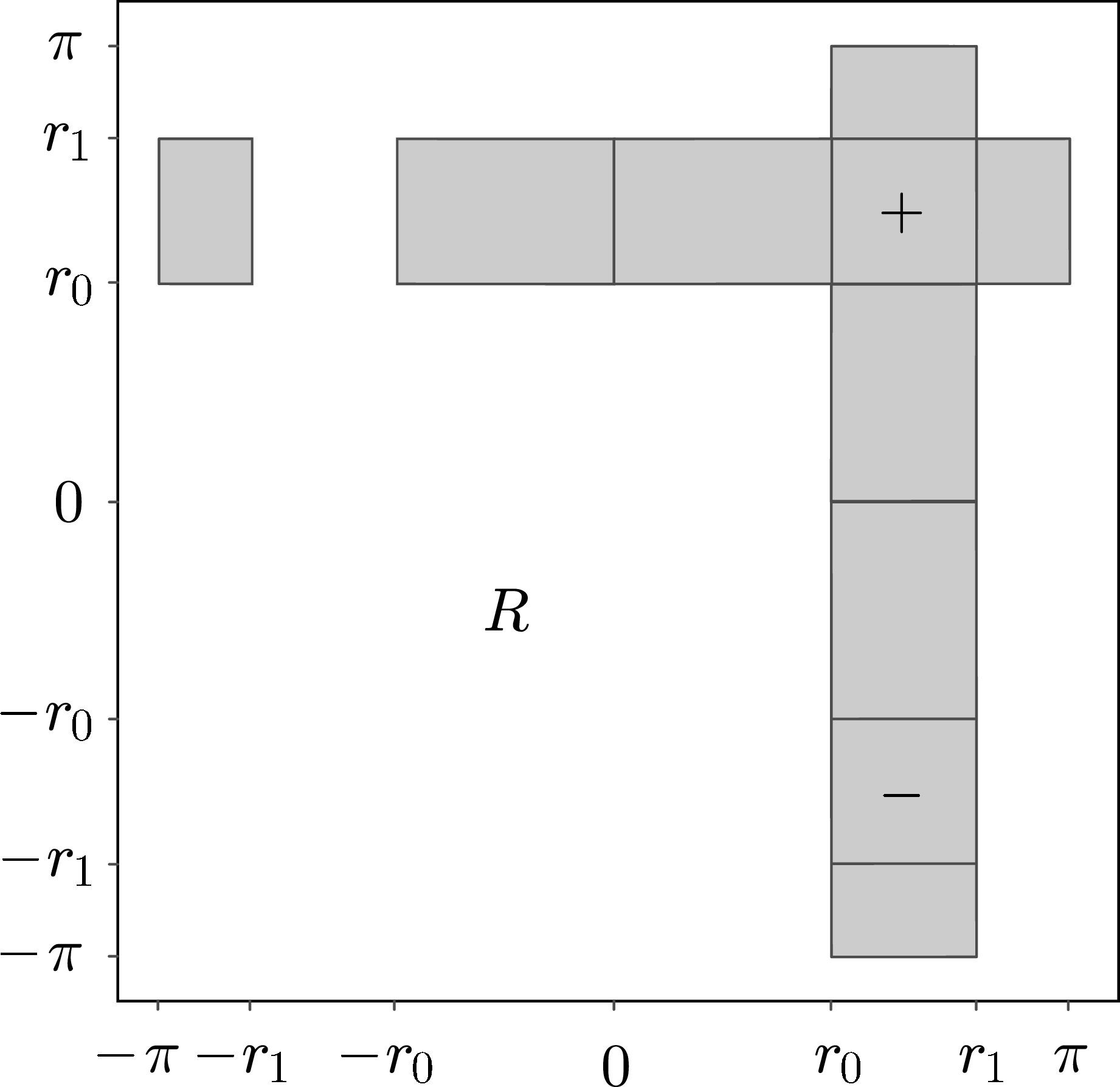}%\quad
(b)\includegraphics[totalheight=0.27\textheight]{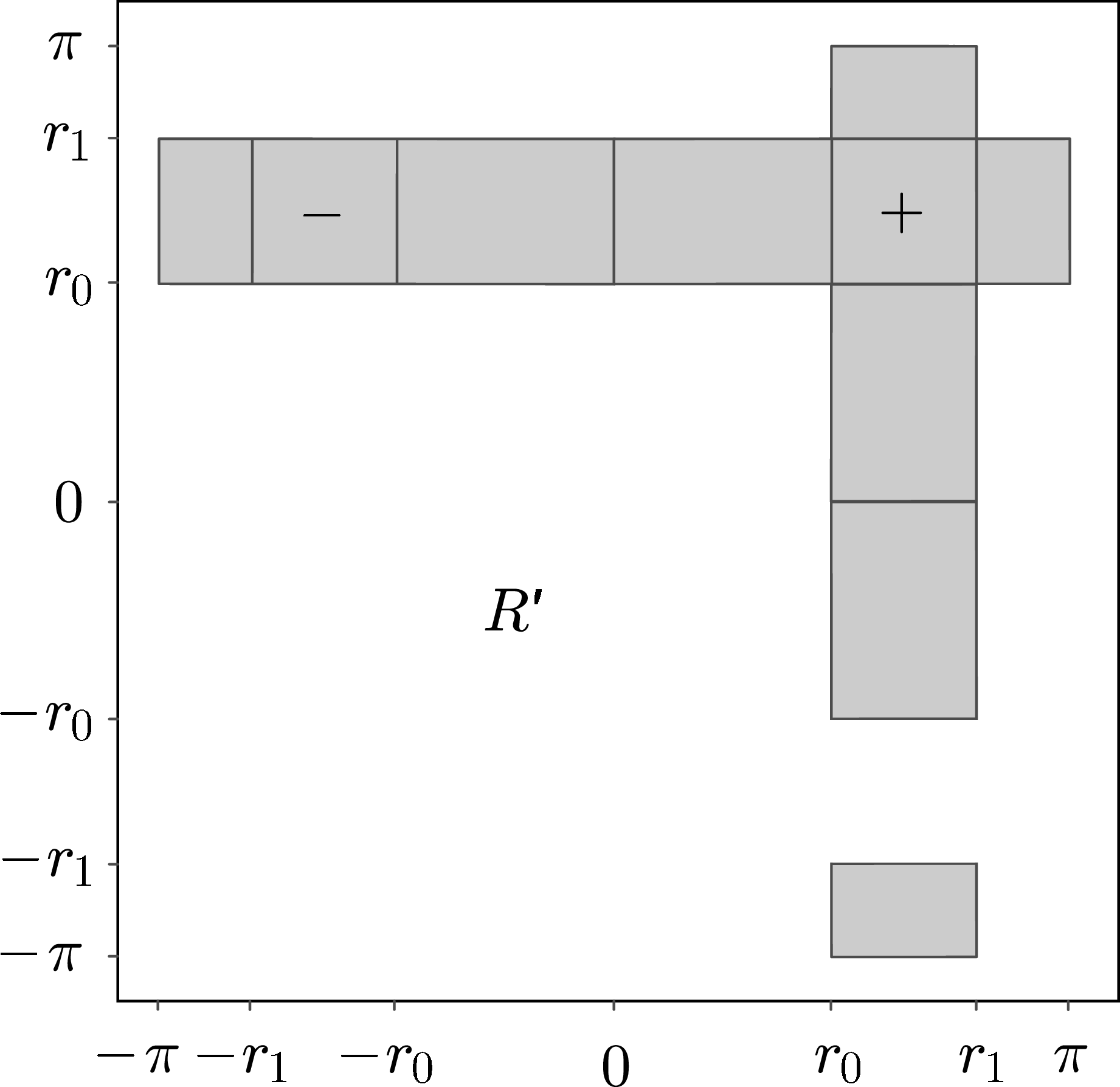}
\caption[The manifolds $R$ and $R'$]{The manifolds $R,R'\subset\mathbb{T}^2=\mathbb{S}^1\times\mathbb{S}^1$, in parts~(a) and~(b) respectively . Each is in one-to-one correspondence with $A$, the difference between the two being that $\Sigma_-$ is represented differently in each. The manifolds have the property that $F$ is a homeomorphism of $R$ and $G$ is a homeomorphism of $R'$.}
\label{fig:R_and_R'}
\end{figure}

We will use notation reminiscent of that we have used for other linked-twist maps. We begin by giving some notation for the manifold $A$ when transformed into the new coordinates. Let
\beq*
R=\left\{\Psi\circ M_+^{-1}(A_+)\right\}\cup\left\{\iota^{-1}\circ\Psi\circ M_-^{-1}(A_-\backslash\Sigma_-)\right\}.
\eeq*
The set $R\subset\mathbb{T}^2$ is shown in Figure~\ref{fig:R_and_R'}(a). There is a one-to-one correspondence between points in $R$ and points in $A$. Let $F:R\to R$ denote the map $\Phi:A_+\to A_+$ in the new coordinates. Recall that $\Phi=M_+\circ\Lambda\circ M_+^{-1}$, so
\beqa*
F &= \Psi\circ M_+^{-1}\circ\Phi\circ M_+\circ\Psi^{-1} \\
	&= \Psi\circ\Lambda\circ\Psi^{-1}.
\eeqa*
We observe that $F$ is a homeomorphism of $R$.

Now let
\beq*
R'=\left\{\Psi\circ M_+^{-1}(A_+\backslash\Sigma_-)\right\}\cup\left\{\iota^{-1}\circ\Psi\circ M_-^{-1}(A_-)\right\}.
\eeq*
This is illustrated in Figure~\ref{fig:R_and_R'}(b). Again, there is a one-to-one correspondence between points in $R'$ and points in $A$. We let $G:R'\to R'$ denote the map $\Gamma:A_-\to A_-$ in the new coordinates, giving
\beqa*
G &= \iota^{-1}\circ\Psi\circ\Lambda^{-1}\circ\Psi^{-1}\circ\iota \\
	&= \iota^{-1}\circ F^{-1}\circ\iota.
\eeqa*
We observe that $G$ is a homeomorphism of $R'$.

Now consider the linked-twist map $\Theta$ expressed in the new coordinates, that is the composition $H=G\circ F$, which we wish to represent as a map of $R$ into itself. $H$ will be a homeomorphism as it is the composition of two homeomorphisms. Recall that when we defined the map $\Theta=\Gamma\circ\Phi$ we changed coordinates between the two mappings; here we do likewise. The difference is that this is now only necessary when the trajectory lands in (the new representation of) $\Sigma_-$. Moreover the coordinate transformation is given simply by $(x,y)\mapsto(-x,-y)$.

Let us make this algorithm more explicit; recall our notation that $\mathcal{R}=[r_0,r_1]$ and $-\mathcal{R}=[-r_1,-r_0]$. Let $(x,y)\in R$, and first apply the map $F$. If and only if $F(x,y)\in\mathcal{R}\times-\mathcal{R}$, reflect both coordinates through the origin. We are now in $R'$. Apply $G$, i.e.\ apply $\iota$, followed by $F^{-1}$, followed by $\iota^{-1}$. Finally, if and only if we are now in $-\mathcal{R}\times\mathcal{R}$, reflect coordinates through the origin again.

Rather than include expressions both for reflections through the origin and the functions $\iota^{\pm}$ we are able to combine the two. We have $H:R\to R$ given by
\beq*
H=\Omega^{-1}\circ F^{-1}\circ\Omega\circ F,
\eeq*
where the function $\Omega$ is given by
\beq*
\Omega^{\pm 1}(x,y)=\left\{
\begin{array}{r@{\quad}l}
\iota^{\mp 1} &\text{if }(x,y)\in\mathcal{\mathcal{R}}\times\mp\mathcal{R} \\
\iota^{\pm 1} &\text{otherwise.}
\end{array}
\right.
\eeq*
Let $S\subset R$ be the image of the `intersection region' $\Sigma$, i.e.\ $S=\left(\mathcal{R}\times\mathcal{R}\right) \cup \left(\mathcal{R}\times-\mathcal{R}\right)$. For $z\in S$ we define the return map $H_S:S\to S$, defined completely analogously to the return map $\Theta_{\Sigma}$.

We now study the derivative of $H$. Let $D_1,D_2$ denote the usual differential operators. We have
\beqa*
F^{\pm 1}(x,y)	&=\Psi\circ\Lambda^{\pm 1}\circ\Psi^{-1}(x,y) \\
								&=\Psi\circ\Lambda^{\pm 1}\left(x,\psi^{-1}(x,y)\right) \\
								&=\Psi\left(x,\psi^{-1}(x,y)\pm c(x-r_0)\right) \\
								&=\left(x,\psi\left(x,\psi^{-1}(x,y)\pm c(x-r_0)\right)\right).
\eeqa*
To simplify the expression we define
\beq*
f_{\pm}(x,y)=\psi\left(x,\psi^{-1}(x,y)\pm c(x-r_0)\right),
\eeq*
then the Jacobians of $F^{\pm 1}$ are given by
\beq*
DF^{\pm 1}=\left( \begin{array}{cc}
1 & 0 \\ D_1f_{\pm}(x,y) & D_2f_{\pm}(x,y)
\end{array}\right).
\eeq*

\begin{figure}[htp]
\centering
\includegraphics[totalheight=0.14\textheight]{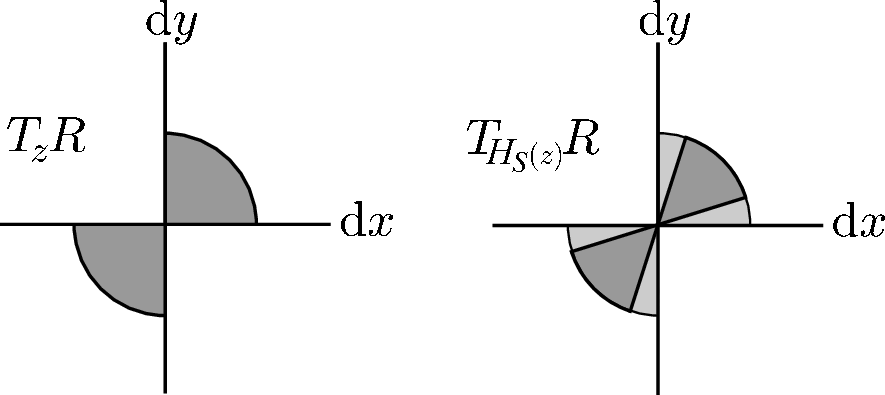}
\caption[The cone $C$]{The invariant cone $C\subset T_zR$ is shown in the left-hand figure. In the right-hand figure is the image of the cone under the differential map $DH_S$. The fact that $C$ is invariant under this differential is immediately implied by Proposition~\ref{prop:DGDF}. Notice that we do \emph{not} claim that this cone is expanded by $DH_S$. }
\label{fig:C}
\end{figure}

We introduce further notation $\tilde{y}_{\pm}=\tilde{y}_{\pm}(x,y)=\psi^{-1}(x,y)\pm c(x-r_0)$, then
\beqal
\label{eqn:D1f}
D_1f_{\pm}(x,y) &= D_1\psi(x,\tilde{y}_{\pm}) + D_2\psi(x,\tilde{y}_{\pm})\left[D_1\psi^{-1}(x,y)\pm c\right], \\
\label{eqn:D2f}
D_2f_{\pm}(x,y) &= D_2\psi(x,\tilde{y}_{\pm}) D_2\psi^{-1}(x,y).
\eeqal
We use these derivatives to prove a result for $DH$. Let $b_1=\text{d}x,b_2=\text{d}y$ give coordinates in the tangent space $T_z\mathbb{T}^2$ to a point $z=(x,y)\in R$, and define the cones
\beq*
C(z)=\left\{(b_1,b_2):b_1b_2>0\right\},\qquad\tilde{C}(z)=\left\{(b_1,b_2):b_1b_2<0\right\}
\eeq*
The cone $C$ is illustrated in Figure~\ref{fig:C}. We define the cone fields
\beq*
C_+=\bigcup_{z\in R}C(z),\qquad C_-=\bigcup_{z\in R}\tilde{C}(z).
\eeq*

\begin{prop}
\label{prop:DGDF}
Let $r_0=2$ and $r_1=\sqrt{7}$. Then $DF$ and $D(\Omega^{-1}\circ F^{-1}\circ\Omega)$ preserve the cone field $C_+$.
\end{prop}

The remainder of this section is devoted to proving Proposition~\ref{prop:DGDF}. We remark that we have had to give explicit sizes for the annuli and so our result is not as general as one might hope for. We discuss this more in Chapter~\ref{concl}.

\subsection{Proof of Proposition~\ref{prop:DGDF}}

Both $D\Omega$ and $D\Omega^{-1}$ map the cone $C$ into the cone $\tilde{C}$ and \emph{vice versa}. Consequently it suffices to show that $DF$ preserves $C_+$ and $DF^{-1}$ preserves $C_-$. Let $z=(x,y)\in A$, let $(b_1,b_2)\in T_z\mathbb{T}^2$ and define
\beq*
\left( \begin{array}{cc} b_1' \\ b_2'  \end{array}\right)
= DF^{\pm 1} \left( \begin{array}{cc} b_1 \\ b_2  \end{array}\right)
= \left( \begin{array}{cc} b_1 \\ b_1D_1f_{\pm}(x,y)+b_2D_2f_{\pm}(x,y) \end{array}\right).
\eeq*
We have
\beq*
\frac{b_2'}{b_1'}=D_1f_{\pm}(x,y)+\frac{b_2}{b_1}D_2f_{\pm}(x,y).
\eeq*
So it is enough to show that for every $(x,y)\in \mathcal{R}\times\mathbb{S}^1$ we have
\beql
\label{eqn:Df_bounds}
\pm D_1f_{\pm}(x,y)>0\quad\text{and}\quad D_2f_{\pm}(x,y)>0.
\eeql
In light of (\ref{eqn:D1f}) and (\ref{eqn:D2f}), the condition (\ref{eqn:Df_bounds}) will follow from suitable bounds on the derivatives of $\psi$ and its inverse. Determining such bounds is our strategy for this proof; unfortunately, some extensive calculations will be unavoidable. 

Recall the real-valued function $\tau$, defined on $\Psi\circ M_+^{-1}(\Sigma)=\mathcal{R}\times \mathcal{R}$ by
\beq*
\tau(r,\theta)=\frac{1}{4r}\left(r^2-\theta^2+4\right).
\eeq*
It will play a significant role. Partial derivatives of $\tau$ are given by
\beq*
\frac{\pd}{\pd r}\tau(r,\theta)=\frac{1}{4r^2}\left(r^2+\theta^2-4\right)\quad\text{and}\quad\frac{\pd}{\pd \theta}\tau(r,\theta)=-\frac{\theta}{2r},
\eeq*
so $\tau$ is an increasing function of $r$ and a decreasing function of $\theta$. This observation is enough to prove some bounds on certain functions of $\tau$, which we will use in the rest of this section. We collect these in the following lemma, though we omit a formal proof.
\begin{lem}[Properties of the function $\tau$]
\label{lem:prop_tau}
Let $\tau$ denote $\tau(r,r_0)$. Then
\beq*
\tau\in\left[\frac{1}{2},\frac{\sqrt{7}}{4}\right],\quad
\frac{1}{2}-\frac{\tau}{r}=\frac{1}{4},\quad
\sqrt{1-\tau^2}\in\left[\frac{3}{4},\frac{\sqrt{3}}{2}\right],\quad
\cos^{-1}\tau\in\left(\frac{5}{6},\frac{\pi}{3}\right].
\eeq*
Alternatively, let $\tau$ denote $\tau(r,r_1)$. Then
\beq*
\tau\in\left[\frac{1}{8},\frac{\sqrt{7}}{7}\right],\quad
\frac{1}{2}-\frac{\tau}{r}\in\left(\frac{2}{3},\frac{3}{4}\right),\quad
\sqrt{1-\tau^2}\in\left(\frac{11}{12},1\right),\quad
\cos^{-1}\tau\in\left(\frac{7}{6},\frac{3}{2}\right).
\eeq*
\end{lem}

It will be convenient to denote by $A_+^i$ that part of $A_+$ which lies \emph{inside} the annulus $A_-$ and by $A_+^o$ that part of $A_+$ lying \emph{outside} of $A_-$. Moreover it will be necessary to determine the values $\theta$ may take on each part. On $A_+^i$, for a given $r\in \mathcal{R}$, this range is given by $[0,\cos^{-1}\tau]$, where $\tau$ denotes $\tau(r,r_0)$. Using Lemma~\ref{lem:prop_tau} we calculate 
\beq*
\sup_{(r,\theta)\in A_+^i}\theta =\sup_{r\in \mathcal{R}}\cos^{-1}\tau =\cos^{-1}\inf_{r\in \mathcal{R}}\tau =\cos^{-1}\frac{1}{2} =\frac{\pi}{3}.
\eeq*
Similarly, for a given $r\in \mathcal{R}$ the range of $\theta$ such that $(r,\theta)\in A_+^o$ is given by $[\cos^{-1}\tau,\pi]$, where $\tau=\tau(r,r_1)$. We have
\beq*
\inf_{(r,\theta)\in A_+^o}\theta =\inf_{r\in \mathcal{R}}\cos^{-1}\tau =\cos^{-1}\sup_{r\in \mathcal{R}}\tau =\cos^{-1}\frac{\sqrt{7}}{7} >\frac{7}{6}.
\eeq*

We now prove some bounds on the derivatives of $\psi$.

\begin{lem}
\label{lem:D1psi}
$D_1\psi\in\left[0,\frac{7}{6}\right)$.
\end{lem}
\begin{proof}
We consider the three cases separately. Let $(r,\theta)\in A_+^i$, then
\beq*
D_1\psi(r,\theta)=\frac{r_0\theta\left(\frac{1}{2}-\frac{\tau}{r}\right)}{\left(\cos^{-1}\tau\right)^2\sqrt{1-\tau^2}},
\eeq*
where $\tau=\tau(r,r_0)$. The numerator equates to $\theta/2$, so has range $[0,\pi/6]$. Using Lemma~\ref{lem:prop_tau}, the denominator lies in $[25/54,\pi^2\sqrt{3}/18]$. Thus on $A_+^i$ we have $D_1\psi\in[0,9\pi/25]\subset[0,7/6)$.

Next let $(r,\theta)\in A_+^o$, then
\beq*
D_1\psi(r,\theta)=\frac{\left(\frac{1}{2}-\frac{\tau}{r}\right)(\pi-r_1)(\pi-\theta)}{(\pi-\cos^{-1}\tau)^2\sqrt{1-\tau^2}},
\eeq*
where $\tau=\tau(r,r_1)$. A lower bound of $0$ is attained when $\theta=\pi$; using Lemma~\ref{lem:prop_tau} we see calculate that $3/10$ is an upper bound, so here $D_1\psi\in[0,3/10)$.

Lastly let $(r,\theta)\in\Sigma$, so that
\beq*
D_1\psi(r,\theta)=\frac{r-2\cos\theta}{\psi}.
\eeq*
where $\psi=\psi(r,\theta)$ is as in (\ref{eqn:psi1}). On $\Sigma_+$ recall that $\cos\theta=\tau(r,\psi)$, so the numerator is given by
\beq*
r-\frac{1}{2r}\left(r^2-\psi^2+4\right)=\frac{r}{2}-\frac{2}{r}+\frac{\psi^2}{2r}.
\eeq*
By construction, $\psi(r,\theta)\in[r_0,r_1]$ for $(r,\theta)\in\Sigma_+$. Also $0\leqslant r/2-2/r\leqslant3\sqrt{7}/14$, all of which means that $D_1\psi\in[\sqrt{7}/7,5\sqrt{7}/14]$.
\end{proof}

\begin{lem}
\label{lem:D2psi}
$D_2\psi\in\left[\frac{1}{4},\sqrt{7}\right]$.
\end{lem}
\begin{proof}
As before there are three cases. If $(r,\theta)\in A_+^i$ then
\beq*
D_2\psi(r,\theta)=r_0/\cos^{-1}\tau.
\eeq*
where $\tau=\tau(r,r_0)$. Using Lemma~\ref{lem:prop_tau} we easily deduce that $D_2\psi\in[4/3,12/7]$.

If $(r,\theta)\in A_+^o$ then
\beq*
D_2\psi(r,\theta)=(\pi-r_1)/(\pi-\cos^{-1}\tau).
\eeq*
where $\tau=\tau(r,r_1)$. Lemma~\ref{lem:prop_tau} allows us to deduce the bounds $D_2\psi\in(1/4,1/3)$.

Finally if $(r,\theta)\in\Sigma_+$ then
\beq*
D_2\psi(r,\theta)=\frac{2r\sin\theta}{\psi(r,\theta)}.
\eeq*
A lower bound for $\sin\theta$ is calculated as follows:
\beq*
\inf_{(r,\theta)\in\Sigma_+}\sin\theta
=\inf_{r\in \mathcal{R}}\sin\left(\cos^{-1}\tau(r,r_0)\right)
=\inf_{r\in \mathcal{R}}\sqrt{1-[\tau(r,r_0)]^2}
=\frac{3}{4}.
\eeq*
Trivially, we also have upper bound $1$. From here we easily obtain the bounds $D_2\psi\in(3\sqrt{7}/7,\sqrt{7})$.
\end{proof}

It is easily verified that the inverse to $\psi$ on $\Sigma$ (where $\psi$ is given by (\ref{eqn:psi1})) is given by
\beql
\label{eqn:psi_inv1}
\psi^{-1}(x,y)=\cos^{-1}\left(\tau(x,y)\right).
\eeql
The inverse to $\psi$ on $\mathcal{R}\times [0,r_0]$ (the inverse to (\ref{eqn:psi2})) is given by
\beql
\label{eqn:psi_inv2}
\psi^{-1}(x,y)=y\cos^{-1}(\tau(x,r_0))/r_0
\eeql
and the inverse to $\psi$ on $\mathcal{R}\times [r_1,\pi]$ (the inverse to (\ref{eqn:psi3})) is given by
\beql
\label{eqn:psi_inv3}
\psi^{-1}(x,y)=\cos^{-1}(\tau(x,r_1))+\frac{y-r_1}{\pi-r_1}(\pi-\cos^{-1}(\tau(x,r_1))).
\eeql
We will now prove some bounds on the derivatives of $\psi^{-1}$.

\begin{lem}
\label{lem:D1psi_inv}
$D_1\psi^{-1}\in\left[-\frac{9}{11},0\right]$.
\end{lem}
\begin{proof}
Let $(x,y)\in \mathcal{R}\times[0,r_0]$, so $\psi^{-1}$ is given by (\ref{eqn:psi_inv2}). Then
\beq*
D_1\psi^{-1}(x,y)=\frac{y(\frac{\tau}{x}-\frac{1}{2})}{r_0\sqrt{1-\tau^2}},
\eeq*
where $\tau=\tau(x,r_0)$. Using Lemma~\ref{lem:prop_tau} the numerator equates to $-y/4$, having range $[-1/2,0]$. Similarly the denominator has range $[3/2,\sqrt{3}]$, giving $D_1\psi^{-1}\in[-1/3,0]$.

For $(x,y)\in \mathcal{R}\times[r_0,\pi]$ the function $\psi^{-1}$ is given by (\ref{eqn:psi_inv3}) and so
\beq*
D_1\psi^{-1}(x,y)=\frac{(\frac{1}{2}-\frac{\tau}{x})(y-\pi)}{(\pi-r_1)\sqrt{1-\tau^2}},
\eeq*
where $\tau=\tau(x,r_1)$. We find that the numerator is in the range $[3(\sqrt{7}-\pi)/4,0]$, whereas the denominator is in $[11(\pi-\sqrt{7})/12,\pi-\sqrt{7}]$. Thus $D_1\psi^{-1}\in[-9/11,0]$.

Lastly if $(x,y)\in \mathcal{R}\times \mathcal{R}$ then $\psi^{-1}$ is given by (\ref{eqn:psi_inv1}) and so
\beq*
D_1\psi^{-1}(x,y)=\frac{\frac{\tau}{x}-\frac{1}{2}}{\sqrt{1-\tau^2}},
\eeq*
where $\tau=\tau(x,y)$. The numerator equates to $-1/4+(4-y^2)/4x^2$, which is increasing in $x$ (because $y^2\geqslant 4$) and decreasing in $y$. It has the range $[-7/16,-1/4]$. Using Lemma~\ref{lem:prop_tau} and the fact that $\tau$ is increasing in $r$, the denominator is in $[3/4,1]$, giving $D_1\psi^{-1}\in[-7/12,-1/4]$.
\end{proof}

\begin{lem}
\label{lem:D2psi_inv}
$D_2\psi^{-1}\in\left[\frac{\sqrt{7}}{7},4\right]$.
\end{lem}
\begin{proof}
Let $(x,y)\in \mathcal{R}\times[0,r_0]$, then
\beq*
D_2\psi^{-1}(x,y)=\cos^{-1}\tau/r_0,
\eeq*
where $\tau=\tau(x,r_0)$. From the proof of Lemma~\ref{lem:D2psi} we have $D_2\psi^{-1}\in[7/12,3/4]$.

Next, let $(x,y)\in \mathcal{R}\times[r_1,\pi]$, then
\beq*
D_2\psi^{-1}(x,y)=\frac{\pi-\cos^{-1}\tau}{\pi-r_1},
\eeq*
with $\tau=\tau(x,r_1)$. Comparing with the proof of Lemma~\ref{lem:D2psi} we have $D_2\psi^{-1}\in[3,4]$.

Finally let $(x,y)\in\Sigma_+$, then
\beq*
D_2\psi^{-1}(x,y)=\frac{y}{2x\sqrt{1-\tau^2}}.
\eeq*
where $\tau=\tau(x,y)$. In the proof of Lemma~\ref{lem:D1psi_inv} we established that $\sqrt{1-\tau^2}\in[3/4,1]$ and so the denominator is in the range $[3,2\sqrt{7}]$. Thus we have $D_2\psi^{-1}\in[\sqrt{7}/7,\sqrt{7}/3]$.
\end{proof}

With the bounds established in Lemmas~\ref{lem:D1psi} through~\ref{lem:D2psi_inv} it is a simple matter to check that the conditions (\ref{eqn:Df_bounds}) will always be satisfied, and thus the proposition is proved.

%% file: chapter3/section5.tex
\section{The Bernoulli property}\label{PLTM.Bernoulli}
\setcounter{equation}{0}

In this section we improve our estimate on the direction of unstable manifolds in the following sense. Let $z=(x,y)\in R$. We show that $E^u(z)\subset C\subset T_zR$ and so $\gamma^u(z)$ is aligned within the cone $C$ (we will say precisely what we mean by this in a moment, although the reader probably has an intuitive idea). Using this fact we are able to deduce that the strong form of the manifold intersection property, condition~(\ref{eqn:mip_pltm}), is satisfied.

\subsection{Orientation of the unstable subspace}

Let $\omega\in M_+^{-1}(\Sigma)$ be `typical' in the sense of Section~\ref{PLTM.Manifolds}, which is to say that $\gamma^u(\omega)$ exists, its length is well defined and for any $n\in\mathbb{N}$ the length of $\Theta_{\Sigma}^n(\gamma^u(\omega))$ is well defined also. Recall that in this case the latter grows exponentially with $n$. Let $I$ be an interval and $\alpha:I\to M_+^{-1}(A)$ the smooth curve whose trace is $\gamma^u(\omega)$.

Let $z=\Psi(\omega)$ give such a point in the new coordinates, and consider $\gamma^u(z)=\Psi(\gamma^u(\omega))$.

\begin{defn}[Alignment of local invariant manifold]
We say that $\gamma^u(z)$ is aligned within the cone $C$, or that it has orientation in the cone $C$, if and only if the derivative $D\Psi D\alpha_0\in C\subset T_zR$.
\end{defn}

Suppose that $\gamma^u(z)$ is aligned within the cone $C$. It follows immediately from Proposition~\ref{prop:DGDF} that $\Theta_{\Sigma}^n(\gamma^u(\omega))$ is aligned within the cone $C$ also for any $n\in\mathbb{N}$. This is enough for us to prove the Bernoulli property as below.

Unfortunately, verifying that our assumption holds will require one last digression. For a full-measure set $w\in A$, denoting $\omega=M_+^{-1}(w)$ then the results of \citet{woj} tell us that $\gamma^u(\omega)$ has orientation in the cone $U$. It does \emph{not} necessarily follow from our coordinate transformation that if $z=\Psi(\omega)$ then $\gamma^u(z)$ has orientation in the cone $C$. To overcome this problem we formulate the following proposition:

\begin{prop}
\label{lem:unstable_direction}
Let $X$ be a compact metric space of dimension 2 and $T:X\to X$ a (non-uniformly) hyperbolic transformation, preserving a measure $\mu$ on $X$. Let $Y\subset X$ be a subset to which $\mu$-a.e.\ trajectory returns infinitely many times (i.e.\ $\{T^n(x)\}_{n\in\mathbb{N}}\cap Y$ is infinite) and for which the `first-return map' $T_Y:Y\to Y$ is uniformly hyperbolic.

Suppose that $T$ is differentiable $\mu$-a.e.\ and define cones $K(x)\subset T_xX$ of the form
\beq*
K(x)=\left\{(\eta,\zeta):k_1(x)<\zeta/\eta<k_2(x)\right\},
\eeq*
where $k_1,k_2\in[-\infty,\infty]$. If $DT_x\left(K(x)\right)\subset K(T(x))$ and the boundary of $K(x)$ is mapped to the interior of $K(T(x))$, then $E^u(x)$ must lie in the interior of $K(x)$.
\end{prop}

\begin{proof}
In the tangent space $T_xX$ define unit vectors in the stable and unstable subspaces
\beq*
\left(s_1(x),s_2(x)\right)\in E^s(x)\quad\text{and}\quad\left(u_1(x),u_2(x)\right)\in E^u(x)
\eeq*
respectively. Fix $x\in X$ and let $(\eta_0,\zeta_0)\in K(x)\subset T_xX$. There are unique (and, without loss of generality, non-negative) real constants $\alpha_0$ and $\beta_0$ such that
\beq*
(\eta_0,\zeta_0)=\alpha_0(s_1(x),s_2(x))+\beta_0(u_1(x),u_2(x)).
\eeq*
Similarly, if we fix $n\in\mathbb{N}$ and consider $(\eta_n,\zeta_n)=DT_x^n(\eta_0,\zeta_0)\in K(T^n(x))\subset T_{T^n(x)}X$ then there are unique and non-negative $\alpha_n,\beta_n$ such that
\beq*
(\eta_n,\zeta_n)=\alpha_n(s_1(T^n(x)),s_2(T^n(x))+\beta_n(u_1(T^n(x)),u_2(T^n(x))).
\eeq*
The uniform hyperbolicity of the return map to $Y\subset X$ allows us to estimate the magnitudes of $\alpha_n,\beta_n$. Given $m\in\mathbb{N}$ and $x\in Y$ there exists some $n\in\mathbb{N}$ such that $T_Y^m(x)=T^n(x)$. Under such circumstances we have
\beq*
\alpha_n\leqslant\lambda^{-m}\alpha_0\quad\text{and}\quad\beta_n\geqslant\lambda^m\beta_0,
\eeq*
where $\lambda>1$ is a constant independent of $m,n$ or $x$. It must be the case that $m\to\infty$ as $n\to\infty$ (though we say nothing about the relative rates of divergence) and so as $n\to\infty$ we see that $\alpha_n\to 0$ and $\beta_n\to\infty$.

$K(x)$ cannot be contained in $E^s(x)$, so there is $(\eta_0,\zeta_0)\in K(x)$ for which $\beta_0\neq 0$. It follows that
\beq*
\frac{\zeta_n}{\eta_n} =\frac{u_2(T^n(x))+(\alpha_n/\beta_n)s_2(T^n(x))}{u_1(T^n(x))+(\alpha_n/\beta_n)s_1(T^n(x))}.
\eeq*
where $\alpha_n/\beta_n\leqslant\lambda^{-2m}\alpha_0/\beta_0\to 0$ as $n\to\infty$. Consequently
\beq*
|\zeta_n/\eta_n-u_2(T^n(x))/u_1(T^n(x))|
\eeq*
tends to zero and, because $(\eta_n,\zeta_n)$ is in the interior of $K(T^n(x))$, for all sufficiently large $n$, so is $E^u(T^n(x))$. The fact that $x$ was arbitrary completes the proof.
\end{proof}

\subsection{The Bernoulli property}\label{PLTM.Bernoulli.Final}

It remains to show that the divergence and orientation of manifolds we have established give sufficient conditions for the Bernoulli property. This amounts to showing that the strong form of the condition~(\ref{eqn:mip_pltm}) is satisfied.

Just as we have established that for `almost every' $z\in R$ the unstable manifold $\gamma^u(z)$ grows exponentially when iterated with $H_S$ and its orientation remains within the cone $C$, so it can be shown that for `almost every' $z'\in R$ the stable manifold $\gamma^s(z)$ grows exponentially when iterated with $H_S^{-1}$ and its orientation remains within the cone $\tilde{C}$.

We need to develop a little terminology. In particular we will introduce a covering space for the manifold $R$, onto which we will lift our local invariant manifolds and deduce intersections. Any intersection of the lifted local manifolds must imply an intersection of the local manifolds themselves. We construct the covering space in two stages.

\begin{figure}[htp]
\centering
\includegraphics[totalheight=0.3\textheight]{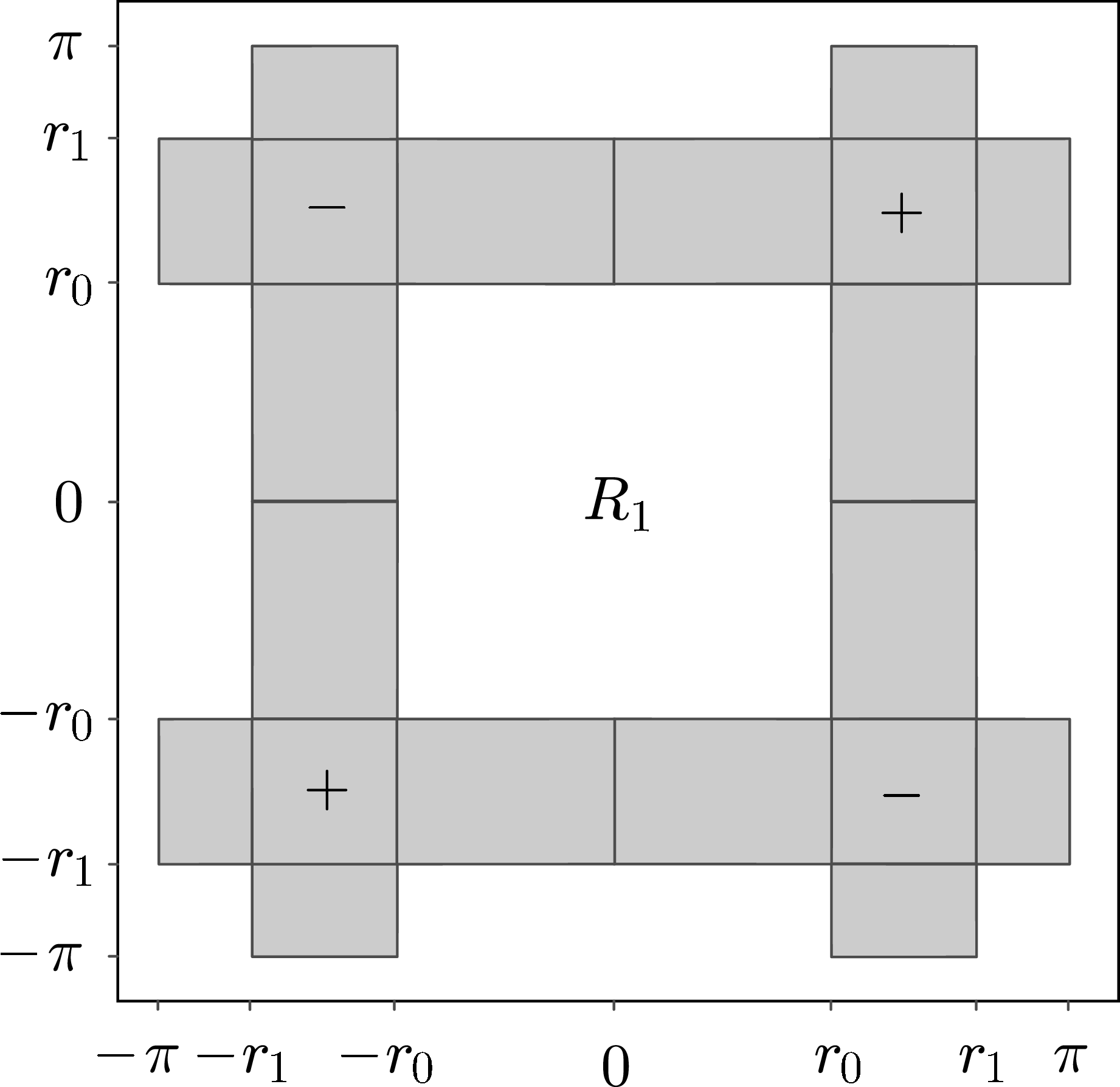}
\caption[The manifold $R_1$]{The manifold $R_1=R\cup R'\subset\mathbb{T}^2$. Together with the map $p':R_1\to R$ this gives a covering space for $R$.}
\label{fig:R_1}
\end{figure}

Let $-R\subset\mathbb{T}^2$ denote those points $(x,y)\subset\mathbb{T}^2$ such that $(-x,-y)\in R$. We notice that $R\cap -R=\emptyset$ and define a manifold $R_1=R\cup -R$. Figure~\ref{fig:R_1} illustrates $R_1$. Let $p':R_1\to R$ be given by $(x,y)\mapsto(x,y)$ if $(x,y)\in R$ and $(x,y)\mapsto(-x,-y)$ otherwise. Then $(R_1,p')$ is a covering space (a double cover, in fact) of $R$. The derivative of $p'$ and its possible inverses preserve the cones $C$ and $\tilde{C}$.

\begin{figure}[htp]
\centering
\includegraphics[totalheight=0.6\textheight]{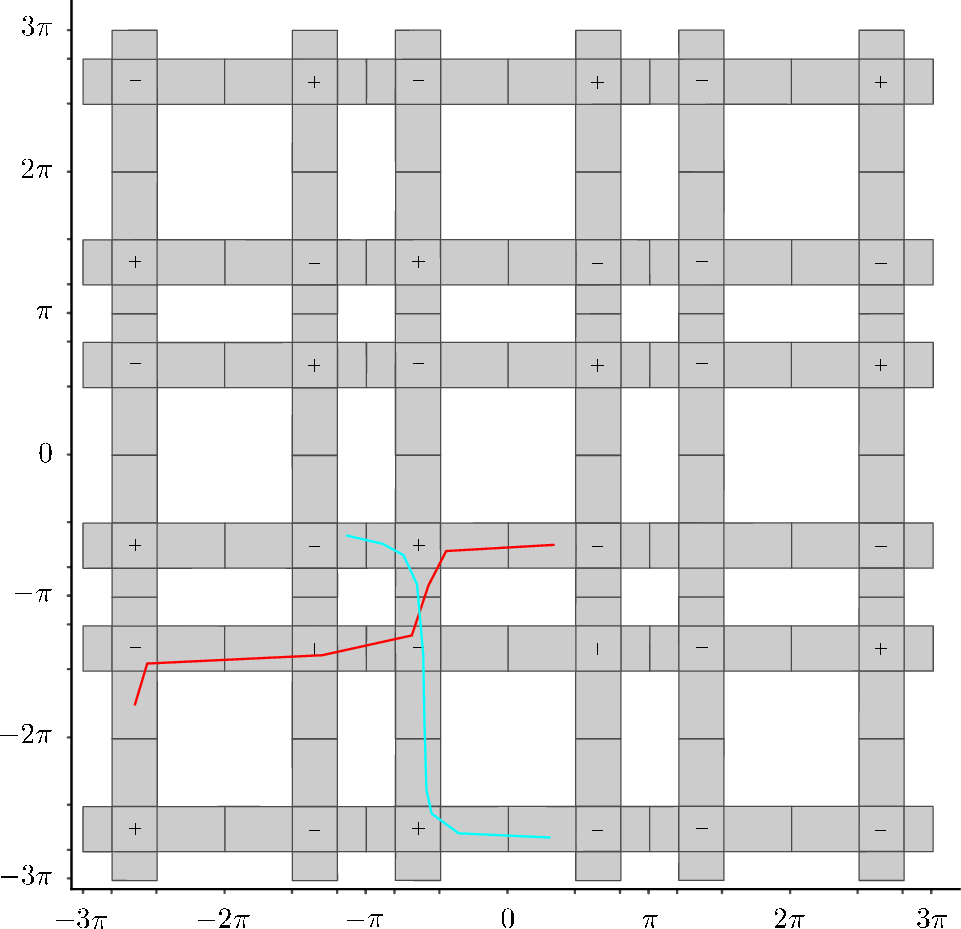}
\caption[The manifold $R_2$]{A portion of the manifold $R_2\subset\mathbb{R}^2$. Together with the map $p:R_2\to R$ this gives a covering space for $R$. In red is a typical piece of the image of a local unstable manifold for some $z\in R$. The gradient at all times is in $C$. Analogously we show a typical piece of the pre-image of a local stable manifold for $z'\in R$. This has gradient in $\tilde{C}$. If they are sufficiently long then they must intersect.}
\label{fig:R_2}
\end{figure}

Recall that $(u,v)$ give Cartesian coordinates in the plane $\mathbb{R}^2$. We define
\beq*
R_2=\left\{(u,v):r_0\leqslant |u-2n\pi|, |v-2m\pi|\leqslant r_1,\text{ for some }m,n\in\mathbb{Z}\right\}.
\eeq*
A portion of $R_2$ is illustrated in Figure~\ref{fig:R_2}. Let $p'':R_2\to R_1$ be the projection which takes each coordinate modulo $\mathbb{S}^1$. Then $(R_2,p'')$ gives a covering space for $R_1$. The derivative of $p''$ and its (locally defined) inverses preserve $C$ and $\tilde{C}$. Thus $(R_2,p=p'\circ p'')$ is a covering space for $R$, and $Dp$, $Dp^{-1}$ preserve $C$ and $\tilde{C}$.

\begin{lem}[Lifting lemma]
Let $z\in R$ and let $\rho$ be a curve into $R$, that is, a continuous map from some interval containing $0$ into $R$. Furthermore let $\rho(0)=z$. Then for any $z''\in R_2$ lying over $z$ (i.e.\ $p(z'')=z$) there is a unique curve $\rho''$ into $R_2$ lying over $\rho$ (i.e.\ $p(\rho'')=\rho$).
\end{lem}
Most topology texts contain a proof of this result; we recommend the book of \citet{armstrong}. We remark that what we have called curves are often called \emph{paths} in the topology literature.

We are ready to prove the main result.

\begin{proof}[Proof of Theorem~\ref{thm:main_plane}]
Recall Theorem~\ref{thm:ks} due to \citet{ks}. It suffices to show that the strong form of the manifold intersection property~(\ref{eqn:mip_pltm}) is satisfied. Expressed in the new coordinates, the condition states that for $\mu$-a.e.\ pair $w,w'\in A$ (writing $z=\Psi\circ M_+^{-1}(w)$ and $z'=\Psi\circ M_+^{-1}(w')$) and for all sufficiently large integers $m$ and $n$ then
\beql
\label{eqn:mip_pltm_new_coords}
H^n(\gamma^u(z))\cap H^{-m}(\gamma^s(z))\neq\emptyset.
\eeql

It has been established that $H^n(\gamma^u(z))$ has orientation in $C$ and that so do any of its lifts to $R_2$. By choosing $n$ large enough a lift may be as long as we would like. Moreover the lifting lemma tells us that any of these lifts is a line with direction `bottom-left' to `top-right', as we have illustrated the manifold in Figure~\ref{fig:R_2}.

Whenever $\gamma^u(z)$ crosses an intersection region (marked $+$ or $-$) horizontally (respectively, vertically) then the next iteration with $F$ (respec.\ $G$) stretches its image vertically (horizontally) adding at least $2\pi$ to its height (width). This shows that the pathological case of $H^n(\gamma^u(z))$ being oriented in $C$ but essentially horizontal or vertical, cannot occur.

We are left with $H^n(\gamma^u(z))$ having arbitrary length and stretching diagonally as previously described. Such a piece is illustrated in red in Figure~\ref{fig:R_2}. Analogous arguments show that $H^{-m}(\gamma^u(z'))$ must be as illustrated in blue and an intersection is inevitable.
\end{proof}

%% file: chapter4/chapter4.tex
\chapter{The Bernoulli property for a linked-twist map on the two-sphere}
\setcounter{equation}{0}
\setcounter{figure}{0}
\label{chapter4}

In this chapter we prove Theorem~\ref{thm:main_sphere} which says that the linked-twist map on $R\subset\mathbb{S}^2$, defined in Section~\ref{Intro.Theorems.Sphere}, has the Bernoulli property.

In Section~\ref{SLTM.DTLTM} we introduce a \emph{generalised} linked-twist map $H:R\to R$ on a manifold $R\subset\mathbb{T}^2$. We append the word `generalised' because the new map does not fit the definition of an abstract linked-twist map given in Section~\ref{Intro.Background.Abstract}, consisting as it does of twist maps defined on four annuli rather than two. Each of the twists is linear.

We are able to prove, by application of the techniques introduced by \citet{woj} and amended by \citet{sturman} in light of the work of \citet{ks}, that the new map has the Bernoulli property. In Section~\ref{SLTM.Toral} we deal with some technical issues concerning the nature of points at which $H$ is non-differentiable and conclude that Lyapunov exponents exist Lebesgue-almost everywhere. We then describe in some detail the return of points to a certain region $S\subset R$ before using this description to show that these Lyapunov exponents are non-zero.

In Section~\ref{SLTM.Toral2} we deal with the `global' aspects of the argument. We consider the orientation of local invariant manifolds and give a rigorous original proof of their growth. Such a result is lacking in the literature and is similar to the analogous result from Chapter~\ref{chapter3}; with it we are able to prove the Bernoulli property for $H$.

In Section~\ref{SLTM.Main} we show that $H$ is semi-conjugate to the linked-twist map $\Theta$ on $\mathbb{S}^2$ and from this conclude that the latter is also Bernoulli. Our proof relies on a result of \citet{orn2}.

\input chapter4/section1

\input chapter4/section2

\input chapter4/section3
\input chapter4/section4

%% file: chapter4/section1.tex
\section{A generalised linked-twist map on the two-torus}\label{SLTM.DTLTM}
\setcounter{equation}{0}

In this section we define a more general linked-twist map on the two-torus $\mathbb{T}^2$. It differs from the map introduced in Section~\ref{Intro.Theorems.Torus} in that it is constructed by embedding not two but four cylinders into $\mathbb{T}^2$.

We will define the generalised linked-twist map directly on the two-torus; we take it to be elementary that the manifold we construct may be obtained by embedding four cylinders into $\mathbb{T}^2$ and do not give a formal proof.

\subsection{Definition of the map}
\begin{figure}[htp]
\centering
\includegraphics[totalheight=0.3\textheight]{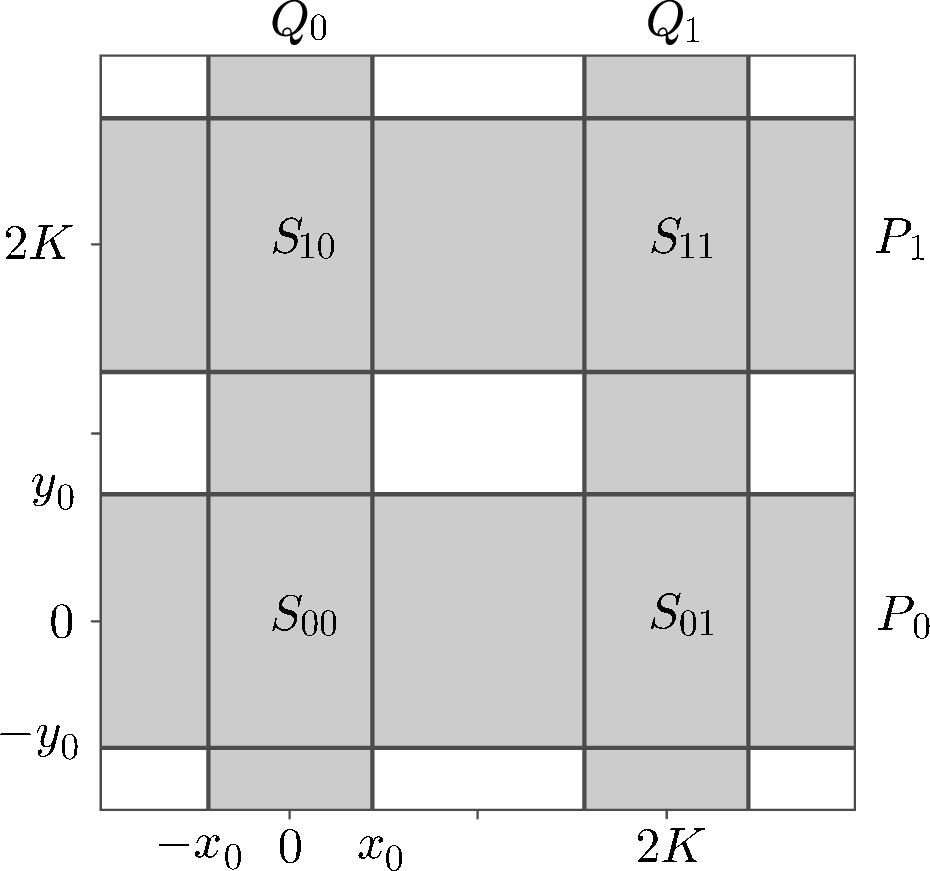}
\caption[The manifold $R\subset\mathbb{T}^2$]{The manifold $R\subset\mathbb{T}^2$ is shaded.}
\label{fig:bigtwist}
\end{figure}

Recall our construction in Section~\ref{Intro.Theorems.Sphere} of the linked-twist map $\Theta:A\to A$, where we described the circle $\mathbb{S}^1$ as having a coordinate that is periodic with period $4K$ and where
\beql\label{eqn:JacobiK}
K=\int_0^{\pi/2}\left(1-\frac{1}{2}\sin^2t\right)^{-1/2}\d t\approx 1.85.
\eeql
Presently it will be most convenient to let this coordinate take values between $-K$ and $3K$, where of course these particular values are identified. Let $0<x_0,y_0<K$ be as previously defined. We define two `horizontal' and two `vertical' annuli on the two-torus:
\beq*
\begin{array}{cc}
P_0=\{(x,y):x\in\mathbb{S}^1,|y|\leqs y_0\},\quad &
P_1=\{(x,y):x\in\mathbb{S}^1,|y-2K|\leqs y_0\}, \\
Q_0=\{(x,y):|x|\leqs x_0,y\in\mathbb{S}^1\},\quad &
Q_1=\{(x,y):|x-2K|\leqs x_0,y\in\mathbb{S}^1\}.
\end{array}
\eeq*
Denote by $P=P_0\cup P_1$ the union of the horizontal annuli and by $Q=Q_0\cup Q_1$ the union of the vertical annuli. We let $R=P\cup Q$ and $S_{hi}=P_h\cap Q_i$ for each pair $h,i\in\{0,1\}$. Let $S$ denote the union of the four regions $S_{hi}$. Figure~\ref{fig:bigtwist} illustrates the manifold $R\subset\mathbb{T}^2$.

We define a twist map on each of the four annuli and group these into two `families'; those maps defined on horizontal annuli are in the first family and those defined on vertical annuli are in the second. Both maps in a given family act simultaneously and the linked-twist map will be the composition of the two families.

\begin{figure}[htp]
\centering
\includegraphics[totalheight=0.2\textheight]{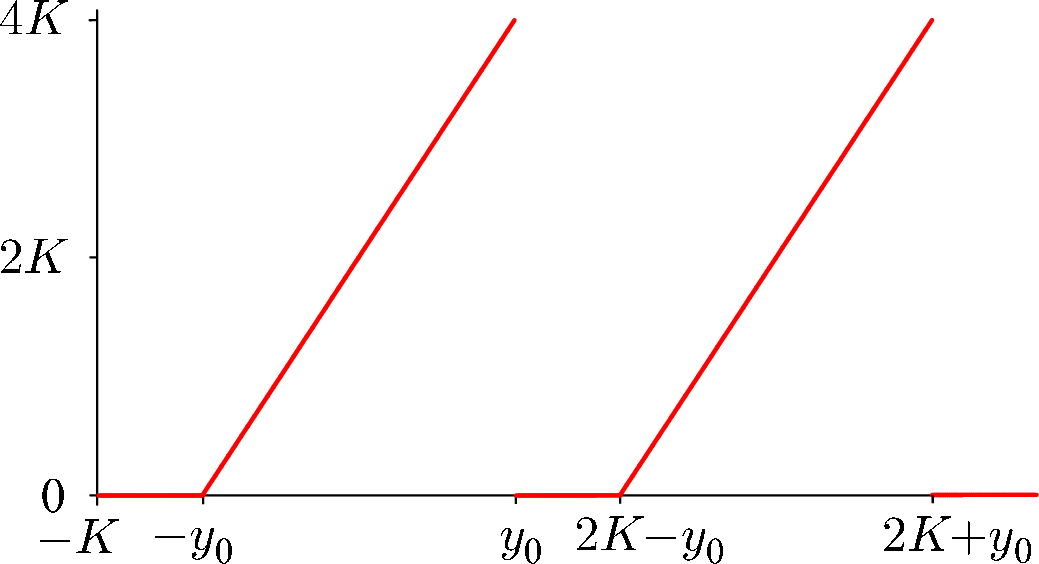}
\caption[The generalised twist function $f$]{The generalised twist function $f:\mathbb{S}^1\to\mathbb{S}^1$ is continuous everywhere and differentiable except at the four points $\pm y_0$ and $2K\pm y_0$. The function $g:\mathbb{S}^1\to\mathbb{S}^1$ is similar.}
\label{fig:f_of_y}
\end{figure}

Let $f:\mathbb{S}^1\to\mathbb{S}^1$ (see Figure~\ref{fig:f_of_y}) and $g:\mathbb{S}^1\to\mathbb{S}^1$ be given by
\beq*
f(y)=\left\{
\begin{array}{r@{\quad}l}
4K(y+y_0)/(2y_0) &\text{if }y\in [-y_0,y_0], \\
4K(y+y_0-2K)/(2y_0) &\text{if }y\in [2K-y_0,2K+y_0], \\
0 & \text{otherwise,} \end{array} \right.,
\eeq*
\beq*
g(x)=\left\{
\begin{array}{r@{\quad}l}
4K(x+x_0)/(2x_0) &\text{if }x\in [-x_0,x_0], \\
4K(x+x_0-2K)/(2x_0) &\text{if }x\in [2K-x_0,2K+x_0], \\
0 & \text{otherwise.} \end{array} \right.
\eeq*

A \emph{family of horizontal twist maps} $F:\mathbb{T}^2\to\mathbb{T}^2$ is given by 
\beq*
F(x,y)=(x+f(y),y)
\eeq*
and a \emph{family of vertical twist maps} $G:\mathbb{T}^2\to\mathbb{T}^2$ is given by
\beq*
G(x,y)=(x,y+g(x)).
\eeq*

We say that the map $H:R\to R$ given by the composition $H=G\circ F$ is a \emph{generalised linked-twist map} on $R$. By analogy with the toral linked-twist map introduced in Section~\ref{Intro.Theorems.Torus} it is easy to see that $H$ preserves the Lebesgue measure on $R$, which we will denote by $\mu$. In Sections~\ref{SLTM.Toral} and~\ref{SLTM.Toral2} we prove the following:
\begin{thm}
\label{thm:Bernoulli_DTLTM}
The map $H:R\to R$ is metrically isomorphic to a Bernoulli shift.
\end{thm}

%% file: chapter4/section2.tex
\section{Non-zero Lyapunov exponents for $H$}\label{SLTM.Toral}
\setcounter{equation}{0}

A major step toward the proof that $H$ has the Bernoulli property is to show that its Lyapunov exponents are non-zero. The present section is divided into three parts. In Section~\ref{SLTM.Toral.Tech} we give some results of a technical nature which satisfy certain hypotheses in the theorem of \citet{ks} (Theorem~\ref{thm:ks}). This allows us to conclude that Lyapunov exponents exist $\mu$-almost everywhere. In Section~\ref{SLTM.Toral.Return} we describe the return of a `typical' trajectory to the region $S$; this behaviour will be crucial in determining hyperbolicity. Finally in Section~\ref{SLTM.Toral.NonZero} we describe the behaviour of the differential $DH$ acting on certain tangent cones, from which the result follows.

\subsection{Technical details}\label{SLTM.Toral.Tech}

Our first consideration is to determine the nature of those points in $R$ at which $H$ and its iterates are non-differentiable. For this we will need to establish some notation for boundaries of the annuli. The set
\beq*
\partial P_0=\{(x,y):x\in\mathbb{S}^1,y\in\{y_0,y_1\}\}
\eeq*
denotes the boundary of annulus $P_0$, with the boundaries of $P_1$, $Q_0$ and $Q_1$ being defined similarly.
\footnote{This notation is necessarily different from that in Chapter~\ref{chapter2} where both our map and our ambitions differed.}
The boundaries of $P$ and $Q$ are denoted by $\partial P=\partial P_0\cup\partial P_1$ and $\partial Q=\partial Q_0\cup\partial Q_1$ respectively.

The twist function $f$ is differentiable provided that $y,y-2K\neq\pm y_0$ and similarly $g$ is differentiable provided that $x,x-2K\neq\pm x_0$. It follows that $F$ is differentiable on $\mathbb{T}^2\backslash\partial P$ and that $G$ is differentiable on $\mathbb{T}^2\backslash\partial Q$. Consequently $H$ is differentiable on $\mathbb{T}^2\backslash s_0$, where
\beq*
s_0=\pd P\cup F^{-1}(\pd Q).
\eeq*

It may be useful for the reader to briefly review Section~\ref{LitReview.Hyp.smws} at the present moment, where we described the theorem of \citet{ks} and the notation necessitated by its statement. We remark that the function $f$ from that section corresponds to the function $H$ from this section; the function $f$ in this section has no analogue in Section~\ref{LitReview.Hyp.smws}.

We define a number of full-measure subsets of $R$ and hope that Table~5.1 is of use to the reader in keeping track of these.

\begin{table}
\centering
\begin{tabular}{|c|l|}\hline
Set& Description\\ \hline
$R$& Union of the four annuli\\
$V$& Open subset; a Riemannian manifold\\
$N$& Open subset on which $H$ is defined\\
$J$& Intersection of all images and pre-images of $N$\\ \hline
\end{tabular}
\caption{Full measure subsets of the complete metric space R. As listed, each set contains those below it.}
\end{table}

Let $\rho$ be a metric on $R$, so that $(R,\rho)$ is a complete metric space. Notice that $V=R\backslash s_0$ is a smooth manifold. $V$ has full $\mu$-measure and is open. Let $N=V\backslash H^{-1}(s_0)$; $N$ also has full $\mu$-measure and is open.

\begin{prop}
\label{prop:smws}
$H|_N:N\to V$ is a smooth map with singularities.
\end{prop}
\begin{proof}
$N$ is open, so given any $z\in N$ the exponential map $\exp_z$ is an injective map from some neighbourhood of $0\in T_z\mathbb{T}^2$ to some neighbourhood $U(z)\subset N$ of $z$. The only restricting factors on the size that this neighbourhood may take are the finite size of $N$ itself and the distance from a given point to the set $\text{sing}(H)=V\backslash N$. Furthermore if $z\in N$ then $z\notin\partial P$ and $F(z)\notin\partial Q$ so $H$ is smooth and injective.
\end{proof}

In our next lemma it will be convenient to consider the Lebesgue measure on $\mathbb{T}^2\supset R$ and we denote this by $\mu_{\mathbb{T}^2}$. For $\mu_{\mathbb{T}^2}$-measurable sets $K\subset\mathbb{T}^2$, $\mu$ (the Lebesgue measure on $R$) is the conditional measure such that
\beq*
\mu(K)=\mu_{\mathbb{T}^2}(K|R)=\frac{\mu_{\mathbb{T}^2}(K\cap R)}{\mu_{\mathbb{T}^2}(R)}.
\eeq*
It follows that for measurable sets $K'\subset R\subset\mathbb{T}^2$ we have
\beql
\label{eqn:relate_measures}
\mu_{\mathbb{T}^2}(K')=\mu(K')\mu_{\mathbb{T}^2}(R).
\eeql

We use the notation $B_{\mathbb{T}^2}(z,\varepsilon)$ to denote the open $\varepsilon$-neighbourhood of $z\in\mathbb{T}^2$ and we use $B(z,\varepsilon)$ to denote the intersection $B_{\mathbb{T}^2}(z,\varepsilon)\cap R$. For a non-empty set $K\subset\mathbb{T}^2$ we let 
\beq*
B_{\mathbb{T}^2}(K,\varepsilon)=\bigcup_{z\in K}B_{\mathbb{T}^2}(z,\varepsilon)
\eeq*
and similarly define $B(K,\varepsilon)$.

\begin{lem}
\label{lem:t_ks1}
$H|_N$ satisfies the condition (\ref{eqn:KS1}).
\end{lem}
We also refer to this as the condition (KS1).
\begin{proof}
We have $\text{sing}(H)=H^{-1}(s_0)=H^{-1}\left(\pd P\cup F^{-1}(\pd Q)\right)$. In $\mathbb{T}^2$, an $\varepsilon$-neighbour-hood of $\partial P$ consists of four strips (rectangles) each of width~$2\varepsilon$ and length~1, thus $\mu_{\mathbb{T}^2}(B_{\mathbb{T}^2}(\partial P,\varepsilon))=8\varepsilon$. $B(\partial P,\varepsilon)$ is a subset of $B_{\mathbb{T}^2}(\partial P,\varepsilon)$ and so its $\mu_{\mathbb{T}^2}$-measure is at most $8\varepsilon$. We are interested in its $\mu$-measure and, by~(\ref{eqn:relate_measures}) above, these are related by
\beq*
\mu_{\mathbb{T}^2}(B(\partial P,\varepsilon))=\mu(B(\partial P,\varepsilon))\mu_{\mathbb{T}^2}(R),
\eeq*
all of which gives $\mu(B(\partial P,\varepsilon))\leqslant 8\varepsilon/\mu_{\mathbb{T}^2}(R)$.

$F^{-1}(\partial Q)$ consists of four piecewise-straight lines, each of length bounded above by some integer $N$. Thus a similar argument applies here. Taking the pre-image with respect to $H$ will `stretch' these strips somewhat, but similar bounds will still apply so that (KS1) will be satisfied with some constant $c\in\mathbb{R}$ and with $a=1$.
\end{proof}

Denote by $DF$, $DG$ and $DH$ the derivatives of $F$, $G$ and $H$ respectively. We use subscripts to denote the point in $R$ at which the derivative is evaluated, so for example at each $z\in N$ we have the familiar chain rule
\beq*
DH_z=DG_{F(z)}DF_z.
\eeq*
Let $J=\bigcap_{i=-\infty}^{\infty}H^i(N)$; we prove that this has full $\mu$-measure:
\begin{lem}
$\mu(J)=1$.
\end{lem}
\begin{proof}
We use the $H$-invariance of $\mu$ and the countable additivity property. From elementary set theory we have $J=R\backslash\{R\backslash J\}$ and so
\beqa*
\mu(J) &= \mu(R)-\mu(R\backslash J) \\
			 &= \mu(R)-\mu\left(\bigcup_{i=-\infty}^{\infty}R\backslash H^i(N)\right) \\
			 &\geqs \mu(R)-\sum_{i=-\infty}^{\infty}\mu\left(R\backslash H^i(N)\right) \\
			 &=\mu(R).
\eeqa*
The fact that $J\subset R$ completes the proof.
\end{proof}

\begin{lem}
\label{lem:t_os}
Lyapunov exponents for $H|_N$ exist for $\mu$-a.e.\ $z\in R$ and every non-zero $v\in T_z\mathbb{T}^2$.
\end{lem}
\begin{proof}Let $z=(x,y)\in J$. The derivatives of $F$ and $G$ are given by
\beq*
DF=\left( \begin{array}{cc}
1 & f'(y) \\ 0 & 1
\end{array}\right) \quad\text{and}\quad
DG=\left( \begin{array}{cc}
1 & 0 \\ g'(x) & 1
\end{array}\right)
\eeq*
respectively. Each takes only finitely many values on $R$ meaning that there are upper and lower bounds on its eigenvalues. Thus
\beq*
\int_R\log^+\|DH^{\pm 1}\|d\mu=\int_J\log^+\|DH^{\pm 1}\|d\mu<\infty
\eeq*
(the norms are operator norms), showing that $H$ satisfies (\ref{eqn:oseledec}). The result follows immediately from the multiplicative ergodic theorem, Theorem~\ref{thm:met}.
\end{proof}

\begin{lem}
\label{lem:t_ks2}
$H|_N$ satisfies the condition (\ref{eqn:KS2}).
\end{lem}
We also refer to this as the condition (KS2).
\begin{proof}
Because our twist functions $f$ and $g$ are piecewise-linear, their derivatives (where they exist) are constant functions. Consequently all second derivatives will be zero and condition (KS2) is satisfied with $b=1$ and any $c_2>0$.
\end{proof}

Lemmas~\ref{lem:t_ks1}, \ref{lem:t_os} and~\ref{lem:t_ks2} and part~(a) of Theorem~\ref{thm:ks} yield two conclusions. First, Lyapunov exponents exist for $\mu$-a.e.\ $z\in R$ and for every non-zero tangent vector at $z$. Second, for those $z$ at which we have a positive exponent, the local unstable manifold $\gamma^u(z)$ exists and is of the form (\ref{eqn:st}). An analogous statement holds for negative exponents and the local stable manifold $\gamma^s(z)$.

Throughout this section we have referred to the restriction $H|_N$ in keeping with the notation of \citet{ks}. We have seen that $N$ has full $\mu$-measure and so statements that are true for a.e.\ $z\in N$ are equally true for a.e.\ $z\in R$. From here on we will speak only of $H$ rather than its restriction to $N$, this leading to cleaner expressions whilst losing none of the accuracy.

\subsection{Return of trajectories to the region $S$}\label{SLTM.Toral.Return}

We will need to describe in some detail the trajectory of a `typical' point $z\in R$ (we will say precisely which points are `typical' soon; for now let us just say that the set of such points will have full measure). In particular we need to analyse first whether, and if so how often, such a point returns to the union of intersection regions $S$. Our main tool for doing this will be the \emph{return map} to the set $S$, which is defined as follows:

\begin{defn}[First (and $n^{\text{th}}$) return map to $S$ for $z$]
For $z\in R$, let $n$ be the smallest (strictly positive) integer such that $H^n(z)\in S$, if such an $n$ exists. If it does then the map \beq*
H_S(z)=H^n(z)
\eeq*
will be called the first return map to $S$ for $z$, often abbreviated to just first return map where it is clear for which $z$ the map is defined. Similarly
\beq*
H^2_S(z)=H_S(H^n(z))\circ H_S(z)
\eeq*
will be called the second return map and the general case follows inductively.
\end{defn}

In the remainder of this section we describe the notation and formalism introduced by \citet{woj} in his study of toral linked-twist maps defined on two annuli. It will enable us to give the description we require. Let $p,q:R\to[0,\infty]$ be defined as follows:
\beq*
p(z)=\left\{
\begin{array}{r@{\quad}l}
0 &\text{if }z\in Q\backslash P, \\
p\in\mathbb{N} &\text{if }z\in P,\quad F^n(z)\notin Q\text{ for }0<n<p\text{ and }F^p(z)\in Q, \\
\infty &\text{if }z\in P\text{ and }F^n(z)\notin Q\text{ for all }n\in\mathbb{N},
\end{array} \right.
\eeq*
\beq*
q(z)=\left\{
\begin{array}{r@{\quad}l}
0 &\text{if }z\in P\backslash Q, \\
q\in\mathbb{N} &\text{if }z\in Q, G^n(z)\notin P\text{ for }0<n<q\text{ and }G^q(z)\in P, \\
\infty &\text{if }z\in Q\text{ and }G^n(z)\notin P\text{ for all }n\in\mathbb{N}.
\end{array} \right.
\eeq*
We give a brief discussion. Let $z\in P$, so $p(z)\neq 0$ and assume also that $p(z)\neq\infty$. Clearly $F^n(z)\in P$ for all natural numbers $n$, and $p(z)$ is the first such natural number for which $F^p(z)\in Q$, i.e.\ the first time the forward trajectory under $F$ lands in $S$. This is of interest because for $0<n<p(z)$ we have $F^n(z)\notin Q$ and so $G(F^n(z))=F^n(z)$.

To phrase this another way, starting with $z$ as above, if we follow its orbit under $H=G\circ F$ then the $p(z)^{\text{th}}$ iterate is the first for which the $G$ component of $H$ is \emph{not} the identity. A similar property holds for the function $q$.

\citet{woj} outlined an argument that the trajectory of $\mu$-a.e.\ $z\in R$ returns to $S$ infinitely many times. The proof is quite elementary and we have discussed an analogous result in Section~\ref{PLTM.Wojtkowski}. We denote by $R_{\tilde{S}}\subset R$ the full measure (and clearly invariant) set of points for which this occurs. 

Define two sequences of functions $p_1:R_{\tilde{S}}\to[0,\infty)$, $p_i:R_{\tilde{S}}\to[1,\infty)$ for $i=2,3,4,...$ and $q_i:R_{\tilde{S}}\to[1,\infty)$ for $i=1,2,3,...$ as follows:\footnote{We will often neglect to express the dependence of $p_i,q_i$ on the point $z$ for the benefit of cleaner notation, where this will not lead to confusion; for the same reason we typically neglect the composition symbol `$\circ$'.}
\beq*
\begin{array}{ll@{\quad}ll}
p_1(z) &= p(z), & q_1(z) &= q(F^{p_1}(z)),\\
p_{i+1}(z) &= p(G^{q_i}F^{p_i}\cdots G^{q_1}F^{p_1}(z)),
& q_{i+1}(z) &= q(F^{q_{i+1}}G^{q_i}F^{p_i}\cdots G^{q_1}F^{p_1}(z)).
\end{array}
\eeq*
Given $z\in R_{\tilde{S}}$ the above scheme defines inductively $p_i$ and $q_i$ for all natural numbers $i$. It is clear that for $z\in R_{\tilde{S}}$ one has $G^{q_1}(F^{p_1}(z))\in S$; for $z\in R_{\tilde{S}}\cap S$ then $F^{p_1}(z)\in S$ also.

Finally define functions $m_1:R_{\tilde{S}}\to[0,\infty)$ and $m_i:R_{\tilde{S}}\to[1,\infty)$ for $i=2,3,4,...$ by
\beq*
m_i(z)=\sum_{n=1}^i\left(p_n(z)+q_n(z)-1\right).
\eeq*
The description of the first-return map that we required is given by the following easy lemma:

\begin{lem}
\label{lem:1st_return_map}
For $z\in R_{\tilde{S}}$, the first return map to $S$ for $z$ is given by
$H_S=H^{m_1(z)}$. In general, for $n\in\mathbb{N}$, the $n^{\text{th}}$ return map to $S$ is given by $H_S^n=H^{m_n(z)}$.
\end{lem}
\begin{proof}
We have
\beq*
\begin{array}{ll}
H^{m_1}&=
\underbrace{\left(G\circ F\right)\circ\cdots\circ\left(G\circ
F\right)}_{q_1-1} \circ
\underbrace{\left(G\circ F\right)\circ\cdots\circ\left(G\circ
F\right)}_{p_1} \\
&=\underbrace{\left(G\circ id\right)\circ\cdots\circ\left(G\circ
id\right)}_{q_1-1} \circ
\underbrace{\left(G\circ F\right)\circ
\left(id\circ F\right)\circ\cdots\circ\left(id\circ
F\right)}_{p_1}
\end{array}
\eeq*
i.e.\ $H^{m_1}=G^{q_1}\circ F^{p_1}$. The first result follows easily from the definitions of $p_i$ and $q_i$ and the general case by induction on $n$.
\end{proof}

\subsection{Lyapunov exponents are non-zero}\label{SLTM.Toral.NonZero}

In this section we use the first return map to show that Lyapunov exponents are non-zero almost everywhere. Our first step is to describe the derivative of the first return map.

Let $z\in R_{\tilde{S}}\cap J$ and denote $z'=F^{p_1(z)}(z)$, then the derivative of the first-return map is given by $DH_S=DH_z^{m_1}=DG_{z'}^{q_1}DF_z^{p_1}$, where $p_1,q_1$ and $m_1$ are to be evaluated at $z$. Similarly for $n\in\mathbb{N}$ the derivative of the $n^{\text{th}}$ return map is given by
\beq*
DH^{m_n}_z=DG^{q_n}_{z_{n-1}'}DF^{p_n}_{z_{n-1}^{}}\cdots
DG^{q_2}_{z_1'}DF^{p_1}_{z_1^{}}DG^{q_1}_{z_0'}DF^{p_1}_{z_0^{}},
\eeq*
where $z_0^{}=z,z_0'=z'$ and in general $z_i^{}=H^{m_i}(z),z_i'=F^{p_{i+1}}(z_i)$ for $i\in\mathbb{N}$.

Let $(u,v)=(\d x,\d y)$ give coordinates in the tangent space $T_z\mathbb{T}^2$, which we identify with $\mathbb{R}^2$. We define the cone
\beq*
C=\{(u,v)\in\mathbb{R}^2:uv>0\}\subset T_z\mathbb{T}^2,
\eeq*
that is the open first and third quadrants of the plane, and introduce the standard Euclidean norm $\|\cdot\|:\mathbb{R}^2\to[0,\infty)$ on $T_z\mathbb{T}^2$.

We show that for $\mu$-a.e.\ $z\in R$ and for each $w\in C\subset T_z\mathbb{T}^2$ the Lyapunov exponent
\beq*
\chi^+(z,w)=\lim_{n\to\infty}\frac{1}{n}\log\left\|DH_z^nw\right\|
\eeq*
is positive and thus there is a local unstable manifold $\gamma^u(z)$ of the form (\ref{eqn:st}). We do not do so explicitly, but the results are easily re-formulated to show that if
\beq*
\tilde{C}=\{(u,v)\in\mathbb{R}^2:uv<0\}\subset T_z\mathbb{T}^2
\eeq*
then for $\mu$-a.e.\ $z\in R$ and for each $\tilde{w}\in\tilde{C}$ the Lyapunov exponent
\beq*
\chi^-(z,\tilde{w})=\lim_{n\to-\infty}\frac{1}{n}\log\left\|DH_z^nw\right\|
\eeq*
is less than zero, and so a local stable manifold $\gamma^s(z)$ exists also.

The first step is to establish some results concerning the growth of a vector $w\in C$ when acted upon by the derivative of $H$ or of $H_S$. Similar results were proven in \citet{woj}. The result illustrates a key feature of the dynamics: the return map is uniformly hyperbolic.

We use the convention that $\mathbb{N}$ denotes the positive integers; in particular $0\notin\mathbb{N}$.

\begin{lem}
\label{lem:t_non-shrink}
Let $z\in R_{\tilde{S}}\cap J$ and let $w\in C\subset T_z\mathbb{T}^2$. Let $n,r\in\mathbb{N}$ with $r\geqslant n$. Then
\begin{enumerate}
\item[\emph{(a)}] $DH^n_zw\in C\subset T_{H^n(z)}\mathbb{T}^2$,
\item[\emph{(b)}] $\left\|DH^r_zw\right\|\geqslant\left\|DH^n_zw\right\|$.
\end{enumerate}
If additionally $z\in S$ then for a constant $\lambda>1$, which is independent of $z$, we have
\begin{enumerate}
\item[\emph{(c)}] $\left\|DH_z^{m_n(z)}w\right\|\geqslant\lambda^n\|w\|$.
\end{enumerate}
\end{lem}
We remark that part~(a) of the lemma says that the cone $C$ is invariant under the derivative of $H^n$ and that part~(b) is a `non-shrinking' condition, saying that the norm of vectors in $C$ does not decrease when mapped by $DH$. Part~(c) expresses the uniform growth resulting from the derivative of the return map. The extra condition that $z\in S$ will not limit our ability to make use of this lemma, because we are interested in cumulative growth along a trajectory and we have already established that $\mu$-a.e.\ trajectory returns infinitely many times to $S$.
\begin{proof}
Let $z=(x,y)\in R_{\tilde{S}}\cap J$ and let $z'=(x',y')=F(z)$. We have
\beq*
DH_z=\left( \begin{array}{cc}
1 & \alpha \\ \beta & 1+\alpha\beta
\end{array} \right),
\eeq*
where we have introduced the notation $\alpha=f'(y)$ and $\beta=g'(x')$. At most one of $\alpha, \beta$ may be zero and the other(s) positive, because $\alpha=0$ requires $(x,y)\notin P$ and $\beta=0$ requires $F(x,y)\notin Q$, which cannot occur in succession.

Now let $w=(u,v)\in C\subset T_z\mathbb{T}^2$; straight-forward calculation shows that $DH_zw\in C$ and that the Euclidean norm of this vector is at least that of $w$. Induction on $n$ yields the result~(a) and that $\left\|DH^n_zw\right\|\geqslant\|w\|$ for $n\in\mathbb{N}$. The result~(b) follows easily.

Now assume that $z\in S$ and redefine $z'=(x',y')=F^{p_1(z)}(z)$. (The functions $p_i,q_i,m_i$ are always to be evaluated at $z$.)  Notice that (the redefined) $\alpha,\beta$ are now both strictly positive. We have
\beq*
DH^{m_1}_z=DG^{q_1}_{z'}DF^{p_1}_z=
\left( \begin{array}{cc}
1 & p_1\alpha \\ q_1\beta &
1+p_1\alpha q_1\beta \end{array} \right)
\eeq*
where $\alpha=f'(y)>0$ is a positive constant, as is $\beta=g'(x')>0$. Let $\kappa=\min\{p_1\alpha,q_1\beta\}>0$ and $\lambda=\sqrt{1+\kappa^2}>1$. Evaluating $\left\|DH_z^{m_1}w\right\|$ gives
\beq*
\sqrt{u^2(1+q_1^2\beta^2)+v^2(p_1^2\alpha^2+(1+p_1\alpha q_1\beta)^2) +uv(2p_1\alpha+2q_1\beta(1+p_1\alpha q_1\beta))}
\eeq*
for which a lower bound is $\sqrt{(u^2+v^2)(1+\kappa^2)} =\lambda\|w\|$. This proves part~(c) in the case $n=1$. The general case follows by induction on $n$.
\end{proof}

The previous lemma shows that the tangent cone $C$ is uniformly expanded by $DH_S$, but not necessarily by $DH$. In order to demonstrate the non-vanishing of Lyapunov exponents we must consider the \emph{frequency} with which a typical trajectory returns to the intersection region. To that end we use the following result of Burton and Easton:
\begin{lem}[\citet{be}]
\label{lem:be}
Let $X$ be a compact metric space, $\mu$ a Borel measure on $X$ and $T:X\to X$ a $\mu$-preserving homeomorphism. Suppose $Y\subset X$ is measurable. For $\mu$-a.e.\ $x\in X$ such that $T^n(x)\in Y$ for \emph{some} $n\in\mathbb{N}$, that forward orbit will return to $Y$ with positive frequency in $n$.
\end{lem}
The result we want is an easy corollary of Lemma~\ref{lem:be} and the fact that $\mu$-a.e.\ $z\in R$ enters~$S$:
\begin{cor}
\label{cor:t_S}
For $\mu$-a.e.\ $z\in R$, the forward orbit $\mathcal{O}^+(z)=\{H^n(z):n\in\mathbb{N}\}$ lands in $S$ with positive frequency in $n$, i.e.\ there is a set $R_S\subset R_{\tilde{S}}\subset R$, $\mu(R_S)=1$, and for each $z\in R_S$ the limit
\beq*
\delta(z)=\lim_{n\to\infty}\frac{n}{m_n(z)}
\eeq*
exists and is strictly positive.
\end{cor}

We remark that Lyapunov exponents (where they exist), being infinite-time limits, are invariant along a given trajectory. That is, if $\chi(z,w)$ exists for some $z\in R$ and non-zero $w\in T_z\mathbb{T}^2$, and if $n\in\mathbb{N}$, then $\chi(H^n(z),DH^n_zw)$ exists and is equal to $\chi(z,w)$. The main result of this section is the following.

\begin{prop}
\label{lem:t_LE}
For $\mu$-a.e.\ $z\in R$ and for every tangent vector $w\in C\subset T_z\mathbb{T}^2$ at $z$, the Lyapunov exponent $\chi^+(z,w)$ for the map $H$ is positive.
\end{prop}
\begin{proof}
Let $z\in R_S\cap J$ and let $w\in C\subset T_z\mathbb{T}^2$. We denote $z'=H^{m_1}(z)$ and
$w'=DH_z^{m_1}w$. Clearly $z'\in R_S\cap J\cap S$ and by Lemma~\ref{lem:t_non-shrink} we have $w'\in C\subset T_{z'}\mathbb{T}^2$. By the invariance of Lyapunov exponents along a given trajectory it is equivalent to show that $\chi^+(z',w')>0$.

We claim that for any $z'\in R_S\cap J\cap S$ there is a positive constant $N$ and if $n>N$ then $n>m_{\lfloor rn\rfloor}$, where $r(z)>0$ is a constant to be determined. This condition means that for the trajectory of $z'$, after some `transition period' given by $N$, there is a positive lower bound on how frequently it hits $S$.

Suppose for a moment that our claim is verified. Lemma~\ref{lem:t_non-shrink} implies that, corresponding to each such return, the tangent cone $C$ is expanded by a factor $\lambda>1$. The situation is as follows:
\beq*
\begin{array}{ll}
\chi^+(z',w')&=\lim_{n\to\infty}\frac{1}{n}\log\left\|DH^n_{z'}w'\right\| \\ &\geqslant\lim_{n\to\infty}\frac{1}{n}\log\left\|DH_{z'}^{m_{\lfloor
rn\rfloor}}w'\right\| \\ &\geqslant\lim_{n\to\infty}\frac{1}{n}\log\lambda^{\lfloor
rn\rfloor}\|w'\|\\ &=r\log\lambda >0
\end{array}
\eeq*
where the second line follows from our claim (we may, of course, assume $n>N$ in the limit $n\to\infty$) and from Lemma~\ref{lem:t_non-shrink}, the third
from the same lemma and the fourth from elementary properties of the logarithm.

It remains only to prove the claim. Fix $z'\in R_S\cap J\cap S$. By Corollary~\ref{cor:t_S} there exists
$\delta(z')=\lim_{n\to\infty}n/m_n(z')>0$. In particular, there exists $N>0$ such that if $n>N$ then
$n/m_n(z')>\delta(z')/2$, or equivalently that $2n/\delta(z')>m_n(z)$. If we make the substitution $n'=2n/\delta(z')$ then it follows that there exists $N'=2N/\delta(z')>0$ such that for all $n'>N'$ we have $n'>m_{\lfloor n'\delta(z')/2\rfloor}$, which is the statement of our claim, with $r(z')=\delta(z')/2>0$.
\end{proof}

We have remarked that it is completely analogous to show that for $\mu$-a.e.\ $z\in R$ and for each non-zero $\tilde{w}\in\tilde{C}\subset T_z\mathbb{T}^2$ the Lyapunov exponents $\chi^-(z,\tilde{w})$ are strictly negative. $R$ has dimension two, so by the elementary theory of Lyapunov exponents (see \citet{bap}) there are at most two distinct Lyapunov exponents associated to $z$ and we have shown them both to be different from zero.

Part~(b) of Theorem~\ref{thm:ks} implies that our system has an ergodic partition and that moreover, associated to $\mu$-a.e.\ $z\in R$, there is a local unstable manifold $\gamma^u(z)$ and a local stable manifold $\gamma^s(z)$.

%% file: chapter4/section3.tex
\section{Global arguments}\label{SLTM.Toral2}
\setcounter{equation}{0}

We conclude the proof of the Bernoulli property by giving the `global' aspects of the argument. We have established that $H:R\to R$ has an ergodic partition and it remains to show that there is just one positive measure component (i.e.\ a full-measure component) to this. We begin by describing the orientation of local invariant manifolds and deducing that their lengths diverge on iteration with $H$ (or $H^{-1}$, as is appropriate), then show that this behaviour is enough for us to conclude that $H$ is Bernoulli.

\subsection{Orientation of local invariant manifolds}

We discuss the length of local unstable manifolds and of their iterates under $H$. The situation is analogous to Section~\ref{PLTM.Manifolds} in which we looked at local unstable manifolds for the planar linked-twist map. Here our task will be simplified by properties of the tangent cones $C$ and $\tilde{C}$.

For $\mu$-a.e.\ $z\in R$ the local unstable manifold $\gamma^u(z)$ exists and is of the form (\ref{eqn:st}). Thus there is an open interval $I\subset\mathbb{R}$ and a diffeomorphism $\alpha:I\to\gamma^u(z)\subset R$. At $z$, the local unstable manifold is tangent to $E^u(z)\subset\mathbb{T}^2$. It is clear (although Proposition~\ref{lem:unstable_direction} may be used for a rigorous proof) that $E^u(z)\subset C$.

In order to define the length of iterates of the local unstable manifold we need the following lemma:

\begin{lem}\label{lem:pltm_lum_countable}
For $\mu$-a.e.\ $z\in R$ and for any $n\in\mathbb{N}$, the set of $z'\in\gamma^u(z)$ at which $DH^n$ does not exist is at most countable.
\end{lem}
\begin{proof}
The proof is a consequence of the orientation of $\gamma^u(z)$ and of the constituent line segments of the set of non-differentiable points for $DH^n$. The former, we have argued, lies in $C$. We claim that the latter is given by the set
\beql\label{eqn:non_diff_H}
s_n=\bigcup_{i=0}^nH^{-i}\left(s_0\right).
\eeql
We justified this previously in the case $n=1$ and the general case follows by induction on $n$.

\begin{figure}[htp]
\centering
\includegraphics[totalheight=0.3\textheight]{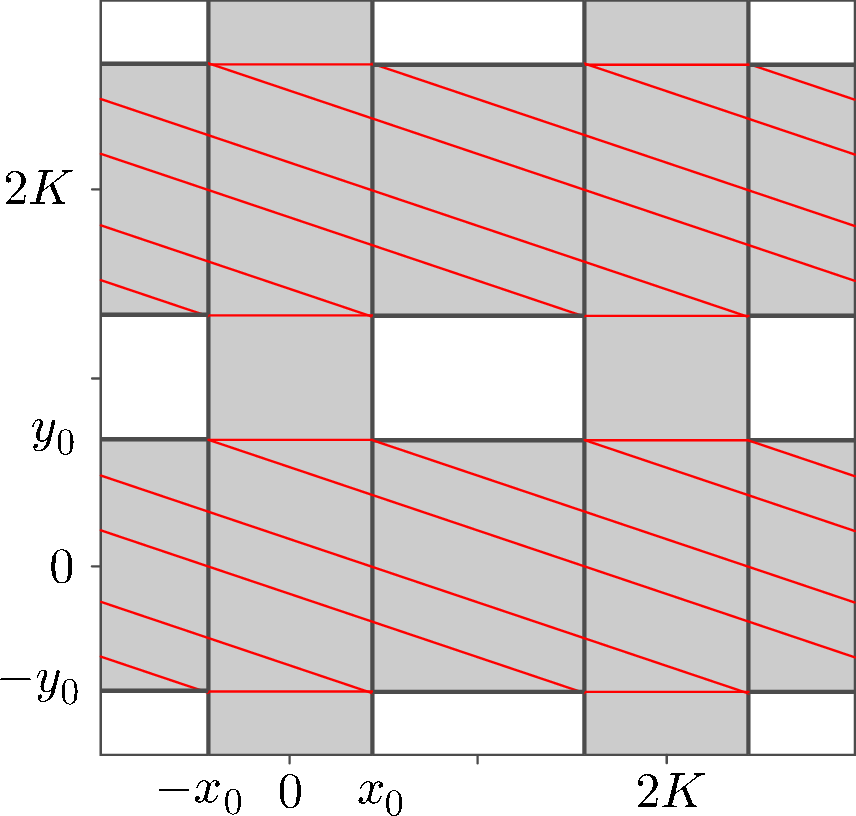}
\caption[Non-differentiable points for $H$]{The manifold $R\subset\mathbb{T}^2$ showing those points (in red) at which $H$ is not differentiable. We illustrate only those points in the interior of the annuli. The map is not differentiable on the boundary either.}
\label{fig:non_diff}
\end{figure}

Figure~\ref{fig:non_diff} illustrates the set $s_0$. We did not construct local \emph{stable} manifolds explicitly, but had we done so we would have formulated a result completely analogous to Lemma~\ref{lem:t_non-shrink}, showing in particular that the cone $\tilde{C}=\{(u,v):uv<0\}$ is preserved by $DH^{-1}$. By referring to the proof of Lemma~\ref{lem:t_non-shrink} it is straight-forward to see that the statement still holds if we replace $\tilde{C}$ with its closure.

We conclude as follows. The non-differentiable set for $H$ consists of finitely many line-segments, each having orientation within the closure of the cone $\tilde{C}$. These orientations are preserved by $DH^{-1}$ so the statement holds true for the non-differentiable set for $H^n$, with $n$ some positive integer. Because $\gamma^u(z)$ has transversal orientation, all intersections
\beq*
\gamma^u(z)\cap s_n
\eeq*
are transversal. It follows from Lemma~\ref{lem:intersection_of_curves} that there are at most countably many such intersections.
\end{proof}

Our proof of the growth of local unstable manifolds is now very similar to that of Theorem~\ref{thm:growth_of_unstable_manifolds}.

\begin{thm}
Let $n\in\mathbb{N}$ and let $\lambda>1$ be the constant given by Lemma~\ref{lem:t_non-shrink}. For $\mu$-a.e.\ $z\in R$ the length of $H_S^n(z)$ is defined and
\beql\label{eqn:planar_growth}
\text{length}\left(H_s^n(\gamma^u(z))\right)\geqs\lambda^n\text{length}\left(\gamma^u(z)\right).
\eeql
\end{thm}
\begin{proof}
For $\mu$-a.e.\ $z\in R$ there is an interval $I\subset\mathbb{R}$ and a diffeomorphism $\alpha:I\to\gamma^u(z)\subset R$. Lemma~\ref{lem:pltm_lum_countable} says that for a.e.\ such $z$ and for any $n\in\mathbb{N}$ the interval $I$ has a countable (at most) decomposition
\beq*
I=(i_0,i_1)\cup\bigcup_{h=1}^N[i_h,i_{h+1})
\eeq*
for some $N\in\mathbb{N}\cup\{\infty\}$, and the restriction of $H_S^n\circ\alpha$ to any interval $(i_h,i_{h+1})$ in this decomposition is a smooth curve. The length is given by our definition above and we have
\beqa*
\text{length}\left(H^n_S(\gamma^u(z))\right) &=\text{length}(H^n_S\circ\alpha) \\
&=\sum_{h=0}^{\infty}\int_{i_h}^{i_{h+1}}\left\|DH^n_SD\alpha_t\right\|\d t \\
&\geqs \sum_{h=0}^{\infty}\int_{i_h}^{i_{h+1}}\lambda^n\left\|D\alpha_t\right\|\d t \\
&=\lambda^n\sum_{h=0}^{\infty}\int_{i_h}^{i_{h+1}}\left\|D\alpha_t\right\|\d t \\
&=\lambda^n\int_I\left\|D\alpha_t\right\|\d t \\
&=\lambda^n\text{length}(\alpha) \\
&=\lambda^n\text{length}\left(\gamma^u(z)\right).
\eeqa*
\end{proof}

We have deduced exponential growth of local invariant manifolds with respect to the first return map $H_S$, and that the orientation of these manifolds is restricted to the cone $C$. The Bernoulli property follows exactly as in Section~\ref{PLTM.Bernoulli.Final}. This completes the proof of Theorem~\ref{thm:Bernoulli_DTLTM}.

%% file: chapter4/section4.tex
\section{Proof of the main result}\label{SLTM.Main}
\setcounter{equation}{0}

In this final section we complete the proof of Theorem~\ref{thm:main_sphere} as follows. The reader should recall the definition of an embedding given in Section~\ref{intro.abs.rev}. We show in Section~\ref{SLTM.Main.Embedding} that $E$ is an embedding of the cylinders $C$ and $C'$ into $\mathbb{S}^2$. In Section~\ref{SLTM.Main.Semi} we demonstrate that $\Theta$ and $H$ are semi-conjugate; this was after all the purpose of introducing the map $H$. We complete the proof in Section~\ref{SLTM.Main.Final} by describing a result due to \citet{orn2} and proving some accompanying results of a technical nature.

Before we begin we must caution the reader as to our notation. In previous sections $(u,v)$ have been reserved for points in a tangent space. From here on $(u,v,w)$ shall represent Cartesian coordinates in $\mathbb{R}^3$. Similarly whereas $C$ has denoted until now a particular tangent cone, from now on it will denote the cylinder in $\mathbb{T}^2$ which is (it will be proven) embedded into $\mathbb{S}^2$.

\subsection{$E$ is an embedding}\label{SLTM.Main.Embedding}

We show that $E:C\to P\subset\mathbb{S}^2$ is an embedding as defined above. By similar arguments one can show that $E\circ N:C'\to Q\subset\mathbb{S}^2$ is also an embedding. We begin by reviewing some properties of the Jacobi elliptic functions which we will require.

\begin{prop}[\citet{meyer}]\label{prop:Jacobi_functions}
Let $\sn,\cn,\dn:\mathbb{R}\to\mathbb{R}$ denote the Jacobi elliptic functions, where the parameter $k$ is fixed at $\sqrt{2}/2$. Let $K$ be as defined by (\ref{eqn:JacobiK}).

The functions $\sn$ and $\cn$ are periodic with period $4K$ and $\dn$ is periodic with period $2K$. For any $t\in\mathbb{R}$ we have the relationships
\beq*
\sn(2K-t)=\sn(t),\quad\cn(2K-t)=-\cn(t)\and\dn(2K-t)=\dn(t);
\eeq*
we have the addition formulae
\beq*
\cn^2(t)+\sn^2(t)=1\and\dn^2(t)+\frac{1}{2}\sn^2(t)=1;
\eeq*
the derivatives are given by
\beq*
\frac{\d}{\d t}\sn(t)=\cn(t)\dn(t),\quad\frac{\d}{\d t}\cn(t)=-\dn(t)\sn(t)\and\frac{\d}{\d t}\dn(t)=-\frac{1}{2}\sn(t)\cn(t);
\eeq*
$\sn$ is a homeomorphism of $(-K,K)$, whereas $\cn$ is a homeomorphism of $(0,2K)$; and finally $\sn$ is positive on $(0,2K)$, negative on $(2K,4K)$, $\cn$ is positive on $(-K,K)$, negative on $(K,3K)$ and $\dn$ is positive for all $t\in\mathbb{R}$.
\end{prop}

Recall that it is convenient to consider $C$ as a subset of $\mathbb{T}^2$. Our first result concerns the image of $C$ with respect to $E$.

\begin{lem}
$E(\mathbb{T}^2)\subset\mathbb{S}^2$.
\end{lem}
\begin{proof}
Let $(x,y)\in\mathbb{T}^2$ and let $(u,v,w)=E(x,y)$. Using Proposition~\ref{prop:Jacobi_functions} we find that
\beqa*
u^2+v^2+w^2 &=\sn^2(x)\dn^2(y)+\cn^2(x)\cn^2(y)+\dn^2(x)\sn^2(y) \\
 &=\sn^2(x)\left(1-\frac{1}{2}\sn^2(y)\right)+\left(1-\sn^2(x)\right)\left(1-\sn^2(y)\right) \\
 &\quad+\left(1-\frac{1}{2}\sn^2(x)\right)\sn^2(y) \\
						 &=1.
\eeqa*
\end{proof}

\begin{lem}
$C\subset\mathbb{T}^2$ and $\mathbb{S}^2\subset\mathbb{R}^3$ are smooth manifolds of dimension 2.
\end{lem}
\begin{proof}
To be precise we should say that $C$ is a smooth manifold \emph{with boundary}, although this introduces extra technical considerations and requires further definitions and we are not in fact interested in the behaviour of the twist map on the boundary, which consists only of fixed points. Thus we show that the \emph{interior} of $C$ is a smooth manifold of dimension 2. We justify this by observing that if we denote the boundary by $\pd C\subset\mathbb{T}^2$ then it is clear that $\mu(C\backslash\pd C)=\mu(C)$.

Let $0<\eps<K$ and define subsets of $\mathbb{R}^2$
\beq*
U_1=\left\{(x,y):x\in(-\eps,3K),y\in(-y_0,y_0)\right\}
\eeq*
and
\beq*
U_2=\left\{(x,y):x\in(K,4K+\eps),y\in(-y_0,y_0)\right\}.
\eeq*
Let $\phi_1:U_1\to C$ and $\phi_2:U_2\to C$ be given by
\beq*
\phi_1(x,y)=\phi_2(x,y)=(x\text{ mod }\mathbb{S}^1,y).
\eeq*
It is not difficult to check that $\bigcup_{i=1}^2\phi_i(U_i)$ covers the interior of $C$ and that the other conditions given in the definition of a smooth manifold of dimension 2 are satisfied.

We now turn to $\mathbb{S}^2$. Let $V\subset\mathbb{R}^2$ be the open unit ball in $\mathbb{R}^2$, i.e.\
\beq*
V=\{(x,y)\in\mathbb{R}^2:x^2+y^2<1\}.
\eeq*
For $t\in\{u,v,w\}$ define $\psi_{\pm}^t:V\to\mathbb{S}^2$ by
\beq*
\psi_{\pm}^u:(x,y)=\left(\pm\sqrt{1-x^2-y^2},x,y\right),
\eeq*
\beq*
\psi_{\pm}^v:(x,y)=\left(x,\pm\sqrt{1-x^2-y^2},y\right),
\eeq*
\beq*
\psi_{\pm}^w:(x,y)=\left(x,y,\pm\sqrt{1-x^2-y^2}\right).
\eeq*
Each pair $\left\{\psi_+^t(V),\psi_-^t(V)\right\}$ covers all but the great circle $t=0$ so the union of all six coordinate neighbourhoods certainly covers $\mathbb{S}^2$. The other properties are easy to establish.
\end{proof}

\begin{lem}
$E$ is a differentiable map.
\end{lem}
\begin{proof}
Let $(x,y)\in C$, choose $i\in\{1,2\}$ so that $(x,y)\in\phi_i(U_i)$ and choose $t\in\{u,v,w\}$ and either $+$ or $-$ so that $E(x,y)\in\psi_{\pm}^t(V)$. Define
\beql\label{eqn:defn_Utpm}
U_{\pm}^t=\phi_i^{-1}\circ E^{-1}\circ\psi_{\pm}^t(V),
\eeql
where $E^{-1}$ denotes the pre-image (as we have yet to discuss the injectivity or otherwise of $E$). We must show that $E(\phi_i(U_{\pm}^t))\subset\psi_{\pm}^t(V)$ and that
\beql\label{eqn:check_diff}
h=\left(\psi_{\pm}^t\right)^{-1}\circ E\circ\phi_i:U_{\pm}^t\to\mathbb{R}^2
\eeql
is differentiable at $\phi_i^{-1}(x,y)$. The former is merely a rearrangement of (\ref{eqn:defn_Utpm}). From the definition of $\phi_i$ and $E$ it is immediate that they are differentiable (in the case of $E$ we need Proposition~\ref{prop:Jacobi_functions}). The differentiability of $h$ thus depends upon the differentiability of $\left(\psi_{\pm}^t\right)^{-1}$. Without loss of generality take $t=u$ and choose $+$. We have $\left(\psi_+^u\right)^{-1}(u,v,w)=(v,w)$, which again is clearly differentiable.
\end{proof}

\begin{lem}\label{lem:E_immersion}
$E$ is an immersion.
\end{lem}
\begin{proof}
Let $z=(x,y)\in C$. We must verify that $DE_z:T_zC\to T_{E(z)}\mathbb{R}^3$ is injective. The Jacobian matrix for $E$ will not be square and so we will not be able to check the injectivity of the differential $DE_z$ by simply computing its determinant (and showing that this is non-zero).

However, the condition is equivalent to asking whether $Dh_{\phi_i^{-1}(x,y)}$ is invertible, where $h$ is as defined in (\ref{eqn:check_diff}), for each admissible choice of $i\in\{1,2\}$, $t\in\{u,v,w\}$ and either $+$ or $-$. The Jacobian for $h$ \emph{will} be square so it will be simple to check its injectivity.

Without loss of generality we take $i=1$, $t=u$ and take $+$. Then
\beqa*
h(\phi_i^{-1})(x,y) &= \left(\psi_{\pm}^u\right)^{-1}\circ E(x,y) \\
										&= \left(\psi_{\pm}^u\right)^{-1}\left(\sn(x)\dn(y),\cn(x)\cn(y),\dn(x)\sn(y)\right) \\
										&= \left(\cn(x)\cn(y),\dn(x)\sn(y)\right).
\eeqa*
The Jacobian is given by
\beq*
Dh=\left( \begin{array}{cc}
-\sn(x)\dn(x)\cn(y) & -\sn(y)\dn(y)\cn(x) \\
-\frac{1}{2}\sn(x)\cn(x)\sn(y) & \cn(y)\dn(y)\dn(x)
\end{array}\right)
\eeq*
which has determinant
\beql\label{eqn:det_h_stolen}
-\sn(x)\dn(y)\left(\dn^2(x)\cn^2(y)+\frac{1}{2}\sn^2(y)\cn^2(x)\right).
\eeql

We require an explicit expression for the domain of $h$:
\beqa*
\phi_1^{-1}\circ E^{-1}\circ\psi_+^u(V) &=\phi_1^{-1}\circ E^{-1}\left(\left\{(u,v,w)\in\mathbb{S}^2:u>0\right\}\right) \\
																					&=\phi_1^{-1}\left(\left\{(x,y)\in C:\sn(x)\dn(y)>0\right\}\right) \\
																					&=\phi_1^{-1}\left(\left\{(x,y)\in C:x\text{ mod }\mathbb{S}^1\in(0,2K)\right\}\right) \\
																					&=(0,2K)\times[-y_0,y_0],
\eeqa*
where in the third line we have used Proposition~\ref{prop:Jacobi_functions}. Now observe that $\sn(x)$ and $\dn(y)$ are strictly positive on the domain considered. The bracketed term in (\ref{eqn:det_h_stolen}) will be positive unless both $\dn^2(x)\cn^2(y)$ and $\sn^2(y)\cn^2(x)$ are zero. In particular we note that $\dn^2(x)\cn^2(y)>0$ for any $x\in\mathbb{R}$ and $y\in(-K,K)\supset[-y_0,y_0]$.
\end{proof}

\begin{lem}
$E$ is a homeomorphism.
\end{lem}
\begin{proof}
As in the proof of Lemma~\ref{lem:E_immersion} it is simpler to work with the composition $h$ defined by (\ref{eqn:check_diff}). From elementary results about compositions of functions we have that $E$ is a homeomorphism if and only if $h$ is a homeomorphism, for each admissible choice of $i,t$ and $\pm$. Without loss of generality we again consider the particular case $i=1$, $t=u$ and $+$. We have
\beq*
h=\left(\psi_+^u\right)^{-1}\circ E\circ\phi_1:\phi_1^{-1}\left(E^{-1}(\psi_+^u(V))\right)\to V.
\eeq*
Recall from the proof of Lemma~\ref{lem:E_immersion} that in this case the domain of $h$ is given by $(0,2K)\times[-y_0,y_0]$. Let
\beq*
s=\cn(x)\cn(y)\and t=\dn(x)\sn(y)
\eeq*
then using Proposition~\ref{prop:Jacobi_functions} we obtain
\beql\label{eqn:bij1}
2t^2=2\left(1-\frac{1}{2}\sn^2(x)\right)\sn^2(y)=\left(1+\cn^2(x)\right)\sn^2(y)
\eeql
which, combined with $s^2=\cn^2(x)\left(1-\sn^2(y)\right)$ gives
\beql\label{eqn:bij2}
2t^2+s^2=\cn^2(x)+\sn^2(y).
\eeql
We will solve these for $\mathcal{S}=\sn^2(y)$ and $\mathcal{C}=\cn^2(x)$. Rearranging (\ref{eqn:bij1}) in terms of $\mathcal{C}$ and substituting into (\ref{eqn:bij2}) we obtain
\beq*
\mathcal{C}^2-\left(2t^2+s^2-1\right)\mathcal{C}-s^2=0
\eeq*
which we solve using the quadratic formula to get
\beql\label{eqn:bij3}
\mathcal{C}=\frac{1}{2}\left(2t^2+s^2-1\pm\sqrt{\left(2t^2+s^2-1\right)^2+4s^2}\right).
\eeql
It is clear that $\sqrt{\left(2t^2+s^2-1\right)^2+4s^2}\geqs 2t^2+s^2-1$ with equality if and only if $s=0$. Clearly $\mathcal{C}=\cn^2(x)$ is non-negative and so in (\ref{eqn:bij3}) we must always take $\pm$ to be $+$. We arrive at
\beql\label{eqn:bij4}
\cn(x)=\pm\sqrt{\frac{1}{2}\left(2t^2+s^2-1+\sqrt{\left(2t^2+s^2-1\right)^2+4s^2}\right)}
\eeql
and again must determine the sign to be taken in (\ref{eqn:bij4}). Notice (see Figure~\ref{fig:elliptic_functions}) that $\cn(y)$ is positive for each $y\in I=[-y_0,y_0]$. Thus $\sgn(\cn(x))=\sgn(\cn(x)\cn(y))=\sgn(s)$ on the whole domain and $\cn(x)$ is well-defined by
\beql\label{eqn:bij5}
\cn(x)=\sgn(s)\sqrt{\frac{1}{2}\left(2t^2+s^2-1+\sqrt{\left(2t^2+s^2-1\right)^2+4s^2}\right)}.
\eeql
We need only observe (Proposition~\ref{prop:Jacobi_functions}) that $\cn$ is bijective on $(0,2K)$ to see that $x$ is uniquely defined by (\ref{eqn:bij5}).

Substituting (\ref{eqn:bij5}) into (\ref{eqn:bij2}) gives
\beql\label{eqn:bij6}
\sn^2(y)=t^2+\frac{1}{2}\left(s^2+1\right)-\frac{1}{2}\sqrt{\left(2t^2+s^2-1\right)^2+4s^2}.
\eeql
Here we need only notice that $\dn(x)$ is strictly positive to see that $\sgn(\sn(y))=\sgn(t)$ giving
\beql\label{eqn:bij7}
\sn(y)=\sgn(t)\sqrt{t^2+\frac{1}{2}\left(s^2+1\right)-\frac{1}{2}\sqrt{\left(2t^2+s^2-1\right)^2+4s^2}}.
\eeql
The fact that $\sn$ is bijective on $[-y_0,y_0]$ completes the proof.
\end{proof}

\subsection{$\Theta$ and $H$ are semi-conjugate}\label{SLTM.Main.Semi}

The cornerstone of our proof that $\Theta:A\to A$ is Bernoulli is the result that it is semi-conjugate to $H$. The formal statement follows and will be proven in this section.
\begin{prop}\label{prop:semi_conj}
The identity $E\circ H=\Theta\circ E$ holds on $R$.
\end{prop}
We first show that $E\circ F=\Phi\circ E$ for $(x,y)\in P$. We deal separately with the cases $(x,y)\in P_0$ and $(x,y)\in P_1$, the former being trivial:

\begin{lem}\label{lem:semi1}
The identity $E\circ F=\Phi\circ E$ holds on $P_0$.
\end{lem}
\begin{proof}
By definition $\Phi:A_+\to A_+$ is given by $E\circ T\circ E^{-1}$, where by $E^{-1}$ we mean the inverse to the restriction of $E$ to $P_0$. $T:C\to C$ is the linear twist map defined in Section~\ref{Intro.Background.Abstract} so it suffices to show that given $(x,y)\in\mathbb{S}^1\times[-y_0,y_0]$ then $T(x,y)=F(x,y)$. This is immediate from the respective definitions.
\end{proof}

Proving that $E\circ F=\Phi\circ E$ on $P_1$ will require a little more work. Let $J:\mathbb{T}^2\to\mathbb{T}^2$ be defined by $J(x,y)=(2K-x,2K-y)$. It is easy to establish that $J$ is a diffeomorphism, is its own inverse and that $J(P_0)=P_1$.
\footnote{$J$ is sometimes referred to as an \emph{involution} and reflects points through the origin in $\mathbb{T}^2$. We are grateful to Prof.\ Robert MacKay for pointing out to us the fact that $\mathbb{T}^2$ with all pairs of points $\{(x,y),J(x,y)\}$ identified is topologically equivalent to $\mathbb{S}^2$. Our map $E$ provides an explicit means by which one might relate a topology on $\mathbb{T}^2$ to the quotient topology on $\mathbb{S}^2$.}

\begin{lem}\label{lem:propJ1}
The identity $E\circ J=E$ holds on $\mathbb{T}^2$.
\end{lem}
\begin{proof}
The proof is an elementary application of Proposition~\ref{prop:Jacobi_functions}. $E(J(x,y))$ is given by
\beqa*
& \left(\sn(2K-x)\dn(2K-y),\cn(2K-x)\cn(2K-y),\dn(2K-x)\sn(2K-y)\right) \\
=& \left(\sn(x)\dn(y),-\cn(x)\cdot-\cn(y),\dn(x)\sn(y)\right) \\
=& E(x,y).
\eeqa*
\end{proof}

\begin{lem}\label{lem:propJ2}
The identity $J\circ F=F\circ J$ holds on $P_0$.
\end{lem}
\begin{proof}
The proof is straight-forward calculation. Let $(x,y)\in P_0$, so $J(x,y)\in P_1$, then:
\beqa*
F\circ J(x,y) &= F(2K-x,2K-y) \\
							&= \left(2K-x+(4Ky_0-4Ky)/2y_0,2K-y\right) \\
							&= \left(2K-x+2K-2Ky/y_0,2K-y\right)
\eeqa*
and conversely
\beqa*
J\circ F(x,y) &= J\left(x+4K(y+y_0)/2y_0,y\right) \\
							&= \left(2K-x-2Ky/y_0-2K,2K-y\right).
\eeqa*
The two are equal because $2K=-2K$ in $\mathbb{S}^1$.
\end{proof}

It is now simple to prove the following:
\begin{lem}\label{lem:semi2}
The identity $E\circ F=\Phi\circ E$ holds on $P_1$.
\end{lem}
\begin{proof}
Let $(\tilde{x},\tilde{y})\in P_1$, let $(x,y)=J(\tilde{x},\tilde{y})$ and observe that $(x,y)\in P_0$. We must show that $E\circ F(\tilde{x},\tilde{y})=\Phi\circ E(\tilde{x},\tilde{y})$. Using Lemma~\ref{lem:propJ1} we have
\beql\label{eqn:conj1}
\Phi\circ E(\tilde{x},\tilde{y})=\Phi\circ E(x,y),
\eeql
and using Lemma~\ref{lem:propJ2} followed by Lemma~\ref{lem:propJ1} we have
\beql\label{eqn:conj2}
E\circ F(\tilde{x},\tilde{y})=E\circ J\circ F(x,y)=E\circ F(x,y).
\eeql
Lemma~\ref{lem:semi1} says that the expressions~(\ref{eqn:conj1}) and~(\ref{eqn:conj2}) are equal.
\end{proof}

In very similar fashion it can be shown that $E\circ G=\Gamma\circ E$ holds on $Q$. The semi-conjugacy follows immediately:
\begin{proof}[Proof of Proposition~\ref{prop:semi_conj}]
Let $(x,y)\in R$, then $E\circ G\circ F(x,y)= \Gamma\circ E\circ F(x,y)=\Gamma\circ\Phi\circ E(x,y)$.
\end{proof}

\subsection{Proof of Theorem~\ref{thm:main_sphere}}\label{SLTM.Main.Final}

Our main result is a consequence of a theorem of \citet{orn2} and essentially follows from the semi-conjugacy just established, although we will need some supplementary results. We need to discuss measure-theoretic factors and we begin by reviewing a few definitions; we have taken these from \citet{kh}.

A measure space $(X,\mathcal{M},\nu)$ with finite measure $\nu$ is called a \emph{Lebesgue space} if it is isomorphic to the union of $[0,1]$ with Lebesgue measure, with at most countably many points of positive measure.

For measure-preserving transformations $T_1:X_1\to X_1$ and $T_2:X_2\to X_2$ of Lebesgue spaces $(X_1,\nu_1)$ and $(X_2,\nu_2)$ respectively, we say that $T_2$ is a \emph{metric factor} of $T_1$, or from now on just a \emph{factor} of $T_1$, if there exists a measure preserving map $\theta:X_1\to X_2$ such that
\beq*
\theta_*\nu_1=\nu_2\and T_2\circ\theta=\theta\circ T_1.
\eeq*
The notation $\theta_*\nu_1$ denotes the \emph{pushforward} measure on $X_2$ obtained from the measure on $X_1$: for each measurable set $B\subset X_2$ this is defined by
\beq*
\theta_*\nu_1(B)=\nu_1(\theta^{-1}(B)).
\eeq*
We will use the following:

\begin{thm}[\citet{orn2}]
\label{thm:factor}
A factor of a Bernoulli map is itself Bernoulli.
\end{thm}

We prove some results concerning $\Theta$ and $(A,\tilde{\mu})$, where $\tilde{\mu}=E_*\mu$. Let $\mathcal{M}$ denote the $\sigma$-algebra of Lebesgue-measurable subsets of $R$ and let
\beq*
\tilde{\mathcal{M}}=\{\tilde{B}\subset A:E^{-1}(\tilde{B})\in\mathcal{M}\},
\eeq*
where $E^{-1}(\tilde{B})$ denotes the pre-image of $\tilde{B}\subset A$ with respect to $E$. It is known (see \citet{rudin3}, Theorem~1.12, p.13) that $\tilde{\mathcal{M}}$ is a $\sigma$-algebra of subsets of $A$.

Recall that $\mu$ denotes the Lebesgue measure on $R$. Let $\tilde{\mu}:\tilde{\mathcal{M}}\to[0,1]$ be given by $\tilde{\mu}=E_*\mu$.

\begin{prop}
\label{prop:measure}
$(A,\tilde{\mu})$ is a Lebesgue space.
\end{prop}
\begin{proof}
The function $\tilde{\mu}$, defined on $\sigma$-algebra $\tilde{\mathcal{M}}$, takes its range in $[0,1]$. We show that it is countably additive and thus a measure. Let $\{\tilde{B}_i:i\in\mathbb{N}\}$ be a disjoint, countable collection
of members of $\tilde{\mathcal{M}}$, so $\{E^{-1}(\tilde{B}_i):i\in\mathbb{N}\}$ is a disjoint, countable collection of members of $\mathcal{M}$. By the countable additivity of $\mu$ we have $\mu\left(\bigcup_{i=1}^{\infty}E^{-1}(\tilde{B}_i)\right)=\sum_{i=1}^{\infty}\mu(E^{-1}(\tilde{B}_i))$. The result then follows from the observation:
\beqa*
\bigcup_{i=1}^{\infty}E^{-1}(\tilde{B}_i)
&=\bigcup_{i=1}^{\infty}\left\{ (u,v)\in R:E(u,v)\in\tilde{B}_i\right\} \\
&=\left\{ (u,v)\in R:E(u,v)\in\bigcup_{i=1}^{\infty}\tilde{B}_i\right\} \\
&=E^{-1}\left(\bigcup_{i=1}^{\infty}\tilde{B}_i\right).
\eeqa*

Next, we appeal to the result that \emph{any Borel probability measure on a separable, locally compact Hausdorff space defines a Lebesgue space} (see \citet{kh}, Theorem~A.6.7, p.734). It is not difficult to see that $(A,\tilde{\mu})$ satisfies these criteria: $A$ is a separable Hausdorff space because $\mathbb{R}^3$ has these properties and by the Heine-Borel theorem it is compact and hence locally compact. It is obvious that $\tilde{\mu}(A)=1$; to see that $\tilde{\mathcal{M}}$ contains all Borel subsets of $A$, let $\tilde{B}\subset A$ be open, then $E^{-1}(\tilde{B})\subset R$ is open (because $E$ is continuous) and so $E^{-1}(\tilde{B})\in\mathcal{M}$. By definition $\tilde{B}\in\tilde{\mathcal{M}}$.
\end{proof}

The result we quoted from \citet{kh} shows equally that $(R,\mu)$ is a Lebesgue space.

\begin{lem}\label{lem:preserve_measure}
$\Theta:A\to A$ preserves $\tilde{\mu}$, i.e.\ if $\tilde{B}\in\tilde{\mathcal{M}}$ then
\beql\label{eqn:preserve_measure}
\tilde{\mu}(\Theta(\tilde{B}))=\tilde{\mu}(\tilde{B}).
\eeql
\end{lem}
\begin{proof}
The proof involves a little manipulation of identities we have established. Let $\tilde{B}\in\tilde{\mathcal{M}}$ and let $B=E^{-1}(\tilde{B})\in\mathcal{M}$. It follows that $E(B)=\tilde{B}$ and
\beql\label{eqn:pm1}
B=E^{-1}\circ E(B).
\eeql
From Lemma~\ref{prop:semi_conj} it follows that $H^{-1}\circ E^{-1}=E^{-1}\circ\Theta^{-1}$ (note that these are pre-images rather than functions). Using this fact and (\ref{eqn:pm1}) we deduce that
\beq*
B=E^{-1}\circ E(B)=E^{-1}\circ\Theta^{-1}\circ\Theta\circ E(B)=H^{-1}\circ E^{-1}\circ E\circ H(B)
\eeq*
i.e.\
\beql\label{eqn:pm2}
H(B)= E^{-1}\circ E\circ H(B).
\eeql
By definition we have
\beql\label{eqn:pm3}
\tilde{\mu}(\tilde{B})=\mu\left(E^{-1}(B)\right)=\mu(B)
\eeql
and using (\ref{eqn:pm1}) followed by (\ref{eqn:pm2}) we have
\beql\label{eqn:pm4}
\tilde{\mu}\left(\Theta(\tilde{B})\right)=\mu\circ E^{-1}\circ\Theta\circ E(B)=\mu\circ E^{-1}\circ E\circ H(B)=\mu\left(H(B)\right).
\eeql
The expressions (\ref{eqn:pm3}) and (\ref{eqn:pm4}) are equal because $H$ preserves $\mu$.
\end{proof}

Our main result follows easily.
\begin{proof}[Proof of theorem~\ref{thm:main_sphere}]
Propositions~\ref{prop:semi_conj} and~\ref{prop:measure} and Lemma~\ref{lem:preserve_measure} show that $\Theta$ is a metric factor of $H$. Theorem~\ref{thm:factor} completes the proof.
\end{proof}

We end with a remark regarding Ornstein's \citeyear{orn2} paper in which Theorem~\ref{thm:factor} is established. Ornstein defines a factor of a Bernoulli shift to be the restriction thereof to an invariant sub-$\sigma$-algebra; we demonstrate briefly that the (more common) definition we have taken is equivalent. Indeed, if we let $\mathcal{M}'\subset\mathcal{M}$ consist of precisely those elements $B\in\mathcal{M}$ for which $B=E^{-1}\circ E(B)$ and let $\mu'$ be the restriction of $\mu$ to $\mathcal{M}'$ then it follows from the results of this section that $E:(R,\mu')\to(A,\tilde{\mu})$ is an isomorphism and thus the result.

%% file: input/concl.tex
%\renewcommand{\thechapter}{\Alph{chapter}}
%\setcounter{chapter}{0}
%\chapter*{Conclusions} 
%\addcontentsline{toc}{chapter}{Conclusions} 
%\chaptermark{Conclusion}
%\setcounter{equation}{0}
%\setcounter{figure}{0}
%\label{concl}

\chapter{Summary and outlook}
\setcounter{equation}{0}
\setcounter{figure}{0}
\label{concl}

We finish by surveying the results we have established and discussing some strengths and weaknesses of our methods. We consider the directions in which productive future work might be undertaken, either as a direct consequence of the present work or otherwise. At the end of the chapter we will propose two conjectures which, we believe, would be an excellent starting point for anyone who wished to generalise our methods to the class of abstract linked-twist maps we have defined.

\section{Summary}

We make some comments about the results we have obtained. 

\subsection{A topological horseshoe in the toral linked-twist map}

In Chapter~\ref{chapter2} we established the existence of a topological horseshoe in the toral linked-twist map defined in Section~\ref{Intro.Theorems.Torus}. Our method was inspired by the work of \citet{d2} who constructed such a horseshoe in the planar linked-twist map defined in Section~\ref{Intro.Theorems.Plane}.

We observed an interesting difference between the two constructions. The conjugacy constructed by Devaney is with a \emph{sub-shift of finite type}, conversely ours is with \emph{full shift on $N$ symbols}.

We believe that this can be explained by the fact that the planar linked-twist map has two distinct intersection regions, so that in this case the invariant Cantor set $\Lambda$ is split between the two. We conjecture, based on our result, that there is an invariant Cantor set $\Lambda_+\subset\Sigma_+$ for the planar map on which the dynamics, as in our example, are conjugate to a full shift on the space of symbol sequences. This set may be constructed by considering only those points that land in $\Sigma_+$ on each iteration, and excluding the other points from Devaney's construction, so in fact it is a proper subset of his invariant set.

Devaney's invariant set has a richer structure than the one that we have constructed, and we believe this to be a consequence of it containing as proper subsets a wealth of other invariant sets whereby there is some restriction on which of $\Sigma_{\pm}$ points return at any given time. An investigation of this structure would certainly be an interesting exercise in its own right, although it is not clear that it would yield any conclusions that might help us to better understand the dynamics on a set of full measure.

\subsection{The Bernoulli planar linked-twist}

The shortcomings of our result are clear: we have established the Bernoulli property for a planar linked-twist map composed of the embeddings of linear twists, but we have had to be explicit about the sizes of the annuli, taking the inner annuli to have size $r_0=2$ and the outer annuli to have size $r_1=\sqrt{7}$. It is clear where these restrictions were required so let us look at this a little more closely.

Key to our proof was Proposition~\ref{prop:DGDF} which states that, with annuli of the sizes specified, then $DF$ preserves the tangent cone $C=\{(u,v):uv>0\}$ and $DF^{-1}$ preserves the tangent cone $\tilde{C}=\{(u,v):uv<0\}$. We do not intend to repeat all of the details here, but an example illustrates the difficulty. We were required to show (see equation~(\ref{eqn:Df_bounds})) that $D_1f_+(x,y)>0$ , or to give the full expression in terms of the function $\psi$, that
\beqa*
D_1\psi\left(x,\psi^{-1}(x,y)+c(x-r_0)\right) \\ +D_2\psi\left(x,\psi^{-1}(x,y)+c(x-r_0)\right)\left[D_1\psi^{-1}(x,y)+\frac{2\pi}{r_1-r_0}\right]>0,
\eeqa*
for each pair $(x,y)\in [r_0,r_1]\times[0,\pi]$.

We took a rather crude approach to this and sought to bound each of the terms $D_1\psi(\cdot,\cdot)$, $D_2\psi(\cdot,\cdot)$ and $D_1\psi^{-1}(\cdot,\cdot)$ individually. We observed a lower bound of $0$ for $D_1\psi$ which (from the proof of Lemma~\ref{lem:D1psi}) seems optimal, but it is quite possible that none of the other bounds established are optimal. Of course, even if we were to obtain optimal bounds on each of the three derivatives individually this would not necessarily give us optimal bounds for $D_1f_+(x,y)$.

Given our crude approach to this problem it is perhaps remarkable that it works for \emph{any} system at all. The fact that we were able to find (and, we should add, with relative ease) choices of $r_0$ and $r_1$ for which the problem is tractable could be interpreted as evidence that the inequalities in fact hold for a much wider choice of annulus size.

If one is motivated to use our method to prove the Bernoulli property for some range of $r_0$ and $r_1$ values then a more sophisticated approach to these inequalities is imperative. Plotting $D_1f_+(x,y)$ over the required domain would be a good start, although even this is non-trivial as one needs a package with sufficiently good programming and graphical capabilities, due to the nature of the functions $\psi$ and $\psi^{-1}$. Such a plot might suggest a way to partition the domain so that tighter bounds can be established on each partition element; this would seem to entail a great deal of work however.

We conclude by discussing how far one might hope to develop the method we have introduced. The ultimate ambition would be to give a `complete description' of the possible dynamics and this, perhaps, might consist of a large open sets of parameter values for $r_0,r_1$ where the Bernoulli property is established, and a complementary set on which it is shown not to occur.

\citet{sturman} provide a number of plots showing numerical simulations of (co-twisting) planar linked-twist maps, some of which appear to exhibit good mixing properties and others which do not. Recall that one of our initial assumptions was of \emph{transversality}, which we believe to be related to our ability to construct the new coordinates (we will say more on this shortly). It would seem from the simulations that transversality is \emph{not} a prerequisite for good mixing, and hence the method we have proposed cannot be expected to provide such a complete description as we have asked for. Of course, this leaves open the possibility that transversality is \emph{sufficient} for a linked-twist map to be Bernoulli.

\subsection{The Bernoulli linked-twist map on the sphere}

It is certainly interesting that we have been able to construct so directly a semi-conjugacy between a map on the torus and a map on the sphere. This is an immediate consequence of the coordinate system we have used. The coordinate transformation would perhaps be of interest to the wider mathematical community given its relatively clean expression and the orthogonal coordinate system it provides. The most interesting development from a dynamical systems perspective might be to use it to construct further examples of Bernoulli maps on full measure subsets of the sphere, by a method analogous to Katok's (\citeyear{katok}). Recall that the starting point for his construction is a hyperbolic toral automorphism with certain points fixed.

Perhaps the strength of the method is also its weakness; it is quite specialised and so it is difficult to see how one might hope to generalise it in order to obtain other results, or indeed what those other results might be. Nevertheless it afforded us the opportunity to use some techniques (the theorem of \citet{orn2}) that perhaps otherwise we would not have discovered.

\section{Ideas for further work}

We conclude by looking at two ways in which one might build upon the results we have established.

\subsection{Decay of correlations}

We have mentioned many applications for which certain linked-twist maps provide a natural model and thereby a means to understand or to predict the behaviour to some extent. This is of particular importance within the nanoscale devices we have mentioned such as the DNA microarray discussed in Section~\ref{LitReview.App.Torus}, because the alternative \emph{trial and error} approach to their design is prohibitively costly. We have mentioned that there is a degree of convergence between those questions that are interesting from a mathematical perspective and those that are interesting from an applications perspective.

The \emph{strength} of mixing is one such question. We have established for two different types of linked-twist map that the Bernoulli property is satisfied. The implication for systems whose dynamics are well approximated (in some sense) by these maps is that they should be expected to mix initial conditions thoroughly.

The \emph{rate} of mixing is another such question which we shall discuss briefly now. Consider for a moment the cornerstone of our proofs; we have spelled this out previously but we re-iterate it now. There is a region of positive area (which we have labeled $\Sigma$ in each case) to which almost every point returns an infinite number of times. The hyperbolicity, which is responsible for the separation of nearby trajectories and thus the strong mixing, is inextricably linked to this behaviour.

Now consider the size (i.e.\ the measure) of this region relative to the size of the whole manifold $A$. Bernoulli systems automatically satisfy the strong mixing property (defined in Section~\ref{LitReview.ErgTh.Erg}) which says that for a measurable set $B\subset A$ having positive measure
\beq*
\lim_{n\to\infty}\frac{\mu\left(\Theta^n(B)\cap\Sigma\right)}{\mu(B)}=\mu(\Sigma),
\eeq*
where we have used the invertability of $\Theta$. We might interpret this as saying that the asymptotic proportion of $B$ in $\Sigma$ is proportional to the size of $\Sigma$. Given the relationship between returns to $\Sigma$ and separation of nearby trajectories, one might conjecture that the greater the size of $\Sigma$, the more `chaotic' the system is in some sense and, importantly for applications, the faster the phase space becomes mixed.

The concept of the \emph{decay of correlations} is the correct mathematical formalism within which to phrase such questions. In essence the idea is to look at the rate of convergence in the ergodic theorem. We don't provide a formal definition but direct the reader to \citet{baladi} for further details.

The seminal work in recent years on the decay of correlations in dynamical systems with some hyperbolicity is \citet{young}. She establishes a framework for studying this decay in a class of systems she characterises as having \emph{`regular returns to sets with good hyperbolic properties'}. In this context we see the importance of the perspective we have taken in analysing linked-twist maps and why we are hopeful that this approach will prove useful in future endeavours.

\subsection{Ergodic properties of abstract linked-twist maps}

Our definion of an abstract linked-twist map invites the question of its ergodic properties. We discuss briefly how the proofs of general results along these lines might be attempted using the techniques introduced in this work. We stress that this should not be considered a `work in progress'; rather these are merely preliminary comments which we hope may be of inspiration to anyone inclined to persue results along these lines, and may at least be of some interest to other readers. It is worth remarking that \citet{p_preprint} appears to have attempted results along these lines. To what extent he has acheived these ambitions is not entirely clear to us.

Recall our definition of an abstract linked-twist map on a smooth manifold $M$ of dimension 2: for $i=1,2$ let $C_i$ be a cylinder, let $E_i:C_i\to A\subset M$ be an embedding, let $T_i:C_i\to C_i$ be a twist map on $C_i$ and let $E_i\circ T_i\circ E_i^{-1}$ be a twist map on $E_i(C_i)\subset R$. If $E_1(C_1)$ and $E_2(C_2)$ are transversal then the composition
\beql\label{eqn:abs_ltm}
\Theta=E_2\circ T_2\circ E_2^{-1}\circ E_1\circ T_1\circ E_1^{-1}
\eeql
is called a linked-twist map.

In analysing the linked-twist map in the plane it was crucial that intersection regions $\Sigma_{\pm}$ (that is, the connected components of $E_1(C_1)\cap E_2(C_2)$) could each be expressed in new coordinates in which it was the Cartesian product of two intervals (a `square', in fact). (In the case of the sphere this followed immediately from our choice of coordinates.) We conjecture the following:
\begin{conj}
Transversality of the embedded cones $E_1(C_1)$ and $E_2(C_2)$ is a sufficient condition for the existence of local coordinates in which each connected component of $E_1(C_1)\cap E_2(C_2)$ is a Cartesian product of two intervals.
\end{conj}
\begin{figure}[htp]
\centering
(a)\includegraphics[totalheight=0.25\textheight]{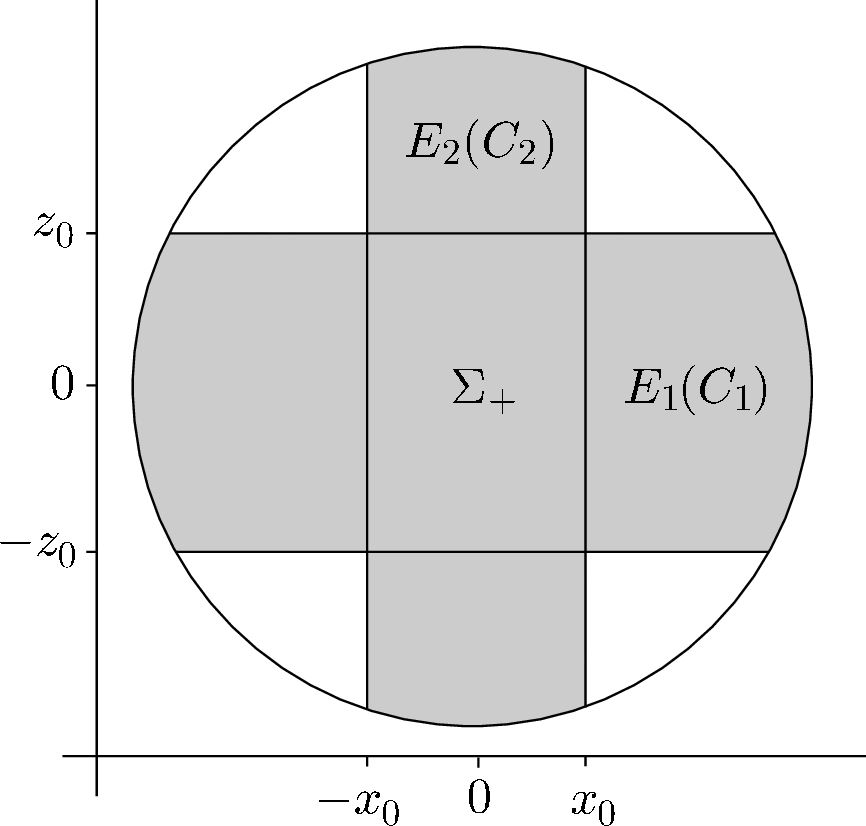}\quad
(b)\includegraphics[totalheight=0.25\textheight]{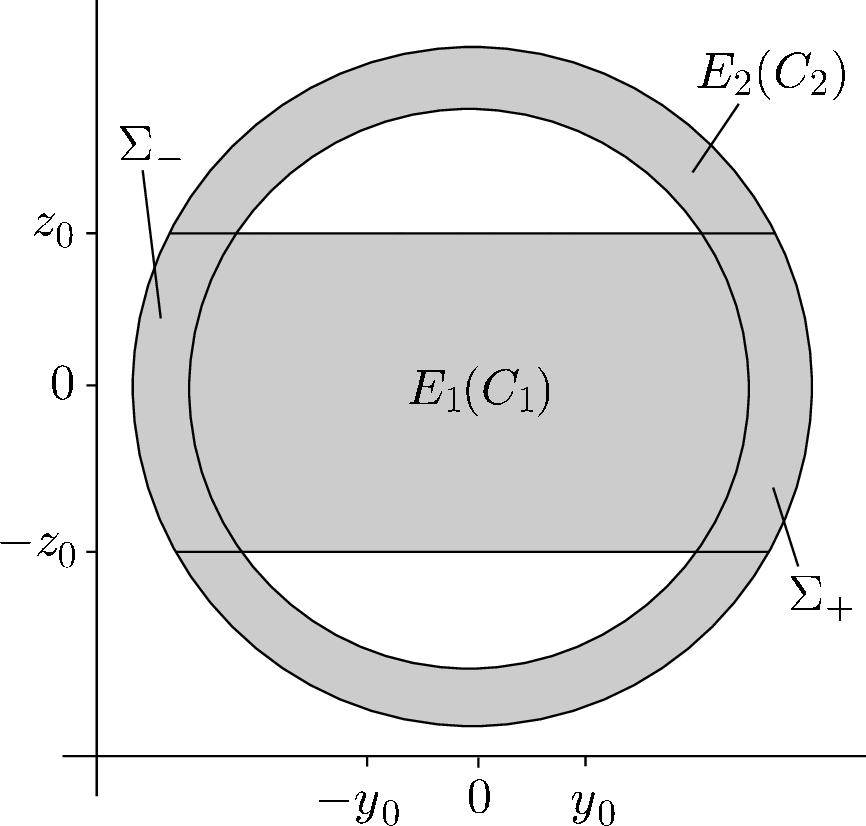}
\caption[An abstract linked-twist map]{An example of an abstract linked-twist map. Part~(a) shows the view parallel to the $y$-axis and part~(b) parallel to the $x$-axis.  The manifold $A$ (shaded) consists of two cylinders $C_1$ and $C_2$ embedded into $\mathbb{S}^2$ with embeddings $E_1$ and $E_2$ respectively.}
\label{fig:abstract_ltm_example}
\end{figure}

We illustrate an example of the situation we have in mind in Figure~\ref{fig:abstract_ltm_example}.
\footnote{We thank Prof.\ Jens Marklof for suggesting this map to us. It is a linked-twist map defined on $\mathbb{S}^2$ and, as opposed to the map studied in Chapter~\ref{chapter4}, points moving under a twist map $E_i\circ T_i\circ E_i$ do so in a plane of constant $x$ or $z$ coordinate (where $(x,y,z)$ are Cartesians in $\mathbb{R}^3$). Moreso than the map of Chapter~\ref{chapter4} this resembles the motion of a `top' undergoing precession. It is therefore possible that this map might embody the essence of certain quantum chaotic motion.}
Part~(a) shows a projection onto the $xz$-plane, whereas part~(b) shows a projection onto the $yz$-plane. The two embedded cylinders $E_1(C_1)$ and $E_2(C_2)$ are bounded by lines of constant $z$ and $x$ coordinate respectively and their union is denoted by $A\subset\mathbb{S}^2$.

\begin{figure}[htp]
\centering
\includegraphics[totalheight=0.3\textheight]{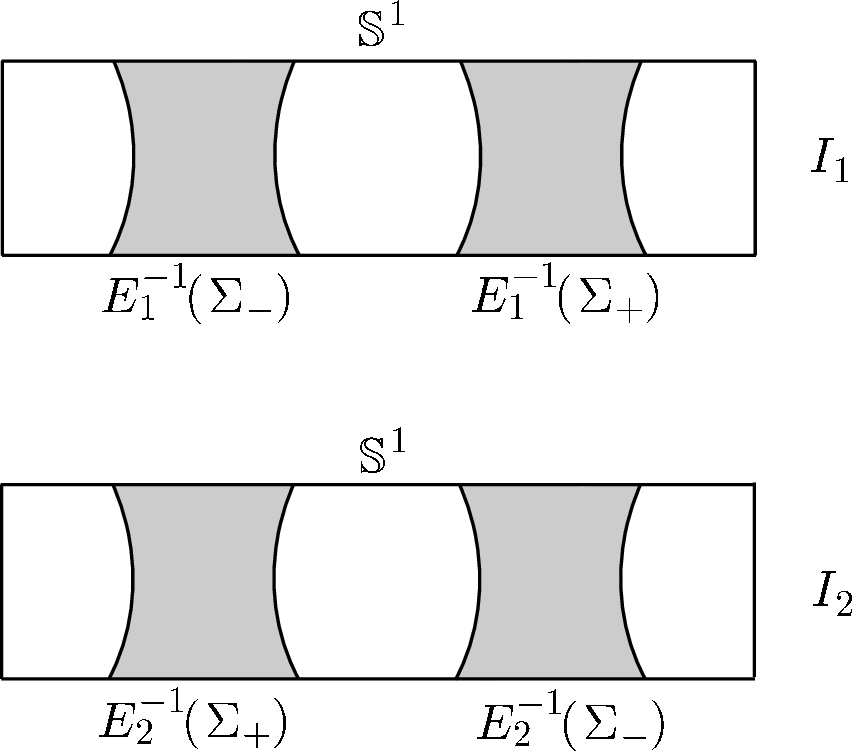}
\caption[Cylinders for an abstract linked-twist map]{The two cylinders with pre-images of the intersection regions shaded.}
\label{fig:abstract_ltm_cyl}
\end{figure}

In Figure~\ref{fig:abstract_ltm_cyl} we show the cylinders themselves with the pre-images of the intersection regions shaded. Let $(s_1,i_1)\in\mathbb{S}^1\times I_1$ give coordinates on $C_1$ and let $(s_2,i_2)\in\mathbb{S}^1\times I_2$ give coordinates on $C_2$.

Let $\mu$ denote Lebesgue measure on $A$; it follows by the usual arguments (assuming that the $T_i$ are sufficiently well behaved) that $\mu$-a.e.\ point returns infinitely many times to $\Sigma_{\pm}$. Define the first return map $\Theta_{\Sigma}:E_1^{-1}(\Sigma)\to E_1^{-1}(\Sigma)$ by
\beq*
\Theta_{\Sigma}=E_1^{-1}\circ E_2\circ T_2^m\circ E_2^{-1}\circ E_1\circ T_1^n
\eeq*
where $n$ and $m$ are positive integers satisfying the usual criteria. Finally, define the usual tangent cone $C$ consisting of the open first and third quadrants. If $D\left(E_2^{-1}\circ E_1\right)$ and $D\left(E_1^{-1}\circ E_2\right)$ preserve and expand the cone $C$ then $D\Theta_{\Sigma}$ will also. In this case we conjecture that $\Theta$ has the Bernoulli property.

If this isn't the case (a situation analogous to the planar linked-twist, where a larger cone $U$ was preserved but $C$ was not) then we might still be able to proceed as before. Using the ideas of Chapter~\ref{chapter3} we can construct new coordinates (if the previous conjecture holds) on $A$, which are equal to $(i_1,i_2)$ for a point $(u,v)\in\Sigma_{\pm}$ such that $E_j^{-1}(u,v)=(i_j,s_j)$, $j=1,2$. The coordinate transformations we have mentioned can be expressed in terms of the embeddings $E_1,E_2$.

\begin{conj}
We can establish sufficient criteria for $\Theta$ to have the Bernoulli property. These criteria consist of a pair of inequalities involving only the derivatives of $E_i^{\pm 1}$ and $T_i$, $i=1,2$.
\end{conj}